\documentclass[thmsa,final,11pt]{article}

\usepackage{eurosym}
\usepackage{amsfonts}
\usepackage{amsmath}
\usepackage{amssymb}
\usepackage{mathrsfs} 
\usepackage{enumerate}
\usepackage{graphicx}
\usepackage{subfigure}
\usepackage{color,xcolor}
\usepackage[onehalfspacing,nodisplayskipstretch]{setspace}

\usepackage{xr}

\newtheorem{theorem}{Theorem}

\newtheorem{lemma}[theorem]{Lemma}

\newtheorem{proposition}[theorem]{Proposition}

\normalsize
\newenvironment{prooft}[1]{\vskip 2mm\noindent {\bf Proof of #1.}}
                    {\hfill $\square$ \vskip 2mm \noindent}

\begin{document}

\title{Asymptotic properties of maximum composite likelihood estimators for max-stable Brown-Resnick random fields over a fixed-domain}
\author{Nicolas CHENAVIER\thanks{%
Universit\'{e} du Littoral C\^{o}te d'Opale, 50 rue F. Buisson 62228 Calais.
nicolas.chenavier@univ-littoral.fr}\quad and\quad Christian Y.\ ROBERT%
\thanks{%
1. Universit\'{e} de
Lyon, Universit\'{e} Lyon 1, Institut de Science Financi\`{e}re et
d'Assurances, 50 Avenue Tony Garnier, F-69007 Lyon, France.
2. Laboratory in Finance and Insurance - LFA CREST - Center for Research in
Economics and Statistics, ENSAE, Palaiseau, France;  christian.robert@univ-lyon1.fr%
}}
\maketitle

\abstract{

Likelihood-based inference for max-stable random fields is challenging, since
finite-dimensional densities are either unavailable in closed form or
computationally intractable in moderate to high dimension. Composite likelihood
methods, based on low-dimensional marginal densities, therefore provide a natural
alternative. In this paper, we study maximum composite likelihood estimation for
spatial Brown--Resnick random fields generated by isotropic fractional Brownian
fields. We work under fixed-domain asymptotics: a single realization of the
max-stable field is observed on an increasingly dense random set of sites, given
by a homogeneous Poisson point process. Pairwise and
triplewise composite likelihoods are constructed by retaining, respectively, the
edges and the triangles of the associated Poisson--Delaunay triangulation.
Our main results establish the consistency of the resulting maximum composite
likelihood estimators of the scale and smoothness parameters, when the other
parameter is known. Their asymptotic behaviour is non-standard: the estimators
converge at rates depending on the smoothness parameter and their centered limits
are non-Gaussian. More precisely, the limiting fluctuations are driven by
aggregated local times associated with the canonical tessellation of the
Brown--Resnick field. These results reveal a fundamental departure from the
classical composite likelihood theory based on increasing domains or independent
replications, and show that Gaussian uncertainty quantification may be
misleading in fixed-domain inference for max-stable spatial extremes.}


\strut

\textit{AMS (2020):} 62G32, 62M30, 60F05, 62H11.

\maketitle

\section{Introduction}
\label{sec:introduction}

Gaussian random fields are widely used in spatial statistics because their
finite-dimensional distributions are completely determined by their mean and
covariance functions. In parametric models, inference therefore reduces to the
estimation of covariance or variogram parameters. When the focus is on extreme
events rather than on average spatial variability, max-stable random fields
provide a natural alternative. They arise as the only possible non-degenerate
limits of normalized pointwise maxima of independent random fields and are now
standard models for spatial extremes.

Likelihood-based inference for max-stable random fields is, however, difficult.
Except in low dimension, the corresponding finite-dimensional densities are
typically unavailable in closed form or computationally intractable. Composite
likelihood methods, based on low-dimensional marginal densities, are therefore
widely used in practice. Most existing asymptotic results for such estimators
are obtained either under increasing-domain asymptotics or under repeated
sampling, where independent temporal replications of the spatial field are
available.

In this paper, we study maximum composite likelihood estimators under a
fixed-domain, or infill, asymptotic framework for Brown--Resnick random fields.
The observation domain is fixed, but the field is observed on an increasingly
dense random set of locations generated by a homogeneous Poisson point process.
Only a single spatial realization of the max-stable field is available. This
setting is substantially different from the replicated framework usually used in
spatial extremes.

Our main result shows that pairwise and triplewise maximum composite likelihood
estimators remain consistent in this fixed-domain setting, when one parameter is
estimated while the other is kept fixed. However, their asymptotic behaviour is
non-standard. The estimators converge at rates depending on the smoothness
parameter, and their centered limits are non-Gaussian. These limits are driven
by aggregated local times associated with the canonical tessellation of the
underlying Brown--Resnick field. This reveals a fundamental departure from the
usual composite likelihood theory and shows that Gaussian approximations may be
misleading for uncertainty quantification from a single spatial realization.

As background, we first review the existing literature on maximum likelihood
and composite likelihood estimation for Gaussian random fields under
fixed-domain asymptotics. We then discuss related results for non-Gaussian
random fields and recall the spectral representation and canonical tessellation
of max-stable random fields.

\subsection{Maximum likelihood estimators for Gaussian random fields under
fixed-domain asymptotics}

The fixed-domain asymptotic framework, also called infill asymptotics
\cite{Stein99,Cressie93}, corresponds to the situation where increasingly many
observations are collected in a fixed bounded sampling domain, usually a subset
of $\mathbb{R}^{d}$, $d\geq 1$. Within this framework, the maximum likelihood
estimators of covariance parameters of Gaussian random fields have been studied
extensively over the last three decades.

A central distinction is between microergodic and non-microergodic parameters.
A parameter is said to be microergodic if, for two different values of this
parameter, the corresponding Gaussian measures are orthogonal
\cite{Ibragimov&Rozanov78,Stein99}. It is non-microergodic if different values
of the parameter lead to equivalent Gaussian measures. Non-microergodic
parameters cannot be estimated consistently under fixed-domain asymptotics.
There is no general theory describing the asymptotic behaviour of maximum
likelihood estimators for all microergodic parameters, and most available
results are model-specific. Examples include the exponential covariance model
\cite{Ying91,Ying93,Vaart96}, the Mat\'ern covariance model
\cite{Matern60,Zhang04,Anderes10,Kaufman&Shaby13}, and the generalized Cauchy
covariance model \cite{Bevilacqua&Faouzi19}.

\subsection{Maximum composite likelihood estimators for Gaussian random fields
under fixed-domain asymptotics}

From a theoretical viewpoint, maximum likelihood is the natural benchmark for
estimating the covariance parameters of a Gaussian random field. Nevertheless,
evaluating the Gaussian likelihood for $n$ observations requires the inversion
of an $n\times n$ covariance matrix and has computational complexity of order
$O(n^{3})$. This becomes prohibitive for large spatial datasets. Composite
likelihood methods replace the full likelihood by an objective function built
from lower-dimensional marginal or conditional likelihoods
\cite{Varin11}. They are particularly useful when the full likelihood is
computationally expensive or difficult to specify, and they often provide a
good compromise between statistical efficiency and computational tractability.

Only a few results are available for maximum composite likelihood estimators of
Gaussian random fields under fixed-domain asymptotics, and they all concern the
one-dimensional case. \cite{Bachoc19} studied the estimation of covariance
parameters for a Gaussian process with exponential covariance function. They
showed that weighted pairwise likelihood estimators of the microergodic
parameter may be either consistent or inconsistent depending on the objective
function. In particular, the weighted pairwise conditional likelihood estimator
is always consistent and asymptotically Gaussian. \cite{Bachoc&Lagnoux20}
considered a Gaussian process whose covariance function is parametrized by a
variance, a scale, and a smoothness parameter. Their composite likelihood
criteria are based on conditional log-likelihoods of observations given a fixed
number of left and right nearest neighbours. They considered both the case where
only the variance parameter is unknown and the case where the variance and
scale parameters are estimated jointly. In the first case, they showed that for
small values of the smoothness parameter the composite likelihood estimator
converges at a sub-optimal rate and has a non-Gaussian asymptotic distribution,
whereas for larger smoothness values the optimal rate is recovered.

\subsection{Fixed-domain asymptotics for non-Gaussian random fields}

There are comparatively few results on maximum likelihood or maximum composite
likelihood estimation for non-Gaussian random fields under fixed-domain
asymptotics. For instance, \cite{Li13} proposed approximate maximum likelihood
estimation for one-dimensional diffusion processes and derived closed-form
asymptotic expansions for transition densities. Such models, however, do not
provide a general framework for spatial random fields in dimension $d\geq 2$.

Other contributions have focused on variogram-based or power-variation-based
estimators. \cite{Chan&Wood04} considered random fields of the form $g(X)$,
where $g:\mathbb{R}\to\mathbb{R}$ is an unknown smooth function and $X$ is a
real-valued stationary Gaussian field on $\mathbb{R}^{d}$, with $d=1$ or $2$,
whose covariance function follows a power law at the origin. They studied the
asymptotic behaviour of variogram-based estimators when $g(X)$, rather than
$X$, is observed under fixed-domain asymptotics. They showed that the limiting
theory is richer in the non-affine case than in the Gaussian case, where $g$ is
affine. Although these estimators are not likelihood-based, this work
illustrates that fixed-domain asymptotics for non-Gaussian random fields may
differ substantially from the Gaussian case.

For max-stable processes, \cite{Robert20} studied realized power variations of
a class of one-dimensional Brown--Resnick max-stable processes whose spectral
processes are continuous exponential martingales. In that framework, a
fixed-domain asymptotic theory is obtained for sums of powers of absolute
increments. The limiting results involve bias terms depending on local times of
differences between the logarithms of the underlying spectral processes. This
already indicates that local times may play a central role in fixed-domain
asymptotics for Brown--Resnick-type models.

\subsection{Max-stable random fields}

Max-stable random fields arise as the only possible non-degenerate limits for
normalized pointwise maxima of independent and identically distributed random
fields with continuous sample paths; see, for instance, \cite{Haan&Ferreira06}.
Their one-dimensional marginal distributions belong to the generalized extreme
value family. Since the present paper focuses on the estimation of dependence
parameters, we restrict attention to max-stable random fields
$\eta=(\eta(x))_{x\in\mathcal{X}}$ on
$\mathcal{X}\subset\mathbb{R}^{d}$ with standard unit Fr\'echet margins, that
is,
\begin{equation*}
\mathbb{P}\{\eta(x)\leq z\}
=
\exp(-z^{-1}),
\qquad x\in\mathcal{X},\ z>0.
\end{equation*}
The max-stability property then takes the form
\begin{equation*}
n^{-1}\bigvee_{i=1}^{n}\eta_i
\overset{d}{=}
\eta,
\end{equation*}
where $(\eta_i)_{1\leq i\leq n}$ are independent copies of $\eta$,
$\bigvee$ denotes the pointwise maximum, and $\overset{d}{=}$ denotes equality
of finite-dimensional distributions.

Max-stable random fields admit a spectral representation
\cite{Haan84,Gine90}. Any stochastically continuous max-stable process can be
written as
\begin{equation}
\eta(x)
=
\bigvee_{i\geq 1} U_i Y_i(x),
\qquad x\in\mathcal{X},
\label{Eq_Spectral_representation}
\end{equation}
where $(U_i)_{i\geq 1}$ is the decreasing enumeration of the points of a
Poisson point process on $(0,\infty)$ with intensity measure
$u^{-2}\mathrm{d}u$, and $(Y_i)_{i\geq 1}$ are independent copies of a
non-negative stochastic process $Y$ such that
$\mathbb{E}\{Y(x)\}=1$ for all $x\in\mathcal{X}$. The sequences
$(U_i)_{i\geq 1}$ and $(Y_i)_{i\geq 1}$ are independent.

The spectral representation \eqref{Eq_Spectral_representation} also induces a
canonical tessellation of $\mathcal{X}$; see \cite{Dombry&Kabluchko18}. The
cell associated with the index $i\geq 1$ is defined by
\[
C_i
=
\{x\in\mathcal{X}: U_iY_i(x)=\eta(x)\}.
\]
It is a possibly empty random closed subset of $\mathcal{X}$. For each fixed
$x\in\mathcal{X}$, the point process $\{U_iY_i(x)\}_{i\geq 1}$ is a Poisson
point process with intensity $u^{-2}\mathrm{d}u$, so that the maximum
$\eta(x)$ is almost surely attained by a unique index. Hence each point belongs
almost surely to a unique cell. The terms \textit{cell} and
\textit{tessellation} are used here in a broader sense than in classical
stochastic geometry: cells need not be convex or connected, and the
tessellation is a random covering by closed sets with pairwise disjoint
interiors.

Likelihood inference for max-stable random fields is challenging because their
finite-dimensional densities are unknown or difficult to compute in moderate to
high dimension. \cite{Padoan10} proposed a composite likelihood approach based
on low-dimensional marginal densities. The asymptotic properties of the
resulting estimators were studied in a replicated framework, where the
observation sites are fixed and many independent temporal replications of the
spatial field are available.

\subsection{Contributions of the paper}

We consider spatial Brown--Resnick random fields in dimension $d=2$, generated
by isotropic fractional Brownian fields as in \cite{Kabluchko09}. The
underlying fractional Brownian field has semi-variogram
\[
\gamma(x)
=
\frac{\sigma^2\|x\|^\alpha}{2},
\qquad x\in\mathbb{R}^{2},
\]
where $\sigma>0$ is a scale parameter and $\alpha\in(0,2)$ is a smoothness
parameter, with Hurst index $H=\alpha/2$.

Since spatial extreme data are rarely observed on regular grids, we work with a
random sampling design. The observation sites are given by a homogeneous
Poisson point process on a fixed bounded window, independent of the max-stable
field. The Poisson--Delaunay triangulation is then used to select the local
pairs and triples entering the composite likelihood. More precisely, the
pairwise criterion only retains Delaunay edges, whereas the triplewise criterion
only retains the vertices of Delaunay triangles.

This choice is natural in the present framework. First, the composite
likelihood is built from pairwise and triplewise marginal densities, so that
local geometric neighbourhoods provide the most informative contributions under
infill asymptotics. Second, the Delaunay triangulation gives an intrinsic way
of selecting local configurations for irregularly spaced observations, without
introducing an additional deterministic cut-off distance. It is also the
triangulation that maximizes the minimum angle of the triangles, and is
therefore one of the most regular triangulations associated with a given point
configuration.

The first contribution of the paper is to establish limit theorems for squared
increment sums of the logarithm of a Brown--Resnick random field over the edges
and triangles of a Poisson--Delaunay triangulation. These results extend the
limit theorems obtained in \cite{Chenavier&Robert25a} for isotropic fractional
Brownian fields and in \cite{Chenavier&Robert25b} for the pointwise maximum of
two independent isotropic fractional Brownian fields. In the Brown--Resnick
case, the limiting quantities are no longer Gaussian. They are expressed in
terms of local times associated with the interfaces of the canonical
max-stable tessellation; see Theorem \ref{prop:BRtrajectories}.

The second contribution is to use these increment limit theorems to derive the
asymptotic properties of pairwise and triplewise maximum composite likelihood
estimators. We prove consistency of the estimators of the scale and smoothness
parameters when the other parameter is known. We also identify their
non-standard rates of convergence and their centered non-Gaussian limits; see
Theorem \ref{Prop:Asym_Prop_CL_Est}. These limits are specific to the
fixed-domain setting and to the local structure of the Brown--Resnick
tessellation.

These results have direct statistical implications. In a replicated framework,
or under increasing-domain asymptotics, composite likelihood estimators are
usually expected to have asymptotically Gaussian fluctuations. By contrast, in
the fixed-domain framework considered here, the limiting distributions are
driven by local times and cannot be consistently estimated from a single
realization of the field. This shows that standard Gaussian uncertainty
quantification may be inappropriate when only one dense spatial observation of
a max-stable process is available.

Throughout the paper, we restrict attention to isotropic fractional Brownian
fields with Hurst index $H=\alpha/2\in(0,1/2)$, equivalently
$\alpha\in(0,1)$. This is the range in which the limit theorems used in the
paper are established. It is also the range most often encountered in empirical
applications of Brown--Resnick models to spatial extremes; see, for instance,
\cite{Davison12,Engelke14,Einmahl15,Fondeville&Davison18}.

The paper is organized as follows. Section \ref{sec:preliminaries} introduces
Brown--Resnick random fields, recalls the main geometric properties of
Poisson--Delaunay triangulations, and defines the local times used in the
asymptotic analysis. Section \ref{sec:increments_BR} studies normalized
increments of the logarithm of the Brown--Resnick field and establishes the
limit theorems for squared increment sums. Section \ref{sec:CLapproach}
introduces the pairwise and triplewise composite likelihood estimators and
states their asymptotic properties. Section \ref{sec:Simulation_study}
presents a simulation study illustrating the finite-sample behaviour of the
estimators. Proofs and auxiliary technical results are deferred in the Supplementary material, Sections
\ref{sec:proofs} and \ref{sec:intermediary_results}.

\section{Preliminaries}
\label{sec:preliminaries}

\subsection{Max-stable Brown--Resnick random fields}

We consider the class of max-stable random fields known as Brown--Resnick
random fields. These fields are constructed from Gaussian random fields with
stationary increments. Recall that a random field
$(W(x))_{x\in\mathbb{R}^{d}}$ is said to have stationary increments if, for
every $x_0\in\mathbb{R}^{d}$, the law of
\[
\bigl(W(x+x_0)-W(x_0)\bigr)_{x\in\mathbb{R}^{d}}
\]
does not depend on the choice of $x_0$.

A central example is the isotropic fractional Brownian field. In this case,
$W(0)=0$ almost surely and the semi-variogram is given by
\begin{equation}
\gamma(x)
=
\frac{\operatorname{Var}(W(x))}{2}
=
\frac{\sigma^2 \|x\|^\alpha}{2},
\qquad x\in\mathbb{R}^{d},
\label{eq:defvar}
\end{equation}
for some $\alpha\in(0,2)$ and $\sigma^2>0$, where $\|\cdot\|$ denotes the
Euclidean norm. The parameter $\sigma$ is a scale parameter, while $\alpha$ is
a smoothness parameter. Equivalently, $H=\alpha/2$ is the Hurst index, which
determines the local regularity of the sample paths. The isotropic fractional
Brownian field is self-similar and has stationary increments, as described for
instance in Definition 3.3.1 of \cite{Cohen&Istas13}. It should not be
confused with the fractional Brownian sheet, which is a self-similar random
field with stationary rectangular increments; see Section 3.3.2 of
\cite{Cohen&Istas13}.

The Brown--Resnick random field introduced in \cite{Kabluchko09} is the
max-stable random field obtained by taking, in the spectral representation
\eqref{Eq_Spectral_representation},
\begin{equation}
Y(x)
=
\exp\{W(x)-\gamma(x)\},
\qquad x\in\mathbb{R}^{2},
\label{Eq_Def_Y}
\end{equation}
where $W$ is an isotropic fractional Brownian field with semi-variogram
\eqref{eq:defvar}. Since
$\mathbb{E}\{\exp(W(x)-\gamma(x))\}=1$, the normalization required in the
spectral representation is satisfied. With this choice, the resulting
max-stable field $\eta$ is stationary, even though the Gaussian field $W$ is
not stationary itself. The stationarity of $\eta$ follows from the stationary
increments of $W$; see \cite{Kabluchko09}.

\subsection{Delaunay triangulation}

We recall the notation and basic facts on Poisson--Delaunay triangulations
used throughout the paper.

Let $P_N$ be a homogeneous Poisson point process with intensity $N$ on
$\mathbb{R}^{2}$. The Delaunay triangulation of $P_N$, denoted by
$\operatorname{Del}(P_N)$, is almost surely the unique triangulation with
vertices in $P_N$ such that the circumdisk of each triangle contains no point
of $P_N$ in its interior; see, for instance, p.~478 in
\cite{Schneider&Weil08}. The Delaunay triangulation is often viewed as the
most regular triangulation associated with a point configuration, in the sense
that it maximizes the minimum angle of the triangles.

We first recall the notion of typical cell for the Delaunay triangulation
associated with a homogeneous Poisson point process $P_1$ of intensity one.
For each cell $C\in\operatorname{Del}(P_1)$, let $z(C)$ denote its
circumcenter. If $\mathbf{B}\subset\mathbb{R}^{2}$ is a Borel set with
Lebesgue measure $|\mathbf{B}|\in(0,\infty)$, the cell intensity
$\beta_2$ of $\operatorname{Del}(P_1)$ is defined by
\[
\beta_2
=
\frac{1}{|\mathbf{B}|}
\mathbb{E}\left[
\left|\{C\in\operatorname{Del}(P_1):z(C)\in\mathbf{B}\}\right|
\right].
\]
It is well known that $\beta_2=2$; see Theorem 10.2.9 in
\cite{Schneider&Weil08}. The typical cell is the random triangle
$\mathcal{C}$ whose distribution is characterized as follows: for every
positive measurable translation-invariant function
$g:\mathcal{K}_2\to\mathbb{R}$,
\[
\mathbb{E}\left[g(\mathcal{C})\right]
=
\frac{1}{\beta_2|\mathbf{B}|}
\mathbb{E}\left[
\sum_{C\in\operatorname{Del}(P_1):z(C)\in\mathbf{B}} g(C)
\right],
\]
where $\mathcal{K}_2$ denotes the set of compact convex subsets of
$\mathbb{R}^{2}$, endowed with the Fell topology; see Section 12.2 in
\cite{Schneider&Weil08}.

Let $\Delta(y_1,y_2,y_3)$ denote the convex hull of the points
$y_1,y_2,y_3\in\mathbb{R}^{2}$, and let
$|\Delta(y_1,y_2,y_3)|$ be its area. The distribution of the typical cell has
the integral representation
\begin{equation}
\mathbb{E}\{g(\mathcal{C})\}
=
\frac{1}{6}
\int_{0}^{\infty}
\int_{(\mathbb{S}^{1})^{3}}
r^{3}e^{-\pi r^{2}}
|\Delta(u_1,u_2,u_3)|
g\bigl(\Delta(ru_1,ru_2,ru_3)\bigr)
\,\lambda(\mathrm{d}u_1)\lambda(\mathrm{d}u_2)\lambda(\mathrm{d}u_3)
\,\mathrm{d}r,
\label{eq:typicalcell}
\end{equation}
where $\mathbb{S}^{1}$ is the unit sphere of $\mathbb{R}^{2}$ and
$\lambda$ denotes the spherical Lebesgue measure on $\mathbb{S}^{1}$,
normalized by $\lambda(\mathbb{S}^{1})=2\pi$. Equivalently,
$\mathcal{C}$ has the same distribution as
$R\Delta(\mathcal{U}_1,\mathcal{U}_2,\mathcal{U}_3)$, where $R$ and
$(\mathcal{U}_1,\mathcal{U}_2,\mathcal{U}_3)$ are independent, with
densities proportional respectively to
$2\pi^2 r^3e^{-\pi r^2}$ and
$|\Delta(u_1,u_2,u_3)|/(12\pi^2)$.

The notion of typical edge is defined similarly. The edge intensity
$\beta_1$ of $\operatorname{Del}(P_1)$ is the mean number of edges per unit
area and satisfies $\beta_1=3$; see again Theorem 10.2.9 in
\cite{Schneider&Weil08}. The length of the typical edge has the same
distribution as
\[
D=R\|\mathcal{U}_1-\mathcal{U}_2\|.
\]
Its distribution is characterized by
\begin{align}
\mathbb{P}\{D\leq \ell\}
&=
\int_{0}^{\ell} f_D(s)\,\mathrm{d}s
\notag\\
&=
\frac{\pi}{3}
\int_{0}^{\infty}
\int_{(\mathbb{S}^{1})^{2}}
r^{3}e^{-\pi r^{2}}
|\Delta(u_1,u_2,e_1)|
\mathbb{I}\{r\|u_1-u_2\|\leq \ell\}
\,\lambda(\mathrm{d}u_1)\lambda(\mathrm{d}u_2)
\,\mathrm{d}r,
\label{eq:typicallength}
\end{align}
where $e_1=(1,0)$ and $\ell>0$.

Following the representation \eqref{eq:typicalcell}, one may also define a
typical pair of distinct Delaunay edges sharing a common vertex. It is
represented by a random vector $(D_1,D_2,\Theta)$, where
$D_1,D_2\geq0$ are the two edge lengths and
$\Theta\in[-\pi/2,\pi/2)$ is the corresponding angle. Its distribution is
given by
\begin{multline*}
\mathbb{P}\{(D_1,D_2,\Theta)\in B\}
=
\frac{1}{6}
\int_{0}^{\infty}
\int_{(\mathbb{S}^{1})^{3}}
r^{3}e^{-\pi r^{2}}
|\Delta(u_1,u_2,u_3)|
\\
\times
\mathbb{I}\Bigl\{
\bigl(
r\|u_3-u_2\|,
r\|u_2-u_1\|,
\arcsin(\cos(\theta_{u_1,u_2}/2))
\bigr)
\in B
\Bigr\}
\,\lambda(\mathrm{d}u_1)\lambda(\mathrm{d}u_2)\lambda(\mathrm{d}u_3)
\,\mathrm{d}r,
\end{multline*}
for every Borel set $B\subset\mathbb{R}_+^2\times[-\pi/2,\pi/2)$, where
$\theta_{u_1,u_2}$ denotes the angle between $u_1$ and $u_2$. In particular,
the marginal distribution of each edge length is the distribution of the
typical edge length $D$ given in \eqref{eq:typicallength}.

Throughout the paper, we identify $\operatorname{Del}(P_N)$ with its
skeleton. When two points $x_1,x_2\in P_N$ are Delaunay neighbours, we write
\[
x_1\sim x_2
\quad\text{in }\operatorname{Del}(P_N).
\]
For a Borel set $\mathbf{B}\subset\mathbb{R}^{2}$, let $E_{N,\mathbf{B}}$
be the set of ordered pairs $(x_1,x_2)$ such that
\[
x_1\sim x_2 \text{ in } \operatorname{Del}(P_N),
\qquad
x_1\in\mathbf{B},
\qquad
x_1\preceq x_2,
\]
where $\preceq$ denotes the lexicographic order. When
$\mathbf{B}=\mathbf{C}:=(-1/2,1/2]^2$, we simply write
\[
E_N:=E_{N,\mathbf{C}}.
\]

Similarly, for a Borel set $\mathbf{B}\subset\mathbb{R}^{2}$, let
$DT_{N,\mathbf{B}}$ be the set of ordered triples $(x_1,x_2,x_3)$ such that
\[
\Delta(x_1,x_2,x_3)\in\operatorname{Del}(P_N),
\qquad
x_1\in\mathbf{B},
\qquad
x_1\preceq x_2\preceq x_3.
\]
When $\mathbf{B}=\mathbf{C}$, we write
\[
DT_N:=DT_{N,\mathbf{C}}.
\]

\subsection{Local time}

Let $W^{(1)}$ and $W^{(2)}$ be two independent copies of the isotropic
fractional Brownian field on $\mathbb{R}^{2}$ with semi-variogram
\eqref{eq:defvar}. We write
\[
W^{(2\setminus 1)}(x)
=
W^{(2)}(x)-W^{(1)}(x),
\qquad x\in\mathbb{R}^{2}.
\]
The local time of this difference field measures, in an occupation-density
sense, the size of the set on which the two fields are close to a prescribed
level. In particular, the local time at level zero will be used later to
describe the contribution of interfaces between competing spectral functions.

Let $\nu^{(2\setminus 1)}$ be the occupation measure of
$W^{(2\setminus 1)}$ over $\mathbf{C}=(-1/2,1/2]^2$, defined by
\[
\nu^{(2\setminus 1)}(A)
=
\int_{\mathbf{C}}
\mathbb{I}\{W^{(2\setminus 1)}(x)\in A\}\,\mathrm{d}x,
\]
for every Borel set $A\subset\mathbb{R}$. For any
$s,t\in\mathbb{R}^{2}$,
\[
\Delta(s,t)
:=
\mathbb{E}\left[
\bigl(W^{(2\setminus 1)}(s)-W^{(2\setminus 1)}(t)\bigr)^2
\right]
=
2\sigma^2\|s-t\|^\alpha.
\]
Since
\[
\int_{\mathbf{C}} \Delta(s,t)^{-1/2}\,\mathrm{d}s <\infty,
\qquad t\in\mathbf{C},
\]
it follows from Section 22 of \cite{Geman&Horowitz80} that the occupation
measure $\nu^{(2\setminus 1)}$ admits a Lebesgue density. The local time at
level $\ell\in\mathbb{R}$ is therefore defined by
\[
L_{W^{(2\setminus 1)}}(\ell)
:=
\frac{\mathrm{d}\nu^{(2\setminus 1)}}{\mathrm{d}\ell}(\ell).
\]

The occupation density formula states that, for every Borel function
$g:\mathbb{R}\to\mathbb{R}$ for which the integrals are well defined,
\[
\int_{\mathbf{C}}
g\bigl(W^{(2\setminus 1)}(x)\bigr)\,\mathrm{d}x
=
\int_{\mathbb{R}}
g(\ell)L_{W^{(2\setminus 1)}}(\ell)\,\mathrm{d}\ell.
\]
Moreover, adapting the proof of Lemma 1.1 in \cite{Jaramillo21}, one obtains
the following $L^2$ representations. For every $\ell\in\mathbb{R}$,
\[
L_{W^{(2\setminus 1)}}(\ell)
=
\lim_{\varepsilon\to0}
\int_{\mathbf{C}}
\frac{1}{\sqrt{2\pi\varepsilon}}
\exp\left\{
-\frac{(W^{(2\setminus 1)}(x)-\ell)^2}{2\varepsilon}
\right\}
\,\mathrm{d}x,
\]
where the limit holds in $L^2$. Equivalently,
\begin{equation}
L_{W^{(2\setminus 1)}}(\ell)
=
\frac{1}{2\pi}
\lim_{M\to\infty}
\int_{-M}^{M}
\int_{\mathbf{C}}
\exp\left\{
\mathrm{i}\xi\bigl(W^{(2\setminus 1)}(x)-\ell\bigr)
\right\}
\,\mathrm{d}x\,\mathrm{d}\xi,
\label{eq:localtimeintegral}
\end{equation}
again with convergence in $L^2$.

\section{Increments of a Brown--Resnick random field}
\label{sec:increments_BR}

\subsection{Asymptotic distributions of normalized increments}

We first study the local behaviour of normalized increments of the logarithm
of a Brown--Resnick random field. These elementary distributional results will
be used later to identify the leading terms in the composite likelihood score
functions and in the squared increment sums.

\subsubsection{Single increments}

Let $x_1,x_2\in\mathbb{R}^{2}$ be two distinct sites and write
\[
d=\|x_2-x_1\|.
\]
For $z_1,z_2>0$, set
\[
a=\sigma d^{\alpha/2},
\qquad
u=\frac{\log(z_2/z_1)}{a},
\qquad
v(t)=\frac{a}{2}+t,\quad t\in\mathbb{R}.
\]
The bivariate distribution function of
$(\eta(x_1),\eta(x_2))$ is given by
\[
F_{x_1,x_2}(z_1,z_2)
=
\mathbb{P}\{\eta(x_1)\leq z_1,\eta(x_2)\leq z_2\}
=
\exp\{-V_{x_1,x_2}(z_1,z_2)\},
\]
where the pairwise exponent function is
\begin{equation}
V_{x_1,x_2}(z_1,z_2)
=
\frac{1}{z_1}\Phi\{v(u)\}
+
\frac{1}{z_2}\Phi\{v(-u)\},
\qquad z_1,z_2>0.
\label{Eq_V_dim_2}
\end{equation}
Here $\Phi$ denotes the distribution function of the standard Gaussian
distribution; see, for instance, \cite{Huser&Davison13}.

We consider the normalized increment of $\log \eta$ between $x_1$ and $x_2$,
defined by
\[
U_{x_1,x_2}^{(\eta)}
=
\frac{1}{\sigma d^{\alpha/2}}
\log\left(\frac{\eta(x_2)}{\eta(x_1)}\right).
\]
The following proposition gives both the conditional and marginal
distributions of this increment. It also shows that, at small spatial scales,
the increment behaves as the corresponding normalized increment of the
underlying fractional Brownian field.

\begin{proposition}
\label{Prop_cond_dist_U}
For any $y>0$, the conditional distribution of
$U_{x_1,x_2}^{(\eta)}$ given $\eta(x_1)=y$ is
\[
\mathbb{P}
\left\{
U_{x_1,x_2}^{(\eta)}\leq u
\,\middle|\,
\eta(x_1)=y
\right\}
=
\exp\left(
-\frac{1}{y}
\left[
V_{x_1,x_2}
\left(1,e^{\sigma d^{\alpha/2}u}\right)
-1
\right]
\right)
\Phi\{v(u)\},
\qquad u\in\mathbb{R}.
\]
Its marginal distribution is
\[
\mathbb{P}\left\{
U_{x_1,x_2}^{(\eta)}\leq u
\right\}
=
\frac{\Phi\{v(u)\}}
{
V_{x_1,x_2}
\left(1,e^{\sigma d^{\alpha/2}u}\right)
},
\qquad u\in\mathbb{R}.
\]
Consequently,
\[
\lim_{d\to0}
\mathbb{P}\left\{
U_{x_1,x_2}^{(\eta)}\leq u
\right\}
=
\Phi(u),
\qquad u\in\mathbb{R}.
\]
\end{proposition}

This Gaussian limit is natural. As $d\to0$, the probability that the two
points $x_1$ and $x_2$ belong to the same cell of the canonical tessellation
of the Brown--Resnick field tends to one. On such an event, the two values of
the max-stable field are generated by the same spectral function, and the
local increment of $\log\eta$ therefore coincides with an increment of an
isotropic fractional Brownian field. Proposition
\ref{Prop_cond_dist_U} extends Proposition 3 of \cite{Robert20} to the
spatial Brown--Resnick setting considered here.

\subsubsection{Pairs of increments}

We now consider three sites $x_1,x_2,x_3\in\mathbb{R}^{2}$. For
$i,j\in\{1,2,3\}$, $i\neq j$, write
\[
d_{i,j}=\|x_j-x_i\|,
\qquad
a_{i,j}=\sigma d_{i,j}^{\alpha/2},
\qquad
u_{i,j}=\frac{\log(z_j/z_i)}{a_{i,j}},
\]
and define
\[
v_{i,j}(t)=\frac{a_{i,j}}{2}+t,
\qquad t\in\mathbb{R}.
\]
The trivariate distribution function of
$(\eta(x_1),\eta(x_2),\eta(x_3))$ is
\[
F_{x_1,x_2,x_3}(z_1,z_2,z_3)
=
\mathbb{P}
\{\eta(x_1)\leq z_1,\eta(x_2)\leq z_2,\eta(x_3)\leq z_3\}
=
\exp\{-V_{x_1,x_2,x_3}(z_1,z_2,z_3)\},
\]
where, see for instance \cite{Huser&Davison13},
\begin{multline}
V_{x_1,x_2,x_3}(z_1,z_2,z_3)
=
\frac{1}{z_1}
\Phi_2
\left(
\begin{pmatrix}
v_{1,2}(u_{1,2})\\
v_{1,3}(u_{1,3})
\end{pmatrix};
\begin{pmatrix}
1 & R_1\\
R_1 & 1
\end{pmatrix}
\right)
\\
+
\frac{1}{z_2}
\Phi_2
\left(
\begin{pmatrix}
v_{1,2}(-u_{1,2})\\
v_{2,3}(u_{2,3})
\end{pmatrix};
\begin{pmatrix}
1 & R_2\\
R_2 & 1
\end{pmatrix}
\right)
\\
+
\frac{1}{z_3}
\Phi_2
\left(
\begin{pmatrix}
v_{1,3}(-u_{1,3})\\
v_{2,3}(-u_{2,3})
\end{pmatrix};
\begin{pmatrix}
1 & R_3\\
R_3 & 1
\end{pmatrix}
\right).
\label{Eq_V_triple}
\end{multline}
Here $\Phi_2(\cdot;\Sigma)$ denotes the distribution function of the centered
bivariate Gaussian distribution with covariance matrix $\Sigma$, and
\[
R_1
=
\frac{
d_{1,2}^{\alpha}+d_{1,3}^{\alpha}-d_{2,3}^{\alpha}
}
{
2(d_{1,2}d_{1,3})^{\alpha/2}
},
\qquad
R_2
=
\frac{
d_{1,2}^{\alpha}+d_{2,3}^{\alpha}-d_{1,3}^{\alpha}
}
{
2(d_{1,2}d_{2,3})^{\alpha/2}
},
\]
\[
R_3
=
\frac{
d_{1,3}^{\alpha}+d_{2,3}^{\alpha}-d_{1,2}^{\alpha}
}
{
2(d_{1,3}d_{2,3})^{\alpha/2}
}.
\]
The coefficient $R_1$ is the correlation between the normalized increments of
the underlying fractional Brownian field from $x_1$ to $x_2$ and from $x_1$
to $x_3$.

We consider the vector of normalized increments
\[
U_{x_1,x_2}^{(\eta)}
=
\frac{1}{\sigma d_{1,2}^{\alpha/2}}
\log\left(\frac{\eta(x_2)}{\eta(x_1)}\right),
\qquad
U_{x_1,x_3}^{(\eta)}
=
\frac{1}{\sigma d_{1,3}^{\alpha/2}}
\log\left(\frac{\eta(x_3)}{\eta(x_1)}\right).
\]
The following proposition is the bivariate analogue of Proposition
\ref{Prop_cond_dist_U}.

\begin{proposition}
\label{Prop_cond_dist_U_1_U_2}
For any $y>0$, the conditional distribution of
$(U_{x_1,x_2}^{(\eta)},U_{x_1,x_3}^{(\eta)})$ given $\eta(x_1)=y$ is
\begin{multline*}
\mathbb{P}
\left\{
U_{x_1,x_2}^{(\eta)}\leq u_2,
U_{x_1,x_3}^{(\eta)}\leq u_3
\,\middle|\,
\eta(x_1)=y
\right\}
\\
=
\exp\left(
-\frac{1}{y}
\left[
V_{x_1,x_2,x_3}
\left(
1,
e^{\sigma d_{1,2}^{\alpha/2}u_2},
e^{\sigma d_{1,3}^{\alpha/2}u_3}
\right)
-1
\right]
\right)
\\
\times
\Phi_2
\left(
\begin{pmatrix}
v_{1,2}(u_2)\\
v_{1,3}(u_3)
\end{pmatrix};
\begin{pmatrix}
1 & R_1\\
R_1 & 1
\end{pmatrix}
\right),
\qquad u_2,u_3\in\mathbb{R}.
\end{multline*}
Its marginal distribution is
\begin{equation*}
\mathbb{P}
\left\{
U_{x_1,x_2}^{(\eta)}\leq u_2,
U_{x_1,x_3}^{(\eta)}\leq u_3
\right\}
=
\frac{
\Phi_2
\left(
\begin{pmatrix}
v_{1,2}(u_2)\\
v_{1,3}(u_3)
\end{pmatrix};
\begin{pmatrix}
1 & R_1\\
R_1 & 1
\end{pmatrix}
\right)
}
{
V_{x_1,x_2,x_3}
\left(
1,
e^{\sigma d_{1,2}^{\alpha/2}u_2},
e^{\sigma d_{1,3}^{\alpha/2}u_3}
\right)
},
\qquad u_2,u_3\in\mathbb{R}.
\end{equation*}
Assume now that
\[
\|x_2-x_1\|=\delta d_{1,2},
\qquad
\|x_3-x_1\|=\delta d_{1,3},
\qquad
\|x_3-x_2\|=\delta d_{2,3},
\]
where $d_{1,2}$, $d_{1,3}$ and $d_{2,3}$ are fixed positive numbers. Then
\[
\lim_{\delta\to0}
\mathbb{P}
\left\{
U_{x_1,x_2}^{(\eta)}\leq u_2,
U_{x_1,x_3}^{(\eta)}\leq u_3
\right\}
=
\Phi_2
\left(
\begin{pmatrix}
u_2\\
u_3
\end{pmatrix};
\begin{pmatrix}
1 & R_1\\
R_1 & 1
\end{pmatrix}
\right),
\qquad u_2,u_3\in\mathbb{R}.
\]
\end{proposition}

The same interpretation as in the one-dimensional case applies. As
$\delta\to0$, the probability that $x_1$, $x_2$ and $x_3$ belong to the same
cell of the canonical tessellation tends to one. Hence the vector
$(U_{x_1,x_2}^{(\eta)},U_{x_1,x_3}^{(\eta)})$ has the same limiting
distribution as the corresponding vector of normalized increments of an
isotropic fractional Brownian field.

\subsection{Limit behaviour of squared increment sums}
\label{subsec:squared_increment_sums_BR}

We now turn to the asymptotic behaviour of squared increment sums of
$\log\eta$ over the edges and triangles of the Poisson--Delaunay
triangulation. Throughout this subsection, let
$(\eta(x))_{x\in\mathbb{R}^{2}}$ be a Brown--Resnick random field with
spectral representation
\[
\eta(x)
=
\bigvee_{i\geq1} U_iY_i(x),
\qquad x\in\mathbb{R}^{2},
\]
where $(U_i)_{i\geq1}$ is the decreasing enumeration of the points of a
Poisson point process on $(0,\infty)$ with intensity measure
$u^{-2}\mathrm{d}u$, and $(Y_i)_{i\geq1}$ are independent copies of
\[
Y(x)
=
\exp\{W(x)-\gamma(x)\},
\qquad x\in\mathbb{R}^{2}.
\]
Here $W$ is an isotropic fractional Brownian field satisfying $W(0)=0$
almost surely and
\[
\gamma(x)
=
\frac{\operatorname{Var}(W(x))}{2}
=
\frac{\sigma^2\|x\|^\alpha}{2},
\qquad x\in\mathbb{R}^{2},
\]
for some $\alpha\in(0,1)$ and $\sigma^2>0$.

For $i\geq1$, set
\[
Z_i(x)
=
\log U_i+\log Y_i(x),
\qquad x\in\mathbb{R}^{2}.
\]
Thus
\[
\log\eta(x)=\bigvee_{i\geq1}Z_i(x).
\]
For $k\neq j$, define the difference field
\[
Z_{k\setminus j}(x)
=
Z_k(x)-Z_j(x),
\qquad x\in\mathbb{R}^{2}.
\]
For $k>j\geq1$, let
\begin{equation}
\mathbf{C}_{k,j}
=
\left\{
x\in\mathbf{C}:
Z_k(x)\wedge Z_j(x)
>
\bigvee_{i\neq j,k}Z_i(x)
\right\}.
\label{eq:def_Ckj}
\end{equation}
The set $\mathbf{C}_{k,j}$ is the region of the observation window
$\mathbf{C}=(-1/2,1/2]^2$ where the two spectral functions $Z_j$ and $Z_k$
are the two largest ones. On this region, changes in the maximizing spectral
index can only occur through crossings of $Z_k-Z_j$. This is why the local
time of $Z_{k\setminus j}$ at level zero appears in the limiting behaviour
of squared increment sums.

If $\mathbf{C}_{k,j}\neq\varnothing$, define the occupation measure of
$Z_{k\setminus j}$ over $\mathbf{C}_{k,j}$ by
\[
\nu^{(k\setminus j)}(A)
=
\int_{\mathbf{C}_{k,j}}
\mathbb{I}\{Z_{k\setminus j}(x)\in A\}\,\mathrm{d}x,
\]
for every Borel set $A\subset\mathbb{R}$. The associated local time at level
zero is denoted by
\[
L_{Z_{k\setminus j}}(0)
=
\frac{\mathrm{d}\nu^{(k\setminus j)}}{\mathrm{d}\ell}(0).
\]
If $\mathbf{C}_{k,j}=\varnothing$, we set
\[
L_{Z_{k\setminus j}}(0)=0.
\]

We consider the following centered squared increment sums:
\begin{align}
V_{2,N}^{(\eta)}
&=
\frac{1}{\sqrt{|E_N|}}
\sum_{(x_1,x_2)\in E_N}
\left[
\left(U_{x_1,x_2}^{(\eta)}\right)^2-1
\right],
\label{eq:def_V2_eta}
\\
V_{3,N}^{(\eta)}
&=
\frac{1}{\sqrt{|DT_N|}}
\sum_{(x_1,x_2,x_3)\in DT_N}
\left[
\begin{pmatrix}
U_{x_1,x_2}^{(\eta)} &
U_{x_1,x_3}^{(\eta)}
\end{pmatrix}
\begin{pmatrix}
1 & R_{x_1,x_2,x_3}\\
R_{x_1,x_2,x_3} & 1
\end{pmatrix}^{-1}
\begin{pmatrix}
U_{x_1,x_2}^{(\eta)}\\
U_{x_1,x_3}^{(\eta)}
\end{pmatrix}
-2
\right],
\label{eq:def_V3_eta}
\end{align}
where
\begin{equation}
R_{x_1,x_2,x_3}
=
\frac{
d_{1,2}^{\alpha}+d_{1,3}^{\alpha}-d_{2,3}^{\alpha}
}
{
2(d_{1,2}d_{1,3})^{\alpha/2}
},
\label{eq:coorU}
\end{equation}
with
\[
d_{1,2}=\|x_2-x_1\|,
\qquad
d_{1,3}=\|x_3-x_1\|,
\qquad
d_{2,3}=\|x_3-x_2\|.
\]
The statistic $V_{2,N}^{(\eta)}$ is based on Delaunay edges, whereas
$V_{3,N}^{(\eta)}$ is based on Delaunay triangles and uses the quadratic form
associated with the limiting covariance matrix of the two normalized
increments.

Strictly speaking, the statistics above are not defined on the events
\(\{|E_N|=0\}\) and \(\{|DT_N|=0\}\). We shall use the convention that
\(V_{2,N}^{(\eta)}=0\) on \(\{|E_N|=0\}\) and \(V_{3,N}^{(\eta)}=0\) on
\(\{|DT_N|=0\}\). This convention is asymptotically immaterial, since these
exceptional events have exponentially small probability as \(N\to\infty\).

The following theorem identifies the leading asymptotic behaviour of
$V_{2,N}^{(\eta)}$ and $V_{3,N}^{(\eta)}$.

\begin{theorem}
\label{prop:BRtrajectories}
Let $\alpha\in(0,1)$. Then, as $N\to\infty$, there exist real constants
$c_{V_2}$ and $c_{V_3}$ such that
\[
\frac{\sqrt{3}}{3}
N^{-(2-\alpha)/4}
V_{2,N}^{(\eta)}
\overset{\mathbb{P}}{\longrightarrow}
c_{V_2}
\sum_{j\geq1}\sum_{k>j}
L_{Z_{k\setminus j}}(0),
\]
and
\[
\frac{\sqrt{2}}{2}
N^{-(2-\alpha)/4}
V_{3,N}^{(\eta)}
\overset{\mathbb{P}}{\longrightarrow}
c_{V_3}
\sum_{j\geq1}\sum_{k>j}
L_{Z_{k\setminus j}}(0).
\]
\end{theorem}

Theorem \ref{prop:BRtrajectories} extends Theorem 2 of
\cite{Chenavier&Robert25b} to the Brown--Resnick field, whose logarithm is
the pointwise maximum of infinitely many shifted spectral functions generated
by independent isotropic fractional Brownian fields. We note that only an
almost surely finite number of the local times
$L_{Z_{k\setminus j}}(0)$ are positive. This follows from the fact that the
canonical tessellation has almost surely only finitely many non-empty cells in
the compact window $\mathbf{C}$.

Finally, using the Slivnyak--Mecke formula, see for instance Theorem 3.2.5 in
\cite{Schneider&Weil08}, together with the same arguments as in the proof of
Proposition 3 of \cite{Robert20}, one obtains
\[
\lim_{N\to\infty}
N^{\alpha/4}
\mathbb{E}
\left[
\frac{1}{|E_N|}
\sum_{(x_1,x_2)\in E_N}
\left\{
\left(U_{x_1,x_2}^{(\eta)}\right)^2-1
\right\}
\right]
=
4\sigma \mathbb{E}\left[D^{\alpha/2}\right]\psi,
\]
where $D$ denotes the typical Delaunay edge length defined in
\eqref{eq:typicallength}, and
\begin{equation}
\psi
=
\int_{0}^{\infty}
u\varphi(u)
\left[
\frac{1}{2}
-
\bar{\Phi}(u)
-
u\bar{\Phi}(u)\frac{\Phi(u)}{\varphi(u)}
\right]
\,\mathrm{d}u
\simeq
-0.094.
\label{eq:constant_psy}
\end{equation}
Here $\varphi$ is the standard Gaussian density and
$\bar{\Phi}=1-\Phi$. Together with the normalization used in Theorem
\ref{prop:BRtrajectories}, this implies that $c_{V_2}<0$.

\section{The Delaunay-tapered composite likelihood approach}
\label{sec:CLapproach}

We assume that the observation sites are given by a realization of a
homogeneous Poisson point process $P_N$ with intensity $N$ on
$\mathbb{R}^{2}$. The point process $P_N$ is assumed to be independent of the
Brown--Resnick random field. The observation window is
\[
\mathbf{C}=(-1/2,1/2]^2.
\]
The Delaunay triangulation associated with $P_N$ is used to select the local
pairs and triples entering the composite likelihood criteria. Equivalently,
the composite likelihoods below may be viewed as weighted composite
likelihoods with weights given by indicators of Delaunay edges and Delaunay
triangles.

\subsection{Composite likelihood objective functions and estimators}

For any two distinct sites $x_1,x_2\in\mathbb{R}^{2}$, the distribution of
$(\eta(x_1),\eta(x_2))$ is absolutely continuous with respect to Lebesgue
measure on $(0,\infty)^2$. We denote its density by
\[
f_{x_1,x_2}(z_1,z_2;\sigma,\alpha)
=
\frac{\partial^2}{\partial z_1\partial z_2}
F_{x_1,x_2}(z_1,z_2;\sigma,\alpha),
\qquad z_1,z_2>0.
\]
For fixed $z_1,z_2>0$, this density is differentiable with respect to
$(\sigma,\alpha)$ on $(0,\infty)\times(0,2)$. The Delaunay-tapered pairwise
composite log-likelihood is defined by
\[
\ell_{2,N}(\sigma,\alpha)
=
\sum_{(x_1,x_2)\in E_N}
\log f_{x_1,x_2}
\bigl(\eta(x_1),\eta(x_2);\sigma,\alpha\bigr).
\]

Similarly, for any three distinct sites $x_1,x_2,x_3\in\mathbb{R}^{2}$, the
distribution of $(\eta(x_1),\eta(x_2),\eta(x_3))$ is absolutely continuous
with respect to Lebesgue measure on $(0,\infty)^3$. We denote its density by
\[
f_{x_1,x_2,x_3}(z_1,z_2,z_3;\sigma,\alpha)
=
\frac{\partial^3}{\partial z_1\partial z_2\partial z_3}
F_{x_1,x_2,x_3}(z_1,z_2,z_3;\sigma,\alpha),
\qquad z_1,z_2,z_3>0.
\]
For fixed $z_1,z_2,z_3>0$, this density is differentiable with respect to
$(\sigma,\alpha)$ on $(0,\infty)\times(0,2)$. The Delaunay-tapered triplewise
composite log-likelihood is defined by
\[
\ell_{3,N}(\sigma,\alpha)
=
\sum_{(x_1,x_2,x_3)\in DT_N}
\log f_{x_1,x_2,x_3}
\bigl(\eta(x_1),\eta(x_2),\eta(x_3);\sigma,\alpha\bigr).
\]

Thus, in the pairwise criterion, only Delaunay edges are retained; in the
triplewise criterion, only triples forming Delaunay triangles are retained.
This selection focuses the composite likelihood on local configurations, which
are the relevant ones under fixed-domain asymptotics. It also avoids the
introduction of an additional deterministic cut-off distance.

The regularity conditions needed to define the maximizers and to control the
objective functions are standard for composite likelihood inference. In
particular, from Section 4.4 of \cite{Dombry18}, there exist families of
positive functions
\[
\bigl(l_{x_1,x_2}\bigr)_{x_1,x_2\in\mathbb{R}^{2}},
\qquad
\bigl(l_{x_1,x_2,x_3}\bigr)_{x_1,x_2,x_3\in\mathbb{R}^{2}},
\]
with $l_{x_1,x_2}:(0,\infty)^2\to\mathbb{R}_{+}$ and
$l_{x_1,x_2,x_3}:(0,\infty)^3\to\mathbb{R}_{+}$, such that the following
Lipschitz-type bounds hold. For any $\sigma_1,\sigma_2>0$ and
$\alpha_1,\alpha_2\in(0,2)$,
\[
\left|
\log
\frac{
f_{x_1,x_2}(z_1,z_2;\sigma_2,\alpha_2)
}{
f_{x_1,x_2}(z_1,z_2;\sigma_1,\alpha_1)
}
\right|
\leq
l_{x_1,x_2}(z_1,z_2)
\bigl(|\sigma_2-\sigma_1|+|\alpha_2-\alpha_1|\bigr),
\]
and
\[
\left|
\log
\frac{
f_{x_1,x_2,x_3}(z_1,z_2,z_3;\sigma_2,\alpha_2)
}{
f_{x_1,x_2,x_3}(z_1,z_2,z_3;\sigma_1,\alpha_1)
}
\right|
\leq
l_{x_1,x_2,x_3}(z_1,z_2,z_3)
\bigl(|\sigma_2-\sigma_1|+|\alpha_2-\alpha_1|\bigr).
\]

Let $(\sigma_0,\alpha_0)$ denote the true parameter value. We assume that
$\sigma_0$ belongs to a compact set $\mathrm{S}_{\sigma}\subset(0,\infty)$
and that $\alpha_0$ belongs to a compact set
$\mathrm{S}_{\alpha}\subset(0,1)$. We consider one-parameter estimation
problems, in which one parameter is estimated while the other one is assumed
to be known.

When $\alpha_0$ is known, the pairwise and triplewise maximum composite
likelihood estimators of $\sigma_0$ are defined by
\[
\widehat{\sigma}_{j,N}
\in
\arg\max_{\sigma\in\mathrm{S}_{\sigma}}
\ell_{j,N}(\sigma,\alpha_0),
\qquad j=2,3.
\]
When $\sigma_0$ is known, the pairwise and triplewise maximum composite
likelihood estimators of $\alpha_0$ are defined by
\[
\widehat{\alpha}_{j,N}
\in
\arg\max_{\alpha\in\mathrm{S}_{\alpha}}
\ell_{j,N}(\sigma_0,\alpha),
\qquad j=2,3.
\]
The asymptotic results below apply to any sequence of maximizers satisfying
the corresponding first-order optimality conditions. Local asymptotic
uniqueness follows from the score expansions derived in Propositions
\ref{Prop_pairwise_pdf} and \ref{Prop_triplewise_pdf}.

\subsection{Asymptotic score functions}

\subsubsection{The pairwise case}

The following proposition gives the local expansion of the pairwise score
contribution as the distance between the two sites tends to zero.

\begin{proposition}
\label{Prop_pairwise_pdf}
Let $u\in\mathbb{R}$ be fixed. Let $x_1,x_2\in\mathbb{R}^{2}$ and
$z_1,z_2>0$ be such that
\[
u
=
\frac{1}{\sigma d^{\alpha/2}}
\log\left(\frac{z_2}{z_1}\right),
\qquad
d=\|x_2-x_1\|>0.
\]
Then, as $d\to0$,
\[
\frac{\partial}{\partial\sigma}
\log f_{x_1,x_2}(z_1,z_2;\sigma,\alpha)
=
\frac{1}{\sigma}(u^2-1)
+
\frac{\omega(u)}{z_1}d^{\alpha/2}
+
o(d^{\alpha/2}),
\]
and
\[
\frac{1}{\log d}
\frac{\partial}{\partial\alpha}
\log f_{x_1,x_2}(z_1,z_2;\sigma,\alpha)
=
\frac{1}{2}(u^2-1)
+
\frac{\sigma \omega(u)}{2z_1}d^{\alpha/2}
+
o(d^{\alpha/2}),
\]
where
\[
\omega(u)
=
u\{2\Phi(u)-1\}
+
\frac{(1-u^2)\Phi(u)\{1-\Phi(u)\}}{\varphi(u)}
-
\varphi(u).
\]
\end{proposition}

The leading pairwise score contribution is therefore proportional to the
second Hermite polynomial $u^2-1$. In the composite likelihood score, $u$ is
replaced by the normalized increment
$U_{x_1,x_2}^{(\eta)}$, which converges in distribution to a standard
Gaussian random variable as $d\to0$ by Proposition
\ref{Prop_cond_dist_U}. The leading score contribution is thus
asymptotically centered.

The first-order correction is also important, since it determines the
centering of the limiting score. Since the score has mean zero under the true
parameter value, for any fixed $x_1,x_2$,
\[
\mathbb{E}
\left[
\frac{\partial}{\partial\sigma}
\log f_{x_1,x_2}
\bigl(\eta(x_1),\eta(x_2);\sigma,\alpha\bigr)
\right]
=
0.
\]
Consequently, applying the exact score identity at the true parameter value gives
\[
0
=
\frac{1}{\sigma}
\mathbb{E}
\left[
\left(U_{x_1,x_2}^{(\eta)}\right)^2-1
\right]
+
d^{\alpha/2}
\mathbb{E}
\left[
\frac{\omega\left(U_{x_1,x_2}^{(\eta)}\right)}{\eta(x_1)}
\right]
+
o(d^{\alpha/2}).
\]
Since \(U_{x_1,x_2}^{(\eta)}\) converges in distribution to a standard Gaussian
random variable \(U\), and since \(1/\eta(x_1)\) has mean one and is asymptotically independent of $U_{x_1,x_2}^{(\eta)}$, it follows that
\[
\mathbb{E}\left[\omega(U)\right]
=
-\frac{1}{\sigma}
\lim_{d\to0}
d^{-\alpha/2}
\mathbb{E}
\left[
\left(U_{x_1,x_2}^{(\eta)}\right)^2-1
\right]
=
-4\psi,
\]
where \(\psi\) is defined in \eqref{eq:constant_psy}.

\subsubsection{The triplewise case}

For $R\in(-1,1)$ and $u_2,u_3\in\mathbb{R}$, define
\[
\Sigma_R
=
\begin{pmatrix}
1 & R\\
R & 1
\end{pmatrix},
\qquad
Q_R(u_2,u_3)
=
\begin{pmatrix}
u_2 & u_3
\end{pmatrix}
\Sigma_R^{-1}
\begin{pmatrix}
u_2\\
u_3
\end{pmatrix}
-2.
\]
The following proposition gives the corresponding expansion for the
triplewise score contribution.

\begin{proposition}
\label{Prop_triplewise_pdf}
Let $u_2,u_3\in\mathbb{R}$ be fixed. Let
$x_1,x_2,x_3\in\mathbb{R}^{2}$ and $z_1,z_2,z_3>0$ be such that
\[
u_2
=
\frac{1}{\sigma\delta^{\alpha/2}d_{1,2}^{\alpha/2}}
\log\left(\frac{z_2}{z_1}\right),
\qquad
u_3
=
\frac{1}{\sigma\delta^{\alpha/2}d_{1,3}^{\alpha/2}}
\log\left(\frac{z_3}{z_1}\right),
\]
where
\[
\|x_2-x_1\|=\delta d_{1,2},
\qquad
\|x_3-x_1\|=\delta d_{1,3},
\qquad
\|x_3-x_2\|=\delta d_{2,3},
\]
with $d_{1,2},d_{1,3},d_{2,3}$ fixed positive numbers satisfying the strict
triangle inequalities. 
Then, as $\delta\to0$, there exists a measurable function
\[
\Omega_{\sigma}(\cdot,\cdot;d_{1,2},d_{1,3},d_{2,3})
\]
such that
\[
\frac{\partial}{\partial\sigma}
\log f_{x_1,x_2,x_3}(z_1,z_2,z_3;\sigma,\alpha)
=
\frac{1}{\sigma}Q_{R_1}(u_2,u_3)
+
\delta^{\alpha/2}
\frac{1}{z_1}
\Omega_{\sigma}(u_2,u_3;d_{1,2},d_{1,3},d_{2,3})
+
o(\delta^{\alpha/2}).
\]

Moreover, if $(U_2,U_3)$ is a centered Gaussian vector with covariance matrix
$\Sigma_{R_1}$, then
\[
\mathbb{E}
\left[
\Omega_{\sigma}(U_2,U_3;d_{1,2},d_{1,3},d_{2,3})
\right]
=
-\frac{1}{\sigma}
\lim_{\delta\to0}
\delta^{-\alpha/2}
\mathbb{E}
\left[
Q_{R_1}
\left(
U_{x_1,x_2}^{(\eta)},
U_{x_1,x_3}^{(\eta)}
\right)
\right].
\]

For the derivative with respect to $\alpha$, there exists a measurable
function
\[
B_\alpha(\cdot,\cdot;d_{1,2},d_{1,3},d_{2,3})
\]
such that
\[
\begin{aligned}
\frac{1}{\log\delta}
\frac{\partial}{\partial\alpha}
\log f_{x_1,x_2,x_3}(z_1,z_2,z_3;\sigma,\alpha)
&=
\frac{1}{2}Q_{R_1}(u_2,u_3)
+
\frac{1}{\log\delta}
B_\alpha(u_2,u_3;d_{1,2},d_{1,3},d_{2,3})
\\
&\quad+
\delta^{\alpha/2}
\frac{\sigma}{2z_1}
\Omega_{\sigma}(u_2,u_3;d_{1,2},d_{1,3},d_{2,3})
\\
&\quad+
o\!\left(\frac{1}{|\log\delta|}\right)
+
o(\delta^{\alpha/2}).
\end{aligned}
\]
In particular,
\[
\frac{1}{\log\delta}
\frac{\partial}{\partial\alpha}
\log f_{x_1,x_2,x_3}(z_1,z_2,z_3;\sigma,\alpha)
=
\frac{1}{2}Q_{R_1}(u_2,u_3)
+
o(1).
\]
\end{proposition}

The leading triplewise score contribution is therefore proportional to the
quadratic form $Q_{R_1}(u_2,u_3)$. When $u_2$ and $u_3$ are replaced by the
normalized increments of $\log\eta$ over a Delaunay triangle, Proposition
\ref{Prop_cond_dist_U_1_U_2} shows that the limiting vector is centered
Gaussian with covariance matrix $\Sigma_{R_1}$. Hence the leading
triplewise score contribution is asymptotically centered.

The function $\Omega$ contains the first-order correction to this Gaussian
approximation. Its explicit closed-form expression involves lengthy
differentiations of the trivariate Brown--Resnick density and of bivariate
Gaussian distribution functions. Since only its order and its mean are used
in the asymptotic analysis of the estimators, we do not display the full
formula.

The leading term is invariant under a relabelling of the vertices. For
example, if one uses $x_2$ as the reference point and defines
\[
\widetilde{u}_1
=
\frac{1}{\sigma\delta^{\alpha/2}d_{1,2}^{\alpha/2}}
\log\left(\frac{z_1}{z_2}\right),
\qquad
\widetilde{u}_3
=
\frac{1}{\sigma\delta^{\alpha/2}d_{2,3}^{\alpha/2}}
\log\left(\frac{z_3}{z_2}\right),
\]
then
\[
\lim_{\delta\to0}
\frac{\partial}{\partial\sigma}
\log f_{x_1,x_2,x_3}(z_1,z_2,z_3;\sigma,\alpha)
=
\frac{1}{\sigma}
Q_{R_2}(\widetilde{u}_1,\widetilde{u}_3).
\]
Thus the first-order triplewise score does not depend on the particular
choice of reference vertex, although the lexicographic ordering remains
useful for defining the sums over Delaunay triangles.

\subsection{Asymptotic properties of the MCLEs and discussion}

We now state the asymptotic properties of the one-parameter maximum
composite likelihood estimators. The key point is that, by Propositions
\ref{Prop_pairwise_pdf} and \ref{Prop_triplewise_pdf}, the leading score
terms are proportional to the squared increment statistics studied in
Theorem \ref{prop:BRtrajectories}.

\begin{theorem}
\label{Prop:Asym_Prop_CL_Est}
Assume that $\sigma_0$ belongs to the interior of a compact subset of
$(0,\infty)$ and that $\alpha_0$ belongs to the interior of a compact subset
of $(0,1)$. Let
\[
L_Z
=
\sum_{j\geq1}\sum_{k>j}
L_{Z_{k\setminus j}}(0).
\]
Then, as $N\to\infty$, the pairwise estimators satisfy
\[
\frac{\sqrt{3}}{3}
\sqrt{|E_N|}
N^{-(2-\alpha_0)/4}
\left(
\widehat{\sigma}_{2,N}^{\,2}-\sigma_0^2
\right)
\overset{\mathbb{P}}{\longrightarrow}
c_{V_2}\sigma_0^2
\left(
L_Z-\mathbb{E}[L_Z]
\right),
\]
when $\alpha_0$ is known, and
\[
\frac{\sqrt{3}}{6}
\sqrt{|E_N|}
N^{-(2-\alpha_0)/4}
\log N
\left(
\widehat{\alpha}_{2,N}-\alpha_0
\right)
\overset{\mathbb{P}}{\longrightarrow}
-
c_{V_2}
\left(
L_Z-\mathbb{E}[L_Z]
\right),
\]
when $\sigma_0$ is known.

Similarly, the triplewise estimators satisfy
\[
\frac{\sqrt{2}}{2}
\sqrt{|DT_N|}
N^{-(2-\alpha_0)/4}
\left(
\widehat{\sigma}_{3,N}^{\,2}-\sigma_0^2
\right)
\overset{\mathbb{P}}{\longrightarrow}
c_{V_3}\sigma_0^2
\left(
L_Z-\mathbb{E}[L_Z]
\right),
\]
when $\alpha_0$ is known, and
\[
\frac{\sqrt{2}}{4}
\sqrt{|DT_N|}
N^{-(2-\alpha_0)/4}
\log N
\left(
\widehat{\alpha}_{3,N}-\alpha_0
\right)
\overset{\mathbb{P}}{\longrightarrow}
-
c_{V_3}
\left(
L_Z-\mathbb{E}[L_Z]
\right),
\]
when $\sigma_0$ is known.
\end{theorem}

Several consequences follow from Theorem
\ref{Prop:Asym_Prop_CL_Est}.

First, the maximum composite likelihood estimators of $\sigma_0^2$ and
$\alpha_0$ are consistent in the fixed-domain framework considered here,
provided that the other parameter is known. Since $|E_N|$ and $|DT_N|$ are
both of order $N$, the convergence rate is of order
$N^{\alpha_0/4}$ for the scale estimators and of order
$\log(N)N^{\alpha_0/4}$ for the smoothness estimators. These rates differ
from the classical rates $N^{1/2}$ and $\log(N)N^{1/2}$ obtained for
isotropic fractional Brownian fields in regular Gaussian settings; see, for
instance, \cite{Zhu&Stein02}.

Second, the limiting fluctuations are not Gaussian. They are driven by the
centered local-time functional $L_Z-\mathbb{E}[L_Z]$, where $L_Z$ aggregates
the local times of the differences between competing spectral functions at
the interfaces of the canonical Brown--Resnick tessellation. These local
times have unknown distributions and cannot be consistently recovered from a
single realization of the max-stable field, since the spectral functions
$(Y_i)_{i\geq1}$ and the Poisson points $(U_i)_{i\geq1}$ are latent.
Consequently, the Gaussian approximation usually used in replicated
composite likelihood inference for spatial extremes, as in \cite{Padoan10},
is not justified in the present fixed-domain framework.

Third, the appearance of local times is closely related to the geometry of
the canonical tessellation. Let $W^{(1)}$ and $W^{(2)}$ be two independent
isotropic fractional Brownian fields with semi-variogram
\eqref{eq:defvar}. According to Corollary 3.4 of \cite{Chen&Xiao12}, the
Hausdorff dimension of the intersection set
\[
\{x\in\mathbb{R}^{2}:W^{(1)}(x)=W^{(2)}(x)\}
\]
is $2-\alpha/2$ for $\alpha\in(0,1)$. We conjecture that an analogous
geometric behaviour holds for the boundaries of the cells of the canonical
Brown--Resnick tessellation. In the limiting case $\alpha=2$, the associated
tessellation is a Laguerre tessellation \cite{Dombry&Kabluchko18}, whose cell
boundaries have Hausdorff dimension $1$.

Fourth, the pairwise and triplewise estimators have the same non-standard
rate and their limits are driven by the same local-time functional, up to
multiplicative constants. This suggests that suitable linear combinations of
pairwise and triplewise estimators could cancel the leading local-time term
and possibly yield an estimator with a faster, Gaussian fluctuation. Such a
construction would require a separate analysis of the next-order terms and is
left for future work.

Finally, joint estimation of $(\sigma_0^2,\alpha_0)$ is substantially more
delicate. Even for fractional Gaussian processes observed at high frequency,
joint estimation of the scale parameter and the Hurst index may lead to
singular Fisher information matrices under standard diagonal normalizations.
In \cite{Brouste&Fukasawa18}, local asymptotic normality is recovered only by
using non-diagonal, parameter-dependent rate matrices. This reflects the fact
that, under infill asymptotics, scale and smoothness effects are
asymptotically intertwined. In the present setting, the difficulty is further
amplified by the max-stable structure, the non-Gaussian local-time limits,
and the random spatial sampling scheme. A full joint asymptotic theory for
$(\sigma_0^2,\alpha_0)$ is therefore beyond the scope of this paper. The
one-parameter results obtained here should be viewed as a necessary first
step toward such a theory.

\section{Simulation study}
\label{sec:Simulation_study}

This section illustrates the finite-sample behaviour of the pairwise maximum
composite likelihood estimators introduced in
Section~\ref{sec:CLapproach} and studied theoretically in
Theorem~\ref{Prop:Asym_Prop_CL_Est}. The purpose of the simulation study is
not to provide a full numerical verification of the convergence rates, but
rather to examine whether the empirical distributions of the estimators display
the qualitative features predicted by the fixed-domain asymptotic theory. In
particular, we focus on the centering of the estimation errors and on possible
departures from Gaussianity.

\begin{figure}[htbp]
\centering
\includegraphics[width=0.75\textwidth]{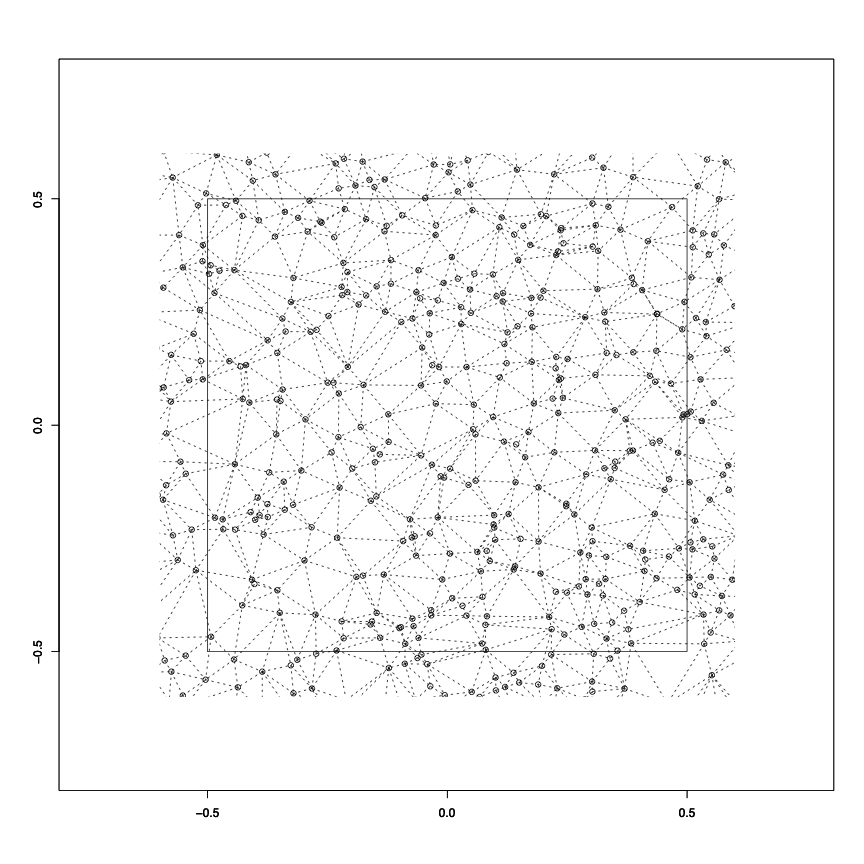}
\caption{A realization of a homogeneous Poisson point process together with its associated Delaunay triangulation. In the
pairwise composite likelihood estimator, the retained pairs are the edges of
the triangulation.}
\label{fig:delaunay}
\end{figure}

In all experiments, the observation sites are generated from a homogeneous
Poisson point process with intensity chosen so
that the expected number of points in $[-1/2,1/2]^2$ is equal to $300$. The associated
Poisson--Delaunay triangulation is then constructed. In the pairwise
composite likelihood estimator considered below, only the edges of this
triangulation are retained. A typical realization of the sampling scheme,
together with its Delaunay triangulation, is displayed in
Figure~\ref{fig:delaunay}.

For each configuration of sites, a Brown--Resnick random field is simulated
using the spectral representation \eqref{Eq_Spectral_representation}. The
spectral representation is truncated to the first $30$ terms. This truncation
level was retained after preliminary numerical checks and was found to provide
a stable approximation over the observation window. For each parameter
configuration, $5\,000$ independent realizations of the Brown--Resnick random
field are generated.

The simulation study considers the one-parameter estimation problems covered by
Theorem~\ref{Prop:Asym_Prop_CL_Est}. When estimating the scale parameter
$\sigma$, the smoothness parameter $\alpha$ is fixed at its true value. When
estimating $\alpha$, the scale parameter $\sigma$ is fixed at its true value.
Thus, the numerical experiment does not address the joint estimation of
$(\sigma,\alpha)$.

We consider several values of the smoothness parameter,
\[
\alpha_0\in\{0.5,0.75,1,1.25\},
\]
in order to cover different dependence regimes. The values
$\alpha_0=0.5$ and $\alpha_0=0.75$ fall within the theoretical range
$\alpha_0\in(0,1)$ used in Theorem~\ref{Prop:Asym_Prop_CL_Est}. The cases
$\alpha_0=1$ and $\alpha_0=1.25$ are included as exploratory numerical
experiments outside, or at the boundary of, the proven asymptotic framework.

\begin{figure}[htbp]
\centering
\includegraphics[width=\textwidth]{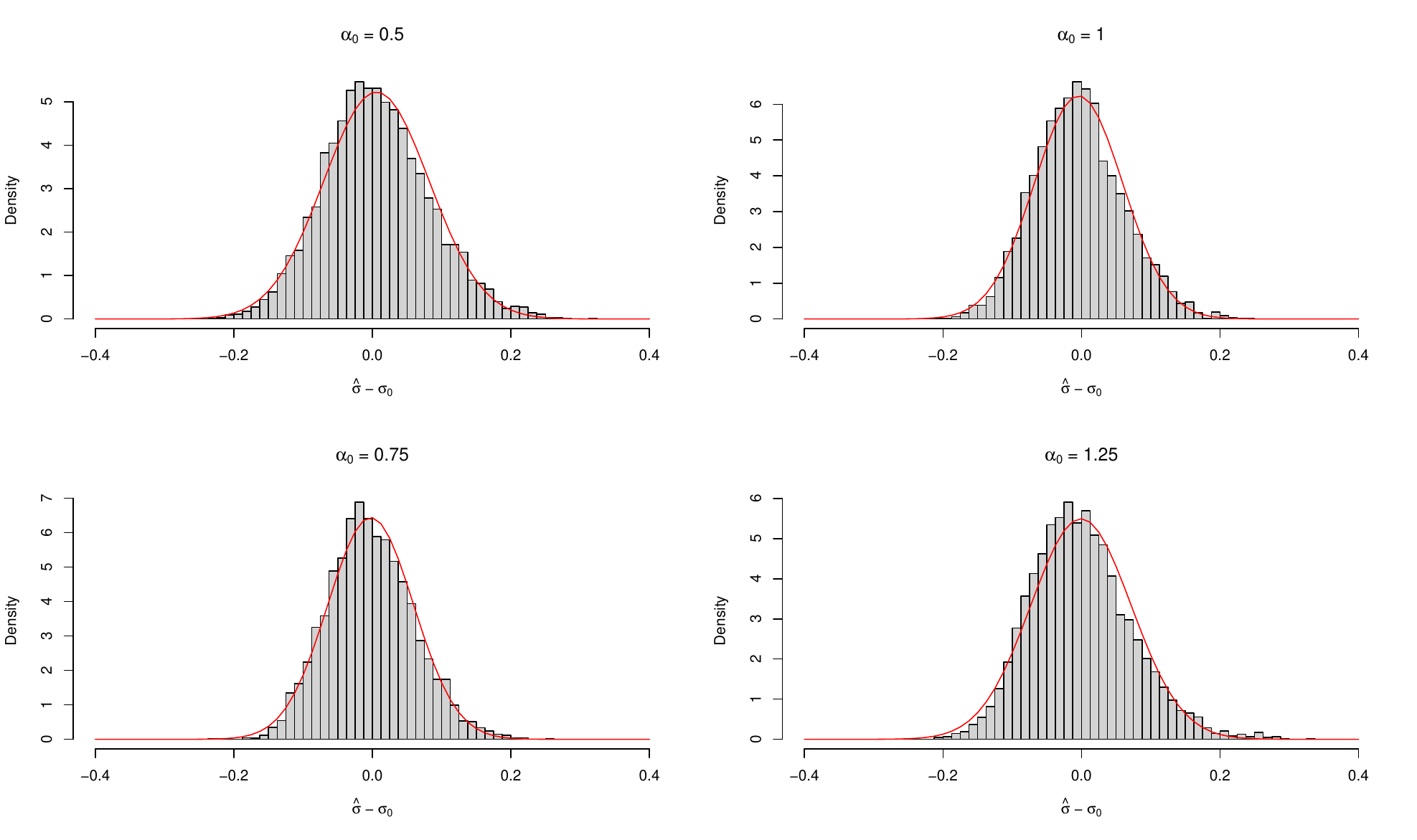}
\caption{Histograms of the estimation error
$\widehat{\sigma}_{2,N}-\sigma_0$ for different values of $\alpha_0$, based
on $5\,000$ simulations with an average of $300$ Poisson points and
truncation level $30$ in the spectral representation. The red curve
corresponds to a Gaussian density with matching empirical mean and variance.}
\label{fig:sigma}
\end{figure}

\begin{figure}[htbp]
\centering
\includegraphics[width=\textwidth]{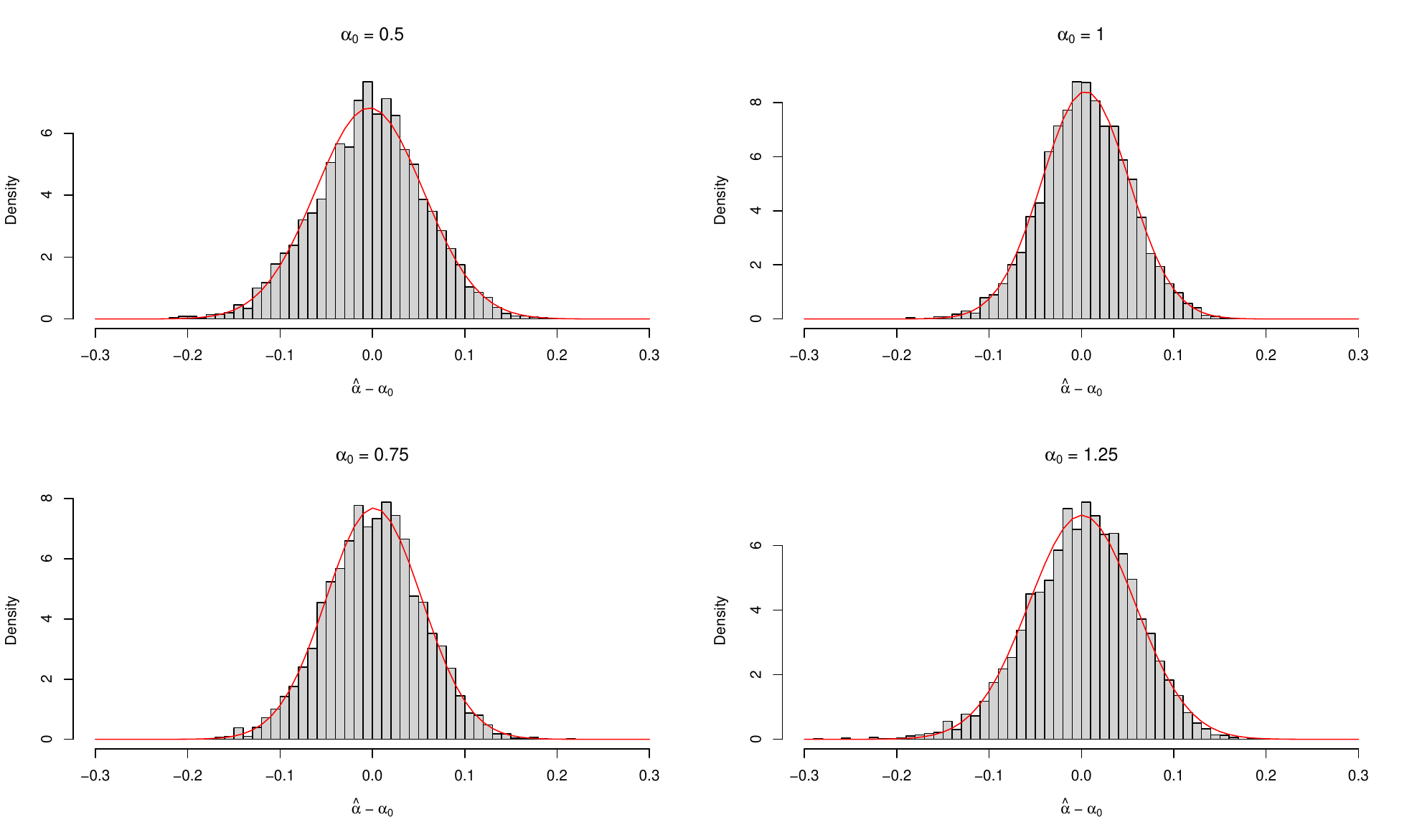}
\caption{Histograms of the estimation error
$\widehat{\alpha}_{2,N}-\alpha_0$ for different values of $\alpha_0$, based
on $5\,000$ simulations with an average of $300$ Poisson points and
truncation level $30$ in the spectral representation. The red curve
corresponds to a Gaussian density with matching empirical mean and variance.}
\label{fig:alpha}
\end{figure}

Figures~\ref{fig:sigma} and~\ref{fig:alpha} display the empirical
distributions of the estimation errors
\[
\widehat{\sigma}_{2,N}-\sigma_0
\qquad\text{and}\qquad
\widehat{\alpha}_{2,N}-\alpha_0,
\]
respectively. Although Theorem~\ref{Prop:Asym_Prop_CL_Est} is stated for
$\widehat{\sigma}_{2,N}^{\,2}-\sigma_0^2$, reporting the error on
$\sigma$ itself is more directly interpretable. The corresponding
asymptotic behaviour follows from the delta method whenever
$\sigma_0>0$.

Several features emerge from these numerical results. First, the empirical
distributions of the estimation errors are centered close to zero for all
parameter configurations considered. This is consistent with the asymptotic
centering predicted by Theorem~\ref{Prop:Asym_Prop_CL_Est} in the range
$\alpha_0\in(0,1)$. The same qualitative behaviour is also observed in the
exploratory cases $\alpha_0=1$ and $\alpha_0=1.25$.

Second, the empirical distributions are close to Gaussian but are not exactly
Gaussian. In particular, mild asymmetry can be observed in some configurations.
This feature is in line with the theoretical results, which show that the
limiting fluctuations are driven by aggregated local times associated with the
canonical Brown--Resnick tessellation, rather than by Gaussian limits. The
departure from normality remains moderate in the present finite-sample
experiments, as illustrated by the Gaussian densities with matching empirical
mean and variance superimposed in Figures~\ref{fig:sigma} and
\ref{fig:alpha}.

These simulations therefore support the qualitative message of the asymptotic
theory. In a fixed-domain framework with a single spatial realization, maximum
composite likelihood estimators may exhibit empirical distributions that are
approximately Gaussian in moderate samples, but their limiting fluctuations are
non-Gaussian. Consequently, Gaussian approximations may provide a useful rough
benchmark in finite samples, but they do not fully describe the asymptotic
uncertainty induced by the max-stable structure and by the local geometry of
the Brown--Resnick tessellation.

\section{Conclusion}

In this paper, we have studied maximum composite likelihood estimation for
spatial Brown--Resnick random fields under fixed-domain asymptotics. The
observation window is fixed, while the field is observed on an increasingly
dense and irregular set of locations generated by a homogeneous Poisson point
process. The composite likelihoods are constructed from local configurations
selected by the Poisson--Delaunay triangulation, using either Delaunay edges
or Delaunay triangles.

Our main contribution is to provide an asymptotic theory for these estimators
when only one spatial realization of the max-stable field is available. We
prove that the pairwise and triplewise maximum composite likelihood estimators
of the scale parameter and of the smoothness parameter are consistent when the
other parameter is known. Their asymptotic behaviour is, however, non-standard.
The convergence rates are slower than the classical rates obtained under
increasing-domain asymptotics or under independent temporal replications, and
the centered limiting distributions are non-Gaussian. They are driven by
aggregated local times associated with the interfaces of the canonical
Brown--Resnick tessellation.

These results show that fixed-domain inference for max-stable random fields is
fundamentally different from the usual replicated framework used in spatial
extremes. In particular, the limiting distributions involve random quantities
whose laws are unknown and which cannot be consistently recovered from a single
realization of the field, since the underlying spectral functions and Poisson
points in the max-stable representation are latent. This is not merely a
technical limitation of the method, but rather a consequence of the local
geometry of the Brown--Resnick field under infill asymptotics. It follows that
standard Gaussian uncertainty quantification may be misleading in this setting.

The simulation study illustrates the finite-sample behaviour of the pairwise
estimators. Although the theoretical limits are non-Gaussian, the empirical
distributions observed in moderate samples remain close to Gaussian in the
experiments considered, with mild departures such as asymmetry. These numerical
results should therefore be interpreted as qualitative support for the
asymptotic theory, rather than as evidence that Gaussian approximations are
asymptotically valid. In applications with a single dense spatial observation,
simulation-based procedures may be preferable for assessing estimation
uncertainty.

Several extensions remain open. First, our theoretical results are established
for isotropic fractional Brownian fields with smoothness parameter
$\alpha\in(0,1)$, equivalently Hurst index $H=\alpha/2\in(0,1/2)$. Extending
the analysis beyond this range would require different arguments, since the
local structure of increments and the associated interface contributions may
change. Second, we have considered one-parameter estimation problems, treating
the other parameter as known. A full joint asymptotic theory for
$(\sigma^2,\alpha)$ is substantially more difficult. Even for Gaussian models
under infill asymptotics, joint estimation of scale and smoothness parameters
may involve singular information matrices and non-standard normalizations. In
the Brown--Resnick setting, these difficulties are compounded by the
max-stable structure, the random sampling design, and the local-time limits
identified in this paper.

Overall, the results clarify what can be inferred from a single realization of
a Brown--Resnick random field observed over a fixed spatial domain. They also
highlight the role of the canonical max-stable tessellation in the asymptotic
behaviour of composite likelihood estimators. We expect that these findings
will be useful for further work on joint estimation, alternative local
composite likelihoods, deterministic or adaptive sampling designs, and
simulation-based inference methods for spatial extremes under fixed-domain
asymptotics.

\bibliographystyle{abbrv}
\bibliography{biblio_infill}

@article {Bevilacqua&Faouzi19,
    AUTHOR = {Bevilacqua, Moreno and Faouzi, Tarik},
     TITLE = {Estimation and prediction of {G}aussian processes using
              generalized {C}auchy covariance model under fixed domain
              asymptotics},
   JOURNAL = {Electron. J. Stat.},
  FJOURNAL = {Electronic Journal of Statistics},
    VOLUME = {13},
      YEAR = {2019},
    NUMBER = {2},
     PAGES = {3025--3048},
      ISSN = {1935-7524},
   MRCLASS = {62M09 (60G15 62F12)},
  MRNUMBER = {4010591},
MRREVIEWER = {Gilles\ Teyssi\`ere},
       DOI = {10.1214/19-ejs1597},
       URL = {https://doi.org/10.1214/19-ejs1597},
}

@article {Brouste&Fukasawa18,
    AUTHOR = {Brouste, Alexandre and Fukasawa, Masaaki},
     TITLE = {Local asymptotic normality property for fractional {G}aussian
              noise under high-frequency observations},
   JOURNAL = {Ann. Statist.},
  FJOURNAL = {The Annals of Statistics},
    VOLUME = {46},
      YEAR = {2018},
    NUMBER = {5},
     PAGES = {2045--2061},
      ISSN = {0090-5364,2168-8966},
   MRCLASS = {62F05 (60G22 62B15 62F12)},
  MRNUMBER = {3845010},
MRREVIEWER = {Danijel\ Grahovac},
       DOI = {10.1214/17-AOS1611},
       URL = {https://doi.org/10.1214/17-AOS1611},
}

@article {Chan&Wood04,
    AUTHOR = {Chan, Grace and Wood, Andrew T. A.},
     TITLE = {Estimation of fractal dimension for a class of non-{G}aussian
              stationary processes and fields},
   JOURNAL = {Ann. Statist.},
  FJOURNAL = {The Annals of Statistics},
    VOLUME = {32},
      YEAR = {2004},
    NUMBER = {3},
     PAGES = {1222--1260},
      ISSN = {0090-5364,2168-8966},
   MRCLASS = {60G10 (28A80)},
  MRNUMBER = {2065204},
MRREVIEWER = {Claude\ Tricot},
       DOI = {10.1214/009053604000000346},
       URL = {https://doi.org/10.1214/009053604000000346},
}

@article {Chen&Xiao12,
    AUTHOR = {Chen, ZhenLong and Xiao, YiMin},
     TITLE = {On intersections of independent anisotropic {G}aussian random
              fields},
   JOURNAL = {Sci. China Math.},
  FJOURNAL = {Science China. Mathematics},
    VOLUME = {55},
      YEAR = {2012},
    NUMBER = {11},
     PAGES = {2217--2232},
      ISSN = {1674-7283,1869-1862},
   MRCLASS = {60G15 (60G17 60G60)},
  MRNUMBER = {2994115},
MRREVIEWER = {Xianmin\ Geng},
       DOI = {10.1007/s11425-012-4521-9},
       URL = {https://doi.org/10.1007/s11425-012-4521-9},
}

@article {Chenavier22,
    AUTHOR = {Chenavier, Nicolas and Henze, Norbert and Otto, Moritz},
     TITLE = {Limit laws for large {$k$}th-nearest neighbor balls},
   JOURNAL = {J. Appl. Probab.},
  FJOURNAL = {Journal of Applied Probability},
    VOLUME = {59},
      YEAR = {2022},
    NUMBER = {3},
     PAGES = {880--894},
      ISSN = {0021-9002,1475-6072},
   MRCLASS = {60F05 (60D05 60G70)},
  MRNUMBER = {4480085},
       DOI = {10.1017/jpr.2021.92},
       URL = {https://doi.org/10.1017/jpr.2021.92},
}

@article {Chenavier&Robert25a,
    AUTHOR = {Chenavier, Nicolas and Robert, Christian Y.},
     TITLE = {Central limit theorems for squared increment sums of fractional {B}rownian fields based on a {D}elaunay triangulation in $2D$},
   JOURNAL = {WP},
  FJOURNAL = {Working paper},
      YEAR = {2026},
}

@article {Chenavier&Robert25b,
    AUTHOR = {Chenavier, Nicolas and Robert, Christian Y.},
     TITLE = {Limit theorems for squared increment sums of the maximum of two isotropic fractional {B}rownian
fields over a fixed-domain},
   JOURNAL = {WP},
  FJOURNAL = {Working paper},
      YEAR = {2026},
}

@book {Cohen&Istas13,
    AUTHOR = {Cohen, Serge and Istas, Jacques},
     TITLE = {Fractional fields and applications},
    SERIES = {Math\'ematiques \& Applications (Berlin) [Mathematics \&
              Applications]},
    VOLUME = {73},
      NOTE = {With a foreword by St\'ephane Jaffard},
 PUBLISHER = {Springer, Heidelberg},
      YEAR = {2013},
     PAGES = {xii+270},
      ISBN = {978-3-642-36738-0; 978-3-642-36739-7},
   MRCLASS = {60G18 (60G22 60G60 62M40 65C50)},
  MRNUMBER = {3088856},
MRREVIEWER = {Ingemar\ Kaj},
       DOI = {10.1007/978-3-642-36739-7},
       URL = {https://doi.org/10.1007/978-3-642-36739-7},
}

@article {Davison12,
    AUTHOR = {Davison, A. C. and Padoan, S. A. and Ribatet, M.},
     TITLE = {Rejoinder: ``{S}tatistical modeling of spatial extremes''
              [MR2963988; MR2963989; MR2963990; MR2963991; MR2963980]},
   JOURNAL = {Statist. Sci.},
  FJOURNAL = {Statistical Science. A Review Journal of the Institute of
              Mathematical Statistics},
    VOLUME = {27},
      YEAR = {2012},
    NUMBER = {2},
     PAGES = {199--201},
      ISSN = {0883-4237,2168-8745},
   MRCLASS = {62M30 (60G70 62G32 62H05 86A32)},
  MRNUMBER = {2963992},
       DOI = {10.1214/12-STS376REJ},
       URL = {https://doi.org/10.1214/12-STS376REJ},
}

@article {Dombry18,
    AUTHOR = {Dombry, Cl\'ement and Engelke, Sebastian and Oesting, Marco},
     TITLE = {Asymptotic properties of the maximum likelihood estimator for multivariate extremes},
   JOURNAL = {Available at arXiv:1612.05178},
      YEAR = {2018},
}

@article {Dombry&Eyi-Minko13,
    AUTHOR = {Dombry, Cl\'ement and Eyi-Minko, Fr\'ed\'eric},
     TITLE = {Regular conditional distributions of continuous max-infinitely
              divisible random fields},
   JOURNAL = {Electron. J. Probab.},
  FJOURNAL = {Electronic Journal of Probability},
    VOLUME = {18},
      YEAR = {2013},
     PAGES = {no. 7, 21},
      ISSN = {1083-6489},
   MRCLASS = {60G70 (60G25)},
  MRNUMBER = {3024101},
MRREVIEWER = {Ou\ Zhao},
       DOI = {10.1214/EJP.v18-1991},
       URL = {https://doi.org/10.1214/EJP.v18-1991},
}

@article {Dombry&Kabluchko18,
    AUTHOR = {Dombry, Cl\'ement and Kabluchko, Zakhar},
     TITLE = {Random tessellations associated with max-stable random fields},
   JOURNAL = {Bernoulli},
  FJOURNAL = {Bernoulli. Official Journal of the Bernoulli Society for
              Mathematical Statistics and Probability},
    VOLUME = {24},
      YEAR = {2018},
    NUMBER = {1},
     PAGES = {30--52},
      ISSN = {1350-7265,1573-9759},
   MRCLASS = {60G60 (60D05 60G52 60G55 60G70)},
  MRNUMBER = {3706749},
MRREVIEWER = {Ilya\ S.\ Molchanov},
       DOI = {10.3150/16-BEJ817},
       URL = {https://doi.org/10.3150/16-BEJ817},
}

@article {Einmahl15,
    AUTHOR = {Einmahl, John H. J. and Kiriliouk, Anna and Krajina, Andrea
              and Segers, Johan},
     TITLE = {An {$M$}-estimator of spatial tail dependence},
   JOURNAL = {J. R. Stat. Soc. Ser. B. Stat. Methodol.},
  FJOURNAL = {Journal of the Royal Statistical Society. Series B.
              Statistical Methodology},
    VOLUME = {78},
      YEAR = {2016},
    NUMBER = {1},
     PAGES = {275--298},
      ISSN = {1369-7412,1467-9868},
   MRCLASS = {62G32 (60F05 60F17 60G70 62G05 62G20 62H11)},
  MRNUMBER = {3453656},
       DOI = {10.1111/rssb.12114},
       URL = {https://doi.org/10.1111/rssb.12114},
}

@article {Engelke14,
    AUTHOR = {Engelke, Sebastian and Malinowski, Alexander and Kabluchko,
              Zakhar and Schlather, Martin},
     TITLE = {Estimation of {H}\"usler-{R}eiss distributions and
              {B}rown-{R}esnick processes},
   JOURNAL = {J. R. Stat. Soc. Ser. B. Stat. Methodol.},
  FJOURNAL = {Journal of the Royal Statistical Society. Series B.
              Statistical Methodology},
    VOLUME = {77},
      YEAR = {2015},
    NUMBER = {1},
     PAGES = {239--265},
      ISSN = {1369-7412,1467-9868},
   MRCLASS = {60G70 (60G52 60G55 62-04 62G32 62M15 62P12)},
  MRNUMBER = {3299407},
       DOI = {10.1111/rssb.12074},
       URL = {https://doi.org/10.1111/rssb.12074},
}

@article {Fondeville&Davison18,
    AUTHOR = {de Fondeville, R. and Davison, A. C.},
     TITLE = {High-dimensional peaks-over-threshold inference},
   JOURNAL = {Biometrika},
  FJOURNAL = {Biometrika},
    VOLUME = {105},
      YEAR = {2018},
    NUMBER = {3},
     PAGES = {575--592},
      ISSN = {0006-3444,1464-3510},
   MRCLASS = {62G32 (60G52 60G70 62P12)},
  MRNUMBER = {3842886},
       DOI = {10.1093/biomet/asy026},
       URL = {https://doi.org/10.1093/biomet/asy026},
}

@article {Geman&Horowitz80,
    AUTHOR = {Geman, Donald and Horowitz, Joseph},
     TITLE = {Occupation densities},
   JOURNAL = {Ann. Probab.},
  FJOURNAL = {The Annals of Probability},
    VOLUME = {8},
      YEAR = {1980},
    NUMBER = {1},
     PAGES = {1--67},
      ISSN = {0091-1798,2168-894X},
   MRCLASS = {60J55 (26A27 60G15 60G17)},
  MRNUMBER = {556414},
MRREVIEWER = {Simeon\ M.\ Berman},
       URL =
              {http://links.jstor.org/sici?sici=0091-1798(198002)8:1<1:OD>2.0.CO;2-M&origin=MSN},
}

@book {Haan&Ferreira06,
    AUTHOR = {de Haan, Laurens and Ferreira, Ana},
     TITLE = {Extreme value theory},
    SERIES = {Springer Series in Operations Research and Financial
              Engineering},
      NOTE = {An introduction},
 PUBLISHER = {Springer, New York},
      YEAR = {2006},
     PAGES = {xviii+417},
      ISBN = {978-0-387-23946-0; 0-387-23946-4},
   MRCLASS = {62-02 (60G70 62G32)},
  MRNUMBER = {2234156},
MRREVIEWER = {L\'aszl\'o\ Viharos},
       DOI = {10.1007/0-387-34471-3},
       URL = {https://doi.org/10.1007/0-387-34471-3},
}

@article {Henze82,
    AUTHOR = {Henze, Norbert},
     TITLE = {The limit distribution for maxima of ``weighted''\
              {$r$}th-nearest-neighbour distances},
   JOURNAL = {J. Appl. Probab.},
  FJOURNAL = {Journal of Applied Probability},
    VOLUME = {19},
      YEAR = {1982},
    NUMBER = {2},
     PAGES = {344--354},
      ISSN = {0021-9002,1475-6072},
   MRCLASS = {62E20 (60F05 62G10)},
  MRNUMBER = {649972},
MRREVIEWER = {W.\ Hoeffding},
       DOI = {10.2307/3213486},
       URL = {https://doi.org/10.2307/3213486},
}

@article {Huser&Davison13,
    AUTHOR = {Huser, R. and Davison, A. C.},
     TITLE = {Composite likelihood estimation for the {B}rown-{R}esnick
              process},
   JOURNAL = {Biometrika},
  FJOURNAL = {Biometrika},
    VOLUME = {100},
      YEAR = {2013},
    NUMBER = {2},
     PAGES = {511--518},
      ISSN = {0006-3444,1464-3510},
   MRCLASS = {62F10 (60G52 62G32 62M30)},
  MRNUMBER = {3068451},
       DOI = {10.1093/biomet/ass089},
       URL = {https://doi.org/10.1093/biomet/ass089},
}

@article {Jaramillo21,
    AUTHOR = {Jaramillo, Arturo and Nourdin, Ivan and Peccati, Giovanni},
     TITLE = {Approximation of fractional local times: zero energy and
              derivatives},
   JOURNAL = {Ann. Appl. Probab.},
  FJOURNAL = {The Annals of Applied Probability},
    VOLUME = {31},
      YEAR = {2021},
    NUMBER = {5},
     PAGES = {2143--2191},
      ISSN = {1050-5164,2168-8737},
   MRCLASS = {60G22 (60F17 60H07 60J55)},
  MRNUMBER = {4332693},
       DOI = {10.1214/20-aap1643},
       URL = {https://doi.org/10.1214/20-aap1643},
}

@article {Kabluchko09,
    AUTHOR = {Kabluchko, Zakhar and Schlather, Martin and de Haan, Laurens},
     TITLE = {Stationary max-stable fields associated to negative definite
              functions},
   JOURNAL = {Ann. Probab.},
  FJOURNAL = {The Annals of Probability},
    VOLUME = {37},
      YEAR = {2009},
    NUMBER = {5},
     PAGES = {2042--2065},
      ISSN = {0091-1798,2168-894X},
   MRCLASS = {60G70 (60G15 60G22 60G55 60G60)},
  MRNUMBER = {2561440},
       DOI = {10.1214/09-AOP455},
       URL = {https://doi.org/10.1214/09-AOP455},
}

@article {Padoan10,
    AUTHOR = {Padoan, S. A. and Ribatet, M. and Sisson, S. A.},
     TITLE = {Likelihood-based inference for max-stable processes},
   JOURNAL = {J. Amer. Statist. Assoc.},
  FJOURNAL = {Journal of the American Statistical Association},
    VOLUME = {105},
      YEAR = {2010},
    NUMBER = {489},
     PAGES = {263--277},
      ISSN = {0162-1459,1537-274X},
   MRCLASS = {99-01},
  MRNUMBER = {2757202},
       DOI = {10.1198/jasa.2009.tm08577},
       URL = {https://doi.org/10.1198/jasa.2009.tm08577},
}

@article {Robert20,
    AUTHOR = {Robert, Christian Y.},
     TITLE = {Power variations for a class of {B}rown-{R}esnick processes},
   JOURNAL = {Extremes},
  FJOURNAL = {Extremes. Statistical Theory and Applications in Science,
              Engineering and Economics},
    VOLUME = {23},
      YEAR = {2020},
    NUMBER = {2},
     PAGES = {215--244},
      ISSN = {1386-1999,1572-915X},
   MRCLASS = {60G70 (60F05 60G17)},
  MRNUMBER = {4102383},
MRREVIEWER = {Peter\ Kern},
       DOI = {10.1007/s10687-020-00373-4},
       URL = {https://doi.org/10.1007/s10687-020-00373-4},
}

@book {Schneider&Weil08,
    AUTHOR = {Schneider, Rolf and Weil, Wolfgang},
     TITLE = {Stochastic and integral geometry},
    SERIES = {Probability and its Applications (New York)},
 PUBLISHER = {Springer-Verlag, Berlin},
      YEAR = {2008},
     PAGES = {xii+693},
      ISBN = {978-3-540-78858-4},
   MRCLASS = {60-02 (52A22 60D05 60G55 62M30)},
  MRNUMBER = {2455326},
MRREVIEWER = {V.\ K.\ Ohanyan},
       DOI = {10.1007/978-3-540-78859-1},
       URL = {https://doi.org/10.1007/978-3-540-78859-1},
}

@article {Varin11,
    AUTHOR = {Varin, Cristiano and Reid, Nancy and Firth, David},
     TITLE = {An overview of composite likelihood methods},
   JOURNAL = {Statist. Sinica},
  FJOURNAL = {Statistica Sinica},
    VOLUME = {21},
      YEAR = {2011},
    NUMBER = {1},
     PAGES = {5--42},
      ISSN = {1017-0405,1996-8507},
   MRCLASS = {62-02 (62F10 62M40)},
  MRNUMBER = {2796852},
}

@article {Zhu&Stein02,
    AUTHOR = {Zhu, Zhengyuan and Stein, Michael L.},
     TITLE = {Parameter estimation for fractional {B}rownian surfaces},
   JOURNAL = {Statist. Sinica},
  FJOURNAL = {Statistica Sinica},
    VOLUME = {12},
      YEAR = {2002},
    NUMBER = {3},
     PAGES = {863--883},
      ISSN = {1017-0405,1996-8507},
   MRCLASS = {28A80 (60G60)},
  MRNUMBER = {1929968},
}

@article{Anderes10,
  author = {Anderes, E.},
  title = {On the consistent separation of scale and variance for Gaussian random fields},
  journal = {The Annals of Statistics},
  volume = {38},
  pages = {870--893},
  year = {2010}
}

@article{Bachoc19,
  author = {Bachoc, F. and Bevilacqua, M. and Velandia, D.},
  title = {Composite likelihood estimation for a Gaussian process under fixed domain asymptotics},
  journal = {Journal of Multivariate Analysis},
  volume = {174},
  pages = {104534},
  year = {2019}
}

@article{Bachoc&Lagnoux20,
  author = {Bachoc, F. and Lagnoux, A.},
  title = {Fixed-domain asymptotic properties of maximum composite likelihood estimators for Gaussian processes},
  journal = {Journal of Statistical Planning and Inference},
  volume = {209},
  pages = {62--75},
  year = {2020}
}

@book{Cressie93,
  author = {Cressie, N.},
  title = {Statistics for Spatial Data},
  publisher = {J. Wiley},
  year = {1993}
}

@article{Gine90,
  author = {Gin\'{e}, E. and Hahn, M.G. and Vatan, P.},
  title = {Max-infinitely divisible and max-stable sample continuous processes},
  journal = {Probability Theory and Related Fields},
  volume = {87},
  pages = {139--165},
  year = {1990}
}

@article{Haan84,
  author = {de Haan, L.},
  title = {A spectral representation for max-stable processes},
  journal = {The Annals of Probability},
  volume = {12},
  pages = {1194--1204},
  year = {1984}
}

@book{Ibragimov&Rozanov78,
  author = {Ibragimov, I. A. and Rozanov, Y. A.},
  title = {Gaussian Random Processes},
  publisher = {Springer-Verlag New York},
  year = {1978}
}

@article{Kaufman&Shaby13,
  author = {Kaufman, C. G. and Shaby, B. A.},
  title = {The role of the range parameter for estimation and prediction in geostatistics},
  journal = {Biometrika},
  volume = {100},
  pages = {473--484},
  year = {2013}
}

@article{Li13,
  author = {Li, C.},
  title = {Maximum-likelihood estimation for diffusion processes via closed-form density expansions},
  journal = {The Annals of Statistics},
  volume = {41},
  pages = {1350--1380},
  year = {2013}
}

@book{Matern60,
  author = {Mat\'ern, B.},
  title = {Spatial variation},
  series = {Meddelanden fran Statens Skogsforskningsinstitut},
  volume = {49},
  pages = {5},
  note = {Second ed. (1986), Lecture Notes in Statistics 36, New York: Springer},
  year = {1960}
}

@book{Stein99,
  author = {Stein, M. L.},
  title = {Interpolation of Spatial Data},
  series = {Springer Series in Statistics},
  publisher = {Springer-Verlag New York},
  year = {1999}
}

@article{Vaart96,
  author = {van der Vaart, A.},
  title = {Maximum likelihood estimation under a spatial sampling scheme},
  journal = {The Annals of Statistics},
  volume = {5},
  pages = {2049--2057},
  year = {1996}
}

@article{Ying91,
  author = {Ying, Z.},
  title = {Asymptotic properties of a maximum likelihood estimator with data from a Gaussian process},
  journal = {Journal of Multivariate Analysis},
  volume = {36},
  pages = {280--296},
  year = {1991}
}

@article{Ying93,
  author = {Ying, Z.},
  title = {Maximum likelihood estimation of parameters under a spatial sampling scheme},
  journal = {The Annals of Statistics},
  volume = {21},
  pages = {1567--1590},
  year = {1993}
}

@article{Zhang04,
  author = {Zhang, H.},
  title = {Inconsistent estimation and asymptotically equivalent interpolations in model-based geostatistics},
  journal = {Journal of the American Statistical Association},
  volume = {99},
  pages = {250--261},
  year = {2004}
}

\newpage

\appendix

\noindent
{\LARGE\bf Supplementary material}

\section{Proofs}
\label{sec:proofs}

\subsection{Proof of Proposition \ref{Prop_cond_dist_U}}

Let
\[
d=\|x_2-x_1\|>0,
\qquad
a=\sigma d^{\alpha/2},
\qquad
u=\frac{\log(z_2/z_1)}{a},
\qquad
v(u)=u+\frac{a}{2}.
\]
For a fixed value $z_1>0$ and for $u\in\mathbb{R}$, set
\[
z_2=z_1 e^{a u}.
\]
Then
\[
\left\{
U_{x_1,x_2}^{(\eta)}\leq u
\right\}
=
\left\{
\eta(x_2)\leq z_2
\right\}
\quad\text{on the event }\{\eta(x_1)=z_1\}.
\]
Therefore,
\[
\mathbb{P}
\left\{
U_{x_1,x_2}^{(\eta)}\leq u
\,\middle|\,
\eta(x_1)=z_1
\right\}
=
\mathbb{P}
\left\{
\eta(x_2)\leq z_2
\,\middle|\,
\eta(x_1)=z_1
\right\}.
\]
By stationarity of the Brown--Resnick field,
\[
\mathbb{P}
\left\{
\eta(x_2)\leq z_2
\,\middle|\,
\eta(x_1)=z_1
\right\}
=
\mathbb{P}
\left\{
\eta(x_2-x_1)\leq z_2
\,\middle|\,
\eta(0)=z_1
\right\}.
\]
Let $h=x_2-x_1$. We now apply the regular conditional distribution formula
for max-stable random fields, see Proposition 4.2 in
\cite{Dombry&Eyi-Minko13}. Since $Y(0)=1$ almost surely, we obtain
\begin{multline*}
\mathbb{P}
\left\{
\eta(h)\leq z_2
\,\middle|\,
\eta(0)=z_1
\right\}
\\
=
\exp\left(
-
\mathbb{E}
\left[
\left(
\frac{Y(h)}{z_2}
-
\frac{Y(0)}{z_1}
\right)_+
\right]
\right)
\mathbb{E}
\left[
\mathbb{I}
\left\{
\frac{Y(h)}{z_2}
\leq
\frac{Y(0)}{z_1}
\right\}
Y(0)
\right].
\end{multline*}

We compute the two terms on the right-hand side. First, since $Y(0)=1$ a.s.,
\begin{align*}
\mathbb{E}
\left[
\mathbb{I}
\left\{
\frac{Y(h)}{z_2}
\leq
\frac{Y(0)}{z_1}
\right\}
Y(0)
\right]
&=
\mathbb{P}
\left\{
\frac{Y(h)}{z_2}
\leq
\frac{1}{z_1}
\right\} \\
&=
\mathbb{P}
\left\{
\exp\{W(h)-\gamma(h)\}
\leq
\frac{z_2}{z_1}
\right\}.
\end{align*}
Since
\[
W(h)\sim \mathcal{N}(0,2\gamma(h)),
\qquad
2\gamma(h)=\sigma^2\|h\|^\alpha=a^2,
\]
we get
\begin{align*}
\mathbb{P}
\left\{
\exp\{W(h)-\gamma(h)\}
\leq
\frac{z_2}{z_1}
\right\}
&=
\mathbb{P}
\left\{
\frac{W(h)}{a}
\leq
\frac{1}{a}\log\left(\frac{z_2}{z_1}\right)
+
\frac{a}{2}
\right\} \\
&=
\Phi\left(u+\frac{a}{2}\right)
=
\Phi(v(u)).
\end{align*}

Second,
\begin{align*}
\mathbb{E}
\left[
\left(
\frac{Y(h)}{z_2}
-
\frac{Y(0)}{z_1}
\right)_+
\right]
&=
\mathbb{E}
\left[
\max\left\{
\frac{Y(h)}{z_2},
\frac{Y(0)}{z_1}
\right\}
-
\frac{Y(0)}{z_1}
\right] \\
&=
V_{0,h}(z_1,z_2)-\frac{1}{z_1}.
\end{align*}
By stationarity of the Brown--Resnick field,
\[
V_{0,h}(z_1,z_2)=V_{x_1,x_2}(z_1,z_2),
\]
and, by homogeneity of the exponent function,
\[
V_{x_1,x_2}(z_1,z_2)
=
\frac{1}{z_1}
V_{x_1,x_2}
\left(
1,
\frac{z_2}{z_1}
\right)
=
\frac{1}{z_1}
V_{x_1,x_2}
\left(
1,
e^{a u}
\right).
\]
Hence
\[
\mathbb{E}
\left[
\left(
\frac{Y(h)}{z_2}
-
\frac{Y(0)}{z_1}
\right)_+
\right]
=
\frac{1}{z_1}
\left[
V_{x_1,x_2}
\left(
1,
e^{a u}
\right)
-1
\right].
\]
Combining the two preceding identities gives
\[
\mathbb{P}
\left\{
U_{x_1,x_2}^{(\eta)}\leq u
\,\middle|\,
\eta(x_1)=z_1
\right\}
=
\exp\left(
-
\frac{1}{z_1}
\left[
V_{x_1,x_2}
\left(
1,
e^{a u}
\right)
-1
\right]
\right)
\Phi(v(u)),
\]
which proves the conditional distribution formula.

We now derive the marginal distribution. Since $\eta(x_1)$ has a standard
unit Fr\'echet distribution, the random variable $1/\eta(x_1)$ has a standard
exponential distribution. Therefore, for
\[
A=
V_{x_1,x_2}
\left(
1,
e^{a u}
\right)
-1,
\]
we have
\[
\mathbb{E}
\left[
\exp\left\{
-\frac{A}{\eta(x_1)}
\right\}
\right]
=
\frac{1}{1+A}
=
\frac{1}{
V_{x_1,x_2}
\left(
1,
e^{a u}
\right)
}.
\]
Using the conditional formula and integrating with respect to $\eta(x_1)$,
we obtain
\[
\mathbb{P}
\left\{
U_{x_1,x_2}^{(\eta)}\leq u
\right\}
=
\frac{\Phi(v(u))}
{
V_{x_1,x_2}
\left(
1,
e^{a u}
\right)
}.
\]

It remains to prove the limiting distribution. Since $a=\sigma d^{\alpha/2}$
tends to $0$ as $d\to0$, we have
\[
\Phi(v(u))
=
\Phi\left(u+\frac{a}{2}\right)
\longrightarrow
\Phi(u).
\]
Moreover, using the explicit expression of the pairwise exponent function,
\[
V_{x_1,x_2}
\left(
1,
e^{a u}
\right)
=
\Phi\left(u+\frac{a}{2}\right)
+
e^{-a u}
\Phi\left(-u+\frac{a}{2}\right),
\]
and therefore
\[
V_{x_1,x_2}
\left(
1,
e^{a u}
\right)
\longrightarrow
\Phi(u)+\Phi(-u)
=
1.
\]
Consequently,
\[
\lim_{d\to0}
\mathbb{P}
\left\{
U_{x_1,x_2}^{(\eta)}\leq u
\right\}
=
\Phi(u),
\qquad u\in\mathbb{R}.
\]
This concludes the proof.

\subsection{Proof of Proposition \ref{Prop_cond_dist_U_1_U_2}}

We follow the same strategy as in the proof of Proposition
\ref{Prop_cond_dist_U}. Let
\[
h_2=x_2-x_1,
\qquad
h_3=x_3-x_1,
\]
and write
\[
d_{1,2}=\|h_2\|,
\qquad
d_{1,3}=\|h_3\|,
\qquad
d_{2,3}=\|h_3-h_2\|.
\]
For fixed $z_1>0$ and $u_2,u_3\in\mathbb{R}$, set
\[
z_2
=
z_1\exp\left\{\sigma d_{1,2}^{\alpha/2}u_2\right\},
\qquad
z_3
=
z_1\exp\left\{\sigma d_{1,3}^{\alpha/2}u_3\right\}.
\]
Then, on the event $\{\eta(x_1)=z_1\}$,
\[
\left\{
U_{x_1,x_2}^{(\eta)}\leq u_2,
U_{x_1,x_3}^{(\eta)}\leq u_3
\right\}
=
\left\{
\eta(x_2)\leq z_2,
\eta(x_3)\leq z_3
\right\}.
\]
Hence
\begin{equation*}
\mathbb{P}
\left\{
U_{x_1,x_2}^{(\eta)}\leq u_2,
U_{x_1,x_3}^{(\eta)}\leq u_3
\,\middle|\,
\eta(x_1)=z_1
\right\}
=
\mathbb{P}
\left\{
\eta(x_2)\leq z_2,
\eta(x_3)\leq z_3
\,\middle|\,
\eta(x_1)=z_1
\right\}.
\end{equation*}
By stationarity of the Brown--Resnick field, this conditional probability is
equal to
\[
\mathbb{P}
\left\{
\eta(h_2)\leq z_2,
\eta(h_3)\leq z_3
\,\middle|\,
\eta(0)=z_1
\right\}.
\]

We now use the regular conditional distribution formula for max-stable random
fields, see Proposition 4.2 in \cite{Dombry&Eyi-Minko13}. Since $Y(0)=1$
almost surely, we obtain
\begin{equation}
\mathbb{P}
\left\{
\eta(h_2)\leq z_2,
\eta(h_3)\leq z_3
\,\middle|\,
\eta(0)=z_1
\right\}
=
\exp\left(
-
\mathbb{E}
\left[
\left(
\max_{i=2,3}\frac{Y(h_i)}{z_i}
-
\frac{1}{z_1}
\right)_+
\right]
\right)
\mathbb{P}
\left\{
\max_{i=2,3}
\frac{Y(h_i)}{z_i}
\leq
\frac{1}{z_1}
\right\}.
\label{eq:cond_triple_DEM}
\end{equation}

We first compute the probability term in \eqref{eq:cond_triple_DEM}. Since
\[
Y(h_i)=\exp\{W(h_i)-\gamma(h_i)\},
\qquad i=2,3,
\]
the event
\[
\left\{
\max_{i=2,3}\frac{Y(h_i)}{z_i}
\leq
\frac{1}{z_1}
\right\}
\]
is equivalent to
\[
\left\{
W(h_2)
\leq
\log\left(\frac{z_2}{z_1}\right)+\gamma(h_2),
\quad
W(h_3)
\leq
\log\left(\frac{z_3}{z_1}\right)+\gamma(h_3)
\right\}.
\]
Using
\[
2\gamma(h_2)=\sigma^2d_{1,2}^{\alpha},
\qquad
2\gamma(h_3)=\sigma^2d_{1,3}^{\alpha},
\]
we get
\[
\frac{\log(z_2/z_1)+\gamma(h_2)}
{\sqrt{2\gamma(h_2)}}
=
u_2+\frac{1}{2}\sigma d_{1,2}^{\alpha/2}
=
v_{1,2}(u_2),
\]
and
\[
\frac{\log(z_3/z_1)+\gamma(h_3)}
{\sqrt{2\gamma(h_3)}}
=
u_3+\frac{1}{2}\sigma d_{1,3}^{\alpha/2}
=
v_{1,3}(u_3).
\]
Moreover, the Gaussian vector
\[
\left(
\frac{W(h_2)}{\sqrt{2\gamma(h_2)}},
\frac{W(h_3)}{\sqrt{2\gamma(h_3)}}
\right)
\]
is centered with covariance matrix
\[
\Sigma_{R_1}
=
\begin{pmatrix}
1 & R_1\\
R_1 & 1
\end{pmatrix},
\]
where
\[
R_1
=
\frac{
d_{1,2}^{\alpha}+d_{1,3}^{\alpha}-d_{2,3}^{\alpha}
}
{
2(d_{1,2}d_{1,3})^{\alpha/2}
}.
\]
Therefore
\[
\mathbb{P}
\left\{
\max_{i=2,3}
\frac{Y(h_i)}{z_i}
\leq
\frac{1}{z_1}
\right\}
=
\Phi_2
\left(
\begin{pmatrix}
v_{1,2}(u_2)\\
v_{1,3}(u_3)
\end{pmatrix};
\Sigma_{R_1}
\right).
\]

We now compute the exponential term in \eqref{eq:cond_triple_DEM}. By the
definition of the exponent function,
\begin{align*}
\mathbb{E}
\left[
\left(
\max_{i=2,3}\frac{Y(h_i)}{z_i}
-
\frac{1}{z_1}
\right)_+
\right]
&=
\mathbb{E}
\left[
\max\left\{
\frac{1}{z_1},
\frac{Y(h_2)}{z_2},
\frac{Y(h_3)}{z_3}
\right\}
-
\frac{1}{z_1}
\right] \\
&=
V_{0,h_2,h_3}(z_1,z_2,z_3)-\frac{1}{z_1}.
\end{align*}
By stationarity of the Brown--Resnick field,
\[
V_{0,h_2,h_3}(z_1,z_2,z_3)
=
V_{x_1,x_2,x_3}(z_1,z_2,z_3).
\]
Using the homogeneity of the exponent function, we obtain
\[
V_{x_1,x_2,x_3}(z_1,z_2,z_3)
=
\frac{1}{z_1}
V_{x_1,x_2,x_3}
\left(
1,
\frac{z_2}{z_1},
\frac{z_3}{z_1}
\right).
\]
Thus
\[
\mathbb{E}
\left[
\left(
\max_{i=2,3}\frac{Y(h_i)}{z_i}
-
\frac{1}{z_1}
\right)_+
\right]
=
\frac{1}{z_1}
\left[
V_{x_1,x_2,x_3}
\left(
1,
e^{\sigma d_{1,2}^{\alpha/2}u_2},
e^{\sigma d_{1,3}^{\alpha/2}u_3}
\right)
-1
\right].
\]
Combining the two preceding computations yields
\begin{multline*}
\mathbb{P}
\left\{
U_{x_1,x_2}^{(\eta)}\leq u_2,
U_{x_1,x_3}^{(\eta)}\leq u_3
\,\middle|\,
\eta(x_1)=z_1
\right\}
\\
=
\exp\left(
-\frac{1}{z_1}
\left[
V_{x_1,x_2,x_3}
\left(
1,
e^{\sigma d_{1,2}^{\alpha/2}u_2},
e^{\sigma d_{1,3}^{\alpha/2}u_3}
\right)
-1
\right]
\right)
\\
\times
\Phi_2
\left(
\begin{pmatrix}
v_{1,2}(u_2)\\
v_{1,3}(u_3)
\end{pmatrix};
\Sigma_{R_1}
\right),
\end{multline*}
which proves the conditional distribution formula.

We next derive the marginal distribution. Since $\eta(x_1)$ has a standard
unit Fr\'echet distribution, the random variable $1/\eta(x_1)$ has a standard
exponential distribution. Therefore, with
\[
A
=
V_{x_1,x_2,x_3}
\left(
1,
e^{\sigma d_{1,2}^{\alpha/2}u_2},
e^{\sigma d_{1,3}^{\alpha/2}u_3}
\right)
-1,
\]
we have
\[
\mathbb{E}
\left[
\exp\left\{-\frac{A}{\eta(x_1)}\right\}
\right]
=
\frac{1}{1+A}.
\]
Integrating the conditional distribution with respect to $\eta(x_1)$ gives
\begin{equation*}
\mathbb{P}
\left\{
U_{x_1,x_2}^{(\eta)}\leq u_2,
U_{x_1,x_3}^{(\eta)}\leq u_3
\right\}
=
\frac{
\Phi_2
\left(
\begin{pmatrix}
v_{1,2}(u_2)\\
v_{1,3}(u_3)
\end{pmatrix};
\Sigma_{R_1}
\right)
}
{
V_{x_1,x_2,x_3}
\left(
1,
e^{\sigma d_{1,2}^{\alpha/2}u_2},
e^{\sigma d_{1,3}^{\alpha/2}u_3}
\right)
}.
\end{equation*}

It remains to prove the limiting distribution. Assume that
\[
\|x_2-x_1\|=\delta d_{1,2},
\qquad
\|x_3-x_1\|=\delta d_{1,3},
\qquad
\|x_3-x_2\|=\delta d_{2,3},
\]
where $d_{1,2},d_{1,3},d_{2,3}$ are fixed positive numbers. Then
\[
v_{1,2}(u_2)
=
u_2+\frac{1}{2}\sigma\delta^{\alpha/2}d_{1,2}^{\alpha/2}
\longrightarrow u_2,
\]
and
\[
v_{1,3}(u_3)
=
u_3+\frac{1}{2}\sigma\delta^{\alpha/2}d_{1,3}^{\alpha/2}
\longrightarrow u_3.
\]
The correlation coefficient $R_1$ is unchanged by the common scaling
$\delta$, since
\[
R_1
=
\frac{
d_{1,2}^{\alpha}+d_{1,3}^{\alpha}-d_{2,3}^{\alpha}
}
{
2(d_{1,2}d_{1,3})^{\alpha/2}
}.
\]
Hence the numerator converges to
\[
\Phi_2
\left(
\begin{pmatrix}
u_2\\
u_3
\end{pmatrix};
\Sigma_{R_1}
\right).
\]

We now show that the denominator converges to one. By the spectral
representation of the exponent function,
\begin{equation*}
V_{0,h_2,h_3}
\left(
1,
e^{\sigma\delta^{\alpha/2}d_{1,2}^{\alpha/2}u_2},
e^{\sigma\delta^{\alpha/2}d_{1,3}^{\alpha/2}u_3}
\right)
=
\mathbb{E}
\left[
\max\left\{
1,
e^{-\sigma\delta^{\alpha/2}d_{1,2}^{\alpha/2}u_2}Y(h_2),
e^{-\sigma\delta^{\alpha/2}d_{1,3}^{\alpha/2}u_3}Y(h_3)
\right\}
\right].
\end{equation*}
As $\delta\to0$, we have $h_2\to0$ and $h_3\to0$. Since the fractional
Brownian field has continuous sample paths,
\[
Y(h_2)\to Y(0)=1,
\qquad
Y(h_3)\to Y(0)=1
\]
almost surely. Moreover, the lognormal variables $Y(h_2)$ and $Y(h_3)$
converge to $1$ in $L^1$, because
\[
\mathbb{E}\{Y(h_i)\}=1,
\qquad
\operatorname{Var}(Y(h_i))
=
\exp\{2\gamma(h_i)\}-1
\longrightarrow 0,
\qquad i=2,3.
\]
It follows that
\[
\max\left\{
1,
e^{-\sigma\delta^{\alpha/2}d_{1,2}^{\alpha/2}u_2}Y(h_2),
e^{-\sigma\delta^{\alpha/2}d_{1,3}^{\alpha/2}u_3}Y(h_3)
\right\}
\longrightarrow 1
\]
in $L^1$. Consequently,
\[
V_{0,h_2,h_3}
\left(
1,
e^{\sigma\delta^{\alpha/2}d_{1,2}^{\alpha/2}u_2},
e^{\sigma\delta^{\alpha/2}d_{1,3}^{\alpha/2}u_3}
\right)
\longrightarrow 1.
\]
By stationarity, the same limit holds for
\[
V_{x_1,x_2,x_3}
\left(
1,
e^{\sigma\delta^{\alpha/2}d_{1,2}^{\alpha/2}u_2},
e^{\sigma\delta^{\alpha/2}d_{1,3}^{\alpha/2}u_3}
\right).
\]
Therefore,
\[
\lim_{\delta\to0}
\mathbb{P}
\left\{
U_{x_1,x_2}^{(\eta)}\leq u_2,
U_{x_1,x_3}^{(\eta)}\leq u_3
\right\}
=
\Phi_2
\left(
\begin{pmatrix}
u_2\\
u_3
\end{pmatrix};
\Sigma_{R_1}
\right),
\]
which concludes the proof.

\subsection{Proof of Theorem \ref{prop:BRtrajectories}}

\subsubsection{Proof for $V_{2,N}^{(\eta)}$}

Let
\[
H_2(u)=u^2-1,
\]
so that
\[
V_{2,N}^{(\eta)}
=
\frac{1}{\sqrt{|E_N|}}
\sum_{(x_1,x_2)\in E_N}
H_2\left(U_{x_1,x_2}^{(\eta)}\right).
\]
For $x_1\neq x_2$, write
\[
d_{1,2}=\|x_2-x_1\|
\]
and define, for each spectral field $W_j$,
\[
U_{x_1,x_2}^{(W_j)}
=
\frac{W_j(x_2)-W_j(x_1)}
{\sigma d_{1,2}^{\alpha/2}}.
\]

We split the proof into four steps.

\medskip
\noindent
\textbf{Step 1. Removing the deterministic variogram correction.}

Set
\[
b_{x_1,x_2}
=
\frac{\gamma(x_2)-\gamma(x_1)}
{\sigma d_{1,2}^{\alpha/2}}.
\]
Since \(H_2(u)=u^2-1\), we have the identity
\begin{equation}
H_2\left(U_{x_1,x_2}^{(\eta)}\right)
=
H_2\left(U_{x_1,x_2}^{(\eta)}+b_{x_1,x_2}\right)
-
2U_{x_1,x_2}^{(\eta)}b_{x_1,x_2}
-
b_{x_1,x_2}^2.
\label{eq:decompositionH2}
\end{equation}
We first show that the last two terms are negligible under the normalization
of Theorem \ref{prop:BRtrajectories}.

By Lemma \ref{Le_Add_Bounds}\textit{(ii)},
\[
\frac{1}{|E_N|}
\sum_{(x_1,x_2)\in E_N}
b_{x_1,x_2}^{2}
=
O_{\mathbb P}\left(N^{-1+\alpha/2}\right).
\]
Since $|E_N|=O_{\mathbb P}(N)$, it follows that
\[
N^{-(2-\alpha)/4}
\frac{1}{\sqrt{|E_N|}}
\sum_{(x_1,x_2)\in E_N}
b_{x_1,x_2}^{2}
=
O_{\mathbb P}\left(N^{3\alpha/4-1}\right)
\longrightarrow 0,
\]
because \(\alpha\in(0,1)\).

We next consider the linear term. By the Cauchy-Schwartz inequality,
\[
\left|
\frac{1}{|E_N|}
\sum_{(x_1,x_2)\in E_N}
U_{x_1,x_2}^{(\eta)}b_{x_1,x_2}
\right|
\leq
\left(
\frac{1}{|E_N|}
\sum_{(x_1,x_2)\in E_N}
\left(U_{x_1,x_2}^{(\eta)}\right)^2
\right)^{1/2}
\left(
\frac{1}{|E_N|}
\sum_{(x_1,x_2)\in E_N}
b_{x_1,x_2}^2
\right)^{1/2}.
\]
The first factor is \(O_{\mathbb P}(1)\), by Proposition
\ref{prop:uniformbounds} and the fact that \(|E_N|=O_{\mathbb P}(N)\). The
second factor is \(O_{\mathbb P}\left(N^{-1/2+\alpha/4}\right)\). Hence
\[
\frac{1}{|E_N|}
\sum_{(x_1,x_2)\in E_N}
U_{x_1,x_2}^{(\eta)}b_{x_1,x_2}
=
O_{\mathbb P}\left(N^{-1/2+\alpha/4}\right).
\]
Consequently,
\[
N^{-(2-\alpha)/4}
\frac{1}{\sqrt{|E_N|}}
\sum_{(x_1,x_2)\in E_N}
U_{x_1,x_2}^{(\eta)}b_{x_1,x_2}
=
O_{\mathbb P}\left(N^{(\alpha-1)/2}\right)
\longrightarrow 0,
\]
again because \(\alpha\in(0,1)\).

Thus it remains to study
\[
\frac{1}{\sqrt{|E_N|}}
\sum_{(x_1,x_2)\in E_N}
H_2\left(U_{x_1,x_2}^{(\eta)}+b_{x_1,x_2}\right).
\]

The next steps are devoted to this remaining contribution. We first decompose it according to the spectral functions that are active at the endpoints of each Delaunay edge; this separates the one-trajectory term, the two-trajectory crossing term, and a remainder term, whose asymptotic contributions are then analysed separately.

\medskip
\noindent
\textbf{Step 2. Decomposition according to the active spectral functions.}

Recall that
\[
Z_i(x)=\log U_i+\log Y_i(x)
=
\log U_i+W_i(x)-\gamma(x),
\]
and therefore
\[
\log\eta(x)=\bigvee_{i\geq1}Z_i(x).
\]
Since
\[
\log\eta(x_2)-\log\eta(x_1)+\gamma(x_2)-\gamma(x_1)
=
\bigvee_{i\geq1}\{\log U_i+W_i(x_2)\}
-
\bigvee_{i\geq1}\{\log U_i+W_i(x_1)\},
\]
the quantity
\[
U_{x_1,x_2}^{(\eta)}+b_{x_1,x_2}
\]
is the normalized increment of the upper envelope of the shifted Gaussian
spectral functions \(\log U_i+W_i\).

Let
\[
C_j=\{x\in\mathbf C: Z_j(x)=\log\eta(x)\}
\]
be the cells of the canonical tessellation. We decompose the contribution
according to whether the same spectral function is active at both endpoints
of the edge, or whether the active spectral function changes between
\(x_1\) and \(x_2\).

If the same index \(j\) is active at both endpoints, then
\[
U_{x_1,x_2}^{(\eta)}+b_{x_1,x_2}
=
U_{x_1,x_2}^{(W_j)}.
\]
If index \(j\) is active at \(x_1\) and index \(k\neq j\) is active at
\(x_2\), then
\[
U_{x_1,x_2}^{(\eta)}+b_{x_1,x_2}
=
U_{x_1,x_2}^{(W_k)}
+
\frac{Z_{k\setminus j}(x_1)}
{\sigma d_{1,2}^{\alpha/2}}.
\]
Consequently,
\begin{multline}
H_2\left(U_{x_1,x_2}^{(\eta)}+b_{x_1,x_2}\right)
=
\sum_{j\geq1}
\mathbb{I}\{x_1\in C_j\}
H_2\left(U_{x_1,x_2}^{(W_j)}\right)
\notag\\
+
\sum_{j\geq1}\sum_{k\neq j}
\mathbb{I}\{x_1\in C_j,\ x_2\in C_k\}
\left[
H_2\left(
U_{x_1,x_2}^{(W_k)}
+
\frac{Z_{k\setminus j}(x_1)}
{\sigma d_{1,2}^{\alpha/2}}
\right)
-
H_2\left(U_{x_1,x_2}^{(W_j)}\right)
\right].
\label{eq:active_decomposition_pair}
\end{multline}

We now reorganize the second term. For \(k>j\), let us recall that
\[
\mathbf C_{k,j}
=
\left\{
x\in\mathbf C:
Z_k(x)\wedge Z_j(x)
>
\bigvee_{i\neq j,k}Z_i(x)
\right\}
\]
is the region where \(Z_j\) and \(Z_k\) are the two largest spectral
functions (see Equation \eqref{eq:def_Ckj}). Define, for any function \(f:\mathbb R\to\mathbb R\),
\[
\Psi_f(x,y,w)
=
\{f(y+w)-f(x)\}\mathbb{I}\{x-y\leq w\leq0\}
+
\{f(x-w)-f(y)\}\mathbb{I}\{0\leq w\leq x-y\}.
\]
This function encodes the two possible orientations of a crossing between
the two spectral functions \(Z_j\) and \(Z_k\).

On the event \(x_1,x_2\in\mathbf C_{k,j}\), the condition that the active
index changes from \(j\) to \(k\) or from \(k\) to \(j\) is equivalent to a
crossing of \(Z_{k\setminus j}\) between \(x_1\) and \(x_2\). More precisely,
if \(j\) is active at \(x_1\) and \(k\) at \(x_2\), then
\[
Z_j(x_1)\geq Z_k(x_1),
\qquad
Z_k(x_2)\geq Z_j(x_2),
\]
which is equivalent to
\[
\frac{Z_{k\setminus j}(x_1)}
{\sigma d_{1,2}^{\alpha/2}}
\leq 0,
\qquad
\frac{Z_{k\setminus j}(x_1)}
{\sigma d_{1,2}^{\alpha/2}}
\geq
U_{x_1,x_2}^{(W_j)}
-
U_{x_1,x_2}^{(W_k)}.
\]
The reverse orientation gives the second term in \(\Psi_f\).

Thus $H_2\left(U_{x_1,x_2}^{(\eta)}+b_{x_1,x_2}\right)$ may be written as
\begin{equation}
H_2\left(U_{x_1,x_2}^{(\eta)}+b_{x_1,x_2}\right)
=
A_{x_1,x_2}
+
B_{x_1,x_2}
+
R_{x_1,x_2},
\label{eq:ABR_decomposition_pair}
\end{equation}
where
\[
A_{x_1,x_2}
=
\sum_{j\geq1}
\mathbb{I}\{x_1\in C_j\}
H_2\left(U_{x_1,x_2}^{(W_j)}\right),
\]
\[
B_{x_1,x_2}
=
\sum_{j\geq1}\sum_{k>j}
\mathbb{I}\{x_1,x_2\in\mathbf C_{k,j}\}
\Psi_{H_2}
\left(
U_{x_1,x_2}^{(W_j)},
U_{x_1,x_2}^{(W_k)},
\frac{Z_{k\setminus j}(x_1)}
{\sigma d_{1,2}^{\alpha/2}}
\right),
\]
and \(R_{x_1,x_2}\) collects the remaining configurations. These remaining
configurations are those for which the active indices at \(x_1\) and \(x_2\)
are \(j\) and \(k\), but at least one endpoint is not in the region
\(\mathbf C_{k,j}\). In other words, they correspond to local
configurations involving at least three competing spectral functions.

\medskip
\noindent
\textbf{Step 3. Negligibility of the one-trajectory and multi-interface
terms.}

We first consider the contribution of \(A_{x_1,x_2}\). Since the compact
window \(\mathbf C\) intersects only finitely many cells of the canonical
tessellation almost surely, the set
\[
\mathcal I_{\mathbf C}
=
\{j\geq1:C_j\cap\mathbf C\neq\varnothing\}
\]
is almost surely finite. Conditionally on the tessellation, the arguments of
Theorem 1 in \cite{Chenavier&Robert25a} yield
\[
\frac{1}{\sqrt{|E_N|}}
\sum_{(x_1,x_2)\in E_N}
A_{x_1,x_2}
=
O_{\mathbb P}(1).
\]
Equivalently, this term has the usual Gaussian order associated with squared
increment sums of fractional Brownian fields inside the cells. Therefore
\[
N^{-(2-\alpha)/4}
\frac{1}{\sqrt{|E_N|}}
\sum_{(x_1,x_2)\in E_N}
A_{x_1,x_2}
\overset{\mathbb P}{\longrightarrow}
0.
\]

We next consider the remainder term \(R_{x_1,x_2}\). This term corresponds
to configurations in which, at the scale of a Delaunay edge, more than two
spectral functions compete for the maximum. Such multi-interface
configurations are negligible at the scale of the two-trajectory crossing
term. More precisely, by the same multi-crossing estimate as in the proof of
Proposition 1 in \cite{Chenavier&Robert25b},
\[
N^{-(2-\alpha)/4}
\frac{1}{\sqrt{|E_N|}}
\sum_{(x_1,x_2)\in E_N}
R_{x_1,x_2}
\overset{\mathbb P}{\longrightarrow}
0.
\]
This estimate uses the almost sure finiteness of the number of spectral
cells intersecting \(\mathbf C\), together with the fact that simultaneous
crossings of three or more spectral functions occur on sets of smaller
effective dimension than the pairwise interfaces.

\medskip
\noindent
\textbf{Step 4. Limit of the two-trajectory crossing contribution.}

It remains to analyse the term \(B_{x_1,x_2}\). Let
\[
\mathcal J_{\mathbf C}
=
\{(j,k): j<k,\ \mathbf C_{k,j}\neq\varnothing\}.
\]
This set is almost surely finite. For each fixed pair \((j,k)\in
\mathcal J_{\mathbf C}\), Proposition 1 in \cite{Chenavier&Robert25b}, applied
on the random set \(\mathbf C_{k,j}\), gives
\begin{equation}
\label{eq:two_traj_limit_pair}
\begin{aligned}
&\frac{\sqrt{3}}{3}
N^{-(2-\alpha)/4}
\frac{1}{\sqrt{|E_N|}}
\sum_{(x_1,x_2)\in E_N}
\mathbb{I}\{x_1,x_2\in\mathbf C_{k,j}\}
\\
&\qquad \times
\Psi_{H_2}
\left(
U_{x_1,x_2}^{(W_j)},
U_{x_1,x_2}^{(W_k)},
\frac{Z_{k\setminus j}(x_1)}
{\sigma d_{1,2}^{\alpha/2}}
\right)
\overset{\mathbb P}{\longrightarrow}
c_{V_2}L_{Z_{k\setminus j}}(0).
\end{aligned}
\end{equation}
Since \(\mathcal J_{\mathbf C}\) is almost surely finite, we may sum
\eqref{eq:two_traj_limit_pair} over all pairs \(j<k\). We obtain
\begin{equation}
\label{eq:B_limit_pair}
\begin{aligned}
&\frac{\sqrt{3}}{3}
N^{-(2-\alpha)/4}
\frac{1}{\sqrt{|E_N|}}
\sum_{(x_1,x_2)\in E_N}
B_{x_1,x_2}
\\
&\qquad\overset{\mathbb P}{\longrightarrow}
c_{V_2}
\sum_{j\geq1}\sum_{k>j}
L_{Z_{k\setminus j}}(0).
\end{aligned}
\end{equation}

Combining \eqref{eq:decompositionH2}, the negligibility of the two
deterministic correction terms in Step 1, the decomposition
\eqref{eq:ABR_decomposition_pair}, the negligibility results of Step 3, and
the limit \eqref{eq:B_limit_pair}, we conclude that
\[
\frac{\sqrt{3}}{3}
N^{-(2-\alpha)/4}
V_{2,N}^{(\eta)}
\overset{\mathbb P}{\longrightarrow}
c_{V_2}
\sum_{j\geq1}\sum_{k>j}
L_{Z_{k\setminus j}}(0).
\]
This proves Theorem \ref{prop:BRtrajectories} for
\(V_{2,N}^{(\eta)}\).

\subsubsection{Proof for $V_{3,N}^{(\eta)}$}

We proceed as in the proof for $V_{2,N}^{(\eta)}$, but now the basic
increment statistic is a quadratic form in two normalized increments. For
$R\in(-1,1)$, define
\begin{equation}
\label{eq:def_H2R_triple}
H_{2,R}(u,v)
=
\begin{pmatrix}
u & v
\end{pmatrix}
\begin{pmatrix}
1 & R\\
R & 1
\end{pmatrix}^{-1}
\begin{pmatrix}
u\\
v
\end{pmatrix}
-2.
\end{equation}
Equivalently,
\[
H_{2,R}(u,v)
=
\frac{u^2+v^2-2Ruv}{1-R^2}-2.
\]
With this notation,
\[
V_{3,N}^{(\eta)}
=
\frac{1}{\sqrt{|DT_N|}}
\sum_{(x_1,x_2,x_3)\in DT_N}
H_{2,R_{x_1,x_2,x_3}}
\left(
U_{x_1,x_2}^{(\eta)},
U_{x_1,x_3}^{(\eta)}
\right),
\]
where \(R_{x_1,x_2,x_3}\) is defined in \eqref{eq:coorU}.

We split the proof into four steps.

\medskip
\noindent
\textbf{Step 1. Removing the deterministic variogram correction.}

For a triangle \(\Delta=(x_1,x_2,x_3)\in DT_N\), write
\[
d_{1,2}=\|x_2-x_1\|,
\qquad
d_{1,3}=\|x_3-x_1\|,
\qquad
R_{x_1,x_2,x_3}=R_{x_1,x_2,x_3}.
\]
Set
\[
b_{1,2}
=
\frac{\gamma(x_2)-\gamma(x_1)}
{\sigma d_{1,2}^{\alpha/2}},
\qquad
b_{1,3}
=
\frac{\gamma(x_3)-\gamma(x_1)}
{\sigma d_{1,3}^{\alpha/2}},
\]
and define the corrected increments
\[
\widetilde U_{1,2}^{(\eta)}
=
U_{x_1,x_2}^{(\eta)}+b_{1,2},
\qquad
\widetilde U_{1,3}^{(\eta)}
=
U_{x_1,x_3}^{(\eta)}+b_{1,3}.
\]
We first show that replacing
\[
\left(U_{x_1,x_2}^{(\eta)},U_{x_1,x_3}^{(\eta)}\right)
\]
by
\[
\left(\widetilde U_{1,2}^{(\eta)},\widetilde U_{1,3}^{(\eta)}\right)
\]
does not change the limit under the normalization of the theorem.

Indeed, since \(H_{2,R}\) is a quadratic form, there exists a constant
\(C_R\), proportional to \((1-R^2)^{-1}\), such that
\begin{align*}
&\left|
H_{2,R}
\left(
U_{x_1,x_2}^{(\eta)},U_{x_1,x_3}^{(\eta)}
\right)
-
H_{2,R}
\left(
\widetilde U_{1,2}^{(\eta)},\widetilde U_{1,3}^{(\eta)}
\right)
\right|
\\
&\qquad\leq
C_R
\left[
\left(
\left|U_{x_1,x_2}^{(\eta)}\right|
+
\left|U_{x_1,x_3}^{(\eta)}\right|
\right)
\left(|b_{1,2}|+|b_{1,3}|\right)
+
b_{1,2}^2+b_{1,3}^2+|b_{1,2}b_{1,3}|
\right].
\end{align*}
The integrability of the factor \(C_R\) with respect to the typical
Poisson--Delaunay triangle follows from the estimates used in the proof of
Theorem 1 in \cite{Chenavier&Robert25a}. Together with Proposition
\ref{prop:uniformbounds} and Lemma \ref{Le_Add_Bounds}, applied to the
Delaunay edges of the triangles, this gives
\begin{equation}
\label{eq:gamma_correction_triple}
N^{-(2-\alpha)/4}
\frac{1}{\sqrt{|DT_N|}}
\sum_{(x_1,x_2,x_3)\in DT_N}
\bigg|
H_{2,R_{x_1,x_2,x_3}}
\left(
U_{x_1,x_2}^{(\eta)},U_{x_1,x_3}^{(\eta)}
\right)
-
H_{2,R_{x_1,x_2,x_3}}
\left(
\widetilde U_{1,2}^{(\eta)},\widetilde U_{1,3}^{(\eta)}
\right)
\bigg|
\overset{\mathbb P}{\longrightarrow}
0.
\end{equation}
It remains to study the corrected statistic.

\medskip
\noindent
\textbf{Step 2. Decomposition according to the active spectral functions.}

Recall that
\[
Z_i(x)
=
\log U_i+\log Y_i(x)
=
\log U_i+W_i(x)-\gamma(x),
\]
and hence
\[
\log\eta(x)=\bigvee_{i\geq1}Z_i(x).
\]
The deterministic correction by \(\gamma\) gives
\[
\widetilde U_{1,r}^{(\eta)}
=
\frac{
\bigvee_{i\geq1}\{\log U_i+W_i(x_r)\}
-
\bigvee_{i\geq1}\{\log U_i+W_i(x_1)\}
}
{\sigma d_{1,r}^{\alpha/2}},
\qquad r=2,3.
\]
For each spectral field \(W_j\), define
\[
U_{x_1,x_r}^{(W_j)}
=
\frac{W_j(x_r)-W_j(x_1)}
{\sigma d_{1,r}^{\alpha/2}},
\qquad r=2,3.
\]

Let
\[
C_j=\{x\in\mathbf C:Z_j(x)=\log\eta(x)\}
\]
be the cells of the canonical tessellation. If the same spectral function
\(Z_j\) is active at \(x_1,x_2,x_3\), then
\[
\left(
\widetilde U_{1,2}^{(\eta)},
\widetilde U_{1,3}^{(\eta)}
\right)
=
\left(
U_{x_1,x_2}^{(W_j)},
U_{x_1,x_3}^{(W_j)}
\right).
\]
This gives the one-trajectory contribution
\[
A_{x_1,x_2,x_3}
=
\sum_{j\geq1}
\mathbb{I}\{x_1\in C_j\}
H_{2,R_{x_1,x_2,x_3}}
\left(
U_{x_1,x_2}^{(W_j)},
U_{x_1,x_3}^{(W_j)}
\right).
\]

We now define the two-trajectory crossing contribution. For any function
\(f:\mathbb R\to\mathbb R\), recall
\[
\Psi_f(x,y,w)
=
\{f(y+w)-f(x)\}\mathbb{I}\{x-y\leq w\leq0\}
+
\{f(x-w)-f(y)\}\mathbb{I}\{0\leq w\leq x-y\}.
\]
Let \(I(u)=u\). For \(k>j\), set
\[
w_{1,2}
=
\frac{Z_{k\setminus j}(x_1)}
{\sigma d_{1,2}^{\alpha/2}},
\qquad
w_{1,3}
=
\frac{Z_{k\setminus j}(x_1)}
{\sigma d_{1,3}^{\alpha/2}}.
\]
For \(a,b,c,d,w_1,w_2\in\mathbb R\) and \(R\in(-1,1)\), define
\begin{align}
\label{eq:def_Lambda_triple}
\Lambda(a,b,c,d,w_1,w_2;R)
&=
\left[
H_{2,R}
\left(
a+\Psi_I(a,c,w_1),
b+\Psi_I(b,d,w_2)
\right)
-
H_{2,R}(a,b)
\right]\mathbb{I}\{w_1<0\}
\notag\\
&\quad+
\left[
H_{2,R}
\left(
c+\Psi_I(a,c,w_1),
d+\Psi_I(b,d,w_2)
\right)
-
H_{2,R}(c,d)
\right]\mathbb{I}\{w_1>0\}.
\end{align}
The sign of \(w_1\) determines which of the two spectral functions is active
at \(x_1\). Since \(w_1\) and \(w_2\) have the same sign, this definition
covers the two possible orientations of the crossing between \(Z_j\) and
\(Z_k\).

The corresponding two-trajectory term is
\[
B_{x_1,x_2,x_3}
=
\sum_{j\geq1}\sum_{k>j}
\mathbb{I}\{x_1,x_2,x_3\in\mathbf C_{k,j}\}
\Lambda
\left(
U_{x_1,x_2}^{(W_j)},
U_{x_1,x_3}^{(W_j)},
U_{x_1,x_2}^{(W_k)},
U_{x_1,x_3}^{(W_k)},
w_{1,2},
w_{1,3};
R_{x_1,x_2,x_3}
\right).
\]
Finally, let \(R_{x_1,x_2,x_3}^{\mathrm{rem}}\) collect all remaining configurations.
These are configurations in which at least one vertex of the triangle is not
in the region \(\mathbf C_{k,j}\) associated with the two spectral functions
that are active at two of the vertices, or configurations involving three or
more competing spectral functions. We therefore have the decomposition
\begin{equation}
\label{eq:ABR_decomposition_triple}
H_{2,R_{x_1,x_2,x_3}}
\left(
\widetilde U_{1,2}^{(\eta)},
\widetilde U_{1,3}^{(\eta)}
\right)
=
A_{x_1,x_2,x_3}+B_{x_1,x_2,x_3}+R_{x_1,x_2,x_3}^{\mathrm{rem}}.
\end{equation}

\medskip
\noindent
\textbf{Step 3. Negligibility of the one-trajectory and remainder terms.}

We first consider the one-trajectory term. The set
\[
\mathcal I_{\mathbf C}
=
\{j\geq1:C_j\cap\mathbf C\neq\varnothing\}
\]
is almost surely finite. Conditionally on the canonical tessellation, the
arguments used in the proof of Theorem 1 in \cite{Chenavier&Robert25a}
imply that
\[
\frac{1}{\sqrt{|DT_N|}}
\sum_{\Delta\in DT_N}
A_{x_1,x_2,x_3}
=
O_{\mathbb P}(1).
\]
This is the usual Gaussian order for squared increment sums of an isotropic
fractional Brownian field over Delaunay triangles. Hence
\begin{equation}
\label{eq:A_negligible_triple}
N^{-(2-\alpha)/4}
\frac{1}{\sqrt{|DT_N|}}
\sum_{\Delta\in DT_N}
A_{x_1,x_2,x_3}
\overset{\mathbb P}{\longrightarrow}
0.
\end{equation}

We now consider the remainder. The term \(R_{x_1,x_2,x_3}^{\mathrm{rem}}\) corresponds
to configurations in which, at the scale of a Delaunay triangle, at least
three spectral functions compete, or in which a crossing along one edge is
not associated with a triangle entirely contained in the corresponding
two-trajectory region \(\mathbf C_{k,j}\). Such configurations are
negligible at the scale of the two-trajectory crossing term. More precisely,
by the same multi-crossing estimates as those used in the proof of
Proposition 1 in \cite{Chenavier&Robert25b}, combined with the triangular
moment estimates from \cite{Chenavier&Robert25a}, we have
\begin{equation}
\label{eq:R_negligible_triple}
N^{-(2-\alpha)/4}
\frac{1}{\sqrt{|DT_N|}}
\sum_{\Delta \in DT_N}
R_{x_1,x_2,x_3}^{\mathrm{rem}}
\overset{\mathbb P}{\longrightarrow}
0.
\end{equation}

\medskip
\noindent
\textbf{Step 4. Limit of the two-trajectory crossing contribution.}

Let
\[
\mathcal J_{\mathbf C}
=
\{(j,k):j<k,\ \mathbf C_{k,j}\neq\varnothing\}.
\]
This set is almost surely finite. For each fixed pair
\((j,k)\in\mathcal J_{\mathbf C}\), the triangular analogue of Proposition 1
in \cite{Chenavier&Robert25b}, applied on the random set
\(\mathbf C_{k,j}\), gives
\begin{equation}
\label{eq:two_traj_limit_triple}
\begin{aligned}
&\frac{\sqrt{2}}{2}
N^{-(2-\alpha)/4}
\frac{1}{\sqrt{|DT_N|}}
\sum_{(x_1,x_2,x_3)\in DT_N}
\mathbb{I}\{x_1,x_2,x_3\in\mathbf C_{k,j}\}
\\
&\qquad\times
\Lambda
\left(
U_{x_1,x_2}^{(W_j)},
U_{x_1,x_3}^{(W_j)},
U_{x_1,x_2}^{(W_k)},
U_{x_1,x_3}^{(W_k)},
\frac{Z_{k\setminus j}(x_1)}
{\sigma d_{1,2}^{\alpha/2}},
\frac{Z_{k\setminus j}(x_1)}
{\sigma d_{1,3}^{\alpha/2}};
R_{x_1,x_2,x_3}
\right)
\\
&\qquad\overset{\mathbb P}{\longrightarrow}
c_{V_3}L_{Z_{k\setminus j}}(0).
\end{aligned}
\end{equation}
Since \(\mathcal J_{\mathbf C}\) is almost surely finite, summing
\eqref{eq:two_traj_limit_triple} over all \(j<k\) yields
\begin{equation}
\label{eq:B_limit_triple}
\frac{\sqrt{2}}{2}
N^{-(2-\alpha)/4}
\frac{1}{\sqrt{|DT_N|}}
\sum_{\Delta \in DT_N}
B_{x_1,x_2,x_3}
\overset{\mathbb P}{\longrightarrow}
c_{V_3}
\sum_{j\geq1}\sum_{k>j}
L_{Z_{k\setminus j}}(0).
\end{equation}

Combining \eqref{eq:gamma_correction_triple},
\eqref{eq:ABR_decomposition_triple},
\eqref{eq:A_negligible_triple},
\eqref{eq:R_negligible_triple}, and \eqref{eq:B_limit_triple}, we obtain
\[
\frac{\sqrt{2}}{2}
N^{-(2-\alpha)/4}
V_{3,N}^{(\eta)}
\overset{\mathbb P}{\longrightarrow}
c_{V_3}
\sum_{j\geq1}\sum_{k>j}
L_{Z_{k\setminus j}}(0).
\]
This proves Theorem \ref{prop:BRtrajectories} for
\(V_{3,N}^{(\eta)}\).

\subsection{Proof of Proposition \protect\ref{Prop_pairwise_pdf}}

Throughout the proof, the observations $z_1,z_2$ and the distance
$d=\|x_2-x_1\|$ are kept fixed when differentiating with respect to
$\sigma$ or $\alpha$. We set
\[
a=\sigma d^{\alpha/2},
\qquad
q=\log(z_2/z_1),
\qquad
u=\frac{q}{a},
\qquad
v_+(u)=u+\frac a2,
\qquad
v_-(u)=-u+\frac a2 .
\]
Thus, in the local regime considered in the proposition, $u$ is fixed and
$a\to0$. However, when differentiating with respect to $\sigma$ or
$\alpha$, the quantity $q$ is fixed and therefore $u=q/a$ depends on the
parameter.

The bivariate exponent function can be written as
\[
V_{x_1,x_2}(z_1,z_2)
=
\frac{1}{z_1}
\left\{
\Phi(v_+(u))+e^{-au}\Phi(v_-(u))
\right\}.
\]
Since
\[
F_{x_1,x_2}(z_1,z_2)
=
\exp\{-V_{x_1,x_2}(z_1,z_2)\},
\]
the density is
\[
f_{x_1,x_2}(z_1,z_2)
=
\exp\{-V_{x_1,x_2}(z_1,z_2)\}
\left(
\frac{\partial V}{\partial z_1}
\frac{\partial V}{\partial z_2}
-
\frac{\partial^2V}{\partial z_1\partial z_2}
\right).
\]
A direct differentiation gives
\[
\frac{\partial V}{\partial z_1}
=
-\frac{1}{z_1^2}A_a(u),
\qquad
\frac{\partial V}{\partial z_2}
=
-\frac{1}{z_2^2}B_a(u),
\qquad
\frac{\partial^2V}{\partial z_1\partial z_2}
=
-\frac{1}{z_1^3}C_a(u),
\]
where
\[
A_a(u)
=
\Phi(v_+(u))
+
\frac{1}{a}\varphi(v_+(u))
-
\frac{1}{a}e^{-au}\varphi(v_-(u)),
\]
\[
B_a(u)
=
\Phi(v_-(u))
+
\frac{1}{a}\varphi(v_-(u))
-
\frac{1}{a}e^{au}\varphi(v_+(u)),
\]
and
\[
C_a(u)
=
\frac{1}{a^2}e^{-au}
\left\{
v_-(u)\varphi(v_+(u))
+
e^{-au}v_+(u)\varphi(v_-(u))
\right\}.
\]
Therefore, using $z_2=z_1e^{au}$,
\[
f_{x_1,x_2}(z_1,z_2)
=
\exp\left(
-\frac{1}{z_1}
\left\{
\Phi(v_+(u))+e^{-au}\Phi(v_-(u))
\right\}
\right)
\left\{
\frac{e^{-2au}}{z_1^4}A_a(u)B_a(u)
+
\frac{1}{z_1^3}C_a(u)
\right\}.
\]
Equivalently,
\begin{equation}
\label{eq:log_pair_density_expansion_start}
\log f_{x_1,x_2}(z_1,z_2)
=
-V_a(u)
-3\log z_1
+
\log C_a(u)
+
\log\left(
1+
\frac{e^{-2au}A_a(u)B_a(u)}{z_1C_a(u)}
\right),
\end{equation}
where
\[
V_a(u)
=
\frac{1}{z_1}
\left\{
\Phi(v_+(u))+e^{-au}\Phi(v_-(u))
\right\}.
\]

We now derive the expansion of the derivative with respect to $\sigma$.
Since $q=\log(z_2/z_1)$ is fixed under differentiation, the differential
operator $\sigma\partial_\sigma$ is
\[
\sigma\frac{\partial}{\partial\sigma}
=
a\frac{\partial}{\partial a}\Big|_{q\ \mathrm{fixed}}
=
a\frac{\partial}{\partial a}
-
u\frac{\partial}{\partial u}
=:\mathcal D.
\]
We shall prove that
\begin{equation}
\label{eq:main_D_logf_pair}
\mathcal D\log f_{x_1,x_2}(z_1,z_2)
=
u^2-1
+
\frac{a}{z_1}\omega(u)
+
o(a),
\end{equation}
where
\[
\omega(u)
=
u\{2\Phi(u)-1\}
+
\frac{(1-u^2)\Phi(u)\{1-\Phi(u)\}}{\varphi(u)}
-
\varphi(u).
\]

We first record the elementary Taylor expansions, valid for fixed
$u\in\mathbb R$ as $a\to0$:
\[
\Phi(v_+(u))
=
\Phi(u)+\frac a2\varphi(u)+O(a^2),
\qquad
\Phi(v_-(u))
=
\Phi(-u)+\frac a2\varphi(u)+O(a^2),
\]
\[
\varphi(v_+(u))
=
\varphi(u)\left(1-\frac a2u+O(a^2)\right),
\qquad
\varphi(v_-(u))
=
\varphi(u)\left(1+\frac a2u+O(a^2)\right).
\]
From these expansions, one obtains
\[
V_a(u)
=
\frac{1}{z_1}
\left[
1+a\{\varphi(u)-u\bar\Phi(u)\}+O(a^2)
\right],
\]
where $\bar\Phi(u)=1-\Phi(u)=\Phi(-u)$. Hence, with
\[
V_1(u)=\varphi(u)-u\bar\Phi(u),
\]
we have
\[
\mathcal D V_a(u)
=
\frac{a}{z_1}\{V_1(u)-uV_1'(u)\}+O(a^2).
\]
Since
\[
V_1'(u)=-\bar\Phi(u),
\]
it follows that
\begin{equation}
\label{eq:D_V_pair}
\mathcal D\{-V_a(u)\}
=
-\frac{a}{z_1}\varphi(u)+O(a^2).
\end{equation}

Next, the same Taylor expansions give
\[
A_a(u)=\Phi(u)+\frac a2\varphi(u)+O(a^2),
\]
\[
B_a(u)=\bar\Phi(u)+\frac a2\varphi(u)+O(a^2),
\]
and
\[
C_a(u)
=
\frac{\varphi(u)}{a}
-\frac32u\varphi(u)
+
O(a).
\]
Consequently,
\[
\log C_a(u)
=
\log\varphi(u)-\log a-\frac32au+O(a^2).
\]
Applying $\mathcal D$ gives
\[
\mathcal D\log C_a(u)
=
u^2-1+O(a^2),
\]
because the first-order term $-\frac32au$ is annihilated by
$\mathcal D$, namely
\[
\mathcal D(au)=a u-u a=0.
\]

Finally,
\[
\frac{e^{-2au}A_a(u)B_a(u)}{z_1C_a(u)}
=
\frac{a}{z_1}
\frac{\Phi(u)\bar\Phi(u)}{\varphi(u)}
+
O(a^2).
\]
Thus
\[
\log\left(
1+
\frac{e^{-2au}A_a(u)B_a(u)}{z_1C_a(u)}
\right)
=
\frac{a}{z_1}
D(u)
+
O(a^2),
\]
where
\[
D(u)=\frac{\Phi(u)\bar\Phi(u)}{\varphi(u)}.
\]
Therefore
\begin{equation}
\label{eq:D_D_pair}
\mathcal D
\log\left(
1+
\frac{e^{-2au}A_a(u)B_a(u)}{z_1C_a(u)}
\right)
=
\frac{a}{z_1}\{D(u)-uD'(u)\}
+
O(a^2).
\end{equation}
Note that the factor \(\sigma\) is absorbed in the definition of the operator \(\mathcal D\): indeed, since \(a=\sigma d^{\alpha/2}\) and \(q=\log(z_2/z_1)\) is kept fixed when differentiating, one has
\[
\mathcal D=\sigma\frac{\partial}{\partial\sigma}
=
a\frac{\partial}{\partial a}-u\frac{\partial}{\partial u}.
\]
Thus no additional factor \(\sigma\) appears when applying \(\mathcal D\); it only reappears when passing back from \(\mathcal D\) to \(\partial/\partial\sigma\).
Since
\[
D'(u)
=
1-2\Phi(u)+uD(u),
\]
we get
\[
D(u)-uD'(u)
=
(1-u^2)D(u)+u\{2\Phi(u)-1\}.
\]
Combining \eqref{eq:D_V_pair} and \eqref{eq:D_D_pair}, we obtain
\[
\mathcal D\log f_{x_1,x_2}(z_1,z_2)
=
u^2-1
+
\frac{a}{z_1}
\left[
-\varphi(u)
+
(1-u^2)\frac{\Phi(u)\bar\Phi(u)}{\varphi(u)}
+
u\{2\Phi(u)-1\}
\right]
+
o(a).
\]
This is precisely \eqref{eq:main_D_logf_pair}.

Then we deduce that
\[
\frac{\partial}{\partial\sigma}
\log f_{x_1,x_2}(z_1,z_2;\sigma,\alpha)
=
\frac{1}{\sigma}(u^2-1)
+
\frac{a}{\sigma z_1}\omega(u)
+
o\left(\frac{a}{\sigma}\right).
\]
As $a=\sigma d^{\alpha/2}$, this becomes
\begin{equation}
\label{eq:pair_score_sigma_final}
\frac{\partial}{\partial\sigma}
\log f_{x_1,x_2}(z_1,z_2;\sigma,\alpha)
=
\frac{1}{\sigma}(u^2-1)
+
\frac{\omega(u)}{z_1}d^{\alpha/2}
+
o(d^{\alpha/2}).
\end{equation}

We now turn to the derivative with respect to $\alpha$. Since
\[
a=\sigma d^{\alpha/2},
\]
we have, again with $z_1,z_2$ and $d$ fixed,
\[
\frac{\partial}{\partial\alpha}
=
\frac{\log d}{2}
a\frac{\partial}{\partial a}\Big|_{q\ \mathrm{fixed}}
=
\frac{\log d}{2}\mathcal D.
\]
Therefore, by \eqref{eq:main_D_logf_pair},
\[
\frac{1}{\log d}
\frac{\partial}{\partial\alpha}
\log f_{x_1,x_2}(z_1,z_2;\sigma,\alpha)
=
\frac12(u^2-1)
+
\frac{a}{2z_1}\omega(u)
+
o(a).
\]
Using again $a=\sigma d^{\alpha/2}$, we obtain
\begin{equation}
\label{eq:pair_score_alpha_final}
\frac{1}{\log d}
\frac{\partial}{\partial\alpha}
\log f_{x_1,x_2}(z_1,z_2;\sigma,\alpha)
=
\frac12(u^2-1)
+
\frac{\sigma\,\omega(u)}{2z_1}d^{\alpha/2}
+
o(d^{\alpha/2}).
\end{equation}
This proves the announced pairwise score expansions, with the correction
terms given in \eqref{eq:pair_score_sigma_final} and
\eqref{eq:pair_score_alpha_final}.

\subsection{Proof of Proposition \protect\ref{Prop_triplewise_pdf}}

Let
\[
\|x_2-x_1\|=\delta d_{1,2},\qquad
\|x_3-x_1\|=\delta d_{1,3},\qquad
\|x_3-x_2\|=\delta d_{2,3},
\]
where \(d_{1,2},d_{1,3},d_{2,3}\) are fixed positive numbers. We set
\[
a_{i,j}
=
\sigma\delta^{\alpha/2}d_{i,j}^{\alpha/2},
\qquad
u_{i,j}
=
\frac{\log(z_j/z_i)}{a_{i,j}},
\qquad
v_{i,j}(t)
=
t+\frac{a_{i,j}}{2}.
\]
In the proposition, \(u_2\) and \(u_3\) are fixed and correspond to
\[
u_{1,2}=u_2,
\qquad
u_{1,3}=u_3.
\]
The third normalized log-ratio is then determined by
\[
u_{2,3}
=
\frac{\log(z_3/z_2)}{a_{2,3}}
=
\frac{a_{1,3}u_3-a_{1,2}u_2}{a_{2,3}}
=
\frac{d_{1,3}^{\alpha/2}u_3-d_{1,2}^{\alpha/2}u_2}
{d_{2,3}^{\alpha/2}}.
\]
Thus \(u_{2,3}\) is independent of \(\delta\).

We use the notation
\[
\Sigma_R=
\begin{pmatrix}
1 & R\\
R & 1
\end{pmatrix},
\qquad
Q_R(u_2,u_3)
=
\begin{pmatrix}
u_2 & u_3
\end{pmatrix}
\Sigma_R^{-1}
\begin{pmatrix}
u_2\\
u_3
\end{pmatrix}
-2.
\]

The trivariate density of a max-stable distribution can be written in terms
of the exponent function \(V=V_{x_1,x_2,x_3}\) as
\begin{equation}
\label{eq:triple_density_partition}
f_{x_1,x_2,x_3}
=
e^{-V}
\left\{
-V_{123}
+
V_1V_{23}
+
V_2V_{13}
+
V_3V_{12}
-
V_1V_2V_3
\right\},
\end{equation}
where \(V_i=\partial V/\partial z_i\),
\(V_{ij}=\partial^2V/\partial z_i\partial z_j\), and
\(V_{123}=\partial^3V/\partial z_1\partial z_2\partial z_3\).

The term \(-V_{123}\) corresponds to the event that the three sites belong
to the same cell of the canonical max-stable tessellation. In the local
regime \(\delta\to0\), this is the dominant term. More precisely, using the
explicit expressions of the derivatives of the Brown--Resnick exponent
function, one obtains
\begin{equation}
\label{eq:triple_density_expansion}
f_{x_1,x_2,x_3}(z_1,z_2,z_3;\sigma,\alpha)
=
e^{-V}
(-V_{123})
\left[
1+
\delta^{\alpha/2}
\frac{\Xi(u_2,u_3;d_{1,2},d_{1,3},d_{2,3})}{z_1}
+
o(\delta^{\alpha/2})
\right],
\end{equation}
for some measurable function \(\Xi\). The explicit form of \(\Xi\) is obtained by inserting into \eqref{eq:triple_density_partition} the trivariate H\"usler--Reiss representation of the Brown--Resnick exponent function displayed in \eqref{Eq_V_triple}, that is, the expression of \(V_{x_1,x_2,x_3}\) as a sum of three bivariate Gaussian distribution functions with arguments \(v_{i,j}(\pm u_{i,j})\) and correlation coefficients \(R_1,R_2,R_3\), together with its derivatives with respect to \(z_1,z_2,z_3\).
 It involves derivatives of the
bivariate Gaussian distribution function and is not used below. What matters
is that \(\Xi\) depends on the local normalized increments and on the
limiting triangle shape, but not on \(\delta\).

We now compute the dominant term. From Appendix B.4 of \cite{Dombry18},
\begin{multline}
\label{eq:V123_BR}
V_{123}
=
-\frac{1}{z_1^2z_2z_3}
\varphi_2
\left(
\begin{pmatrix}
a_{1,2}v_{1,2}(u_{1,2})\\
a_{1,3}v_{1,3}(u_{1,3})
\end{pmatrix};
\begin{pmatrix}
a_{1,2}^2 &
\{a_{1,2}^2+a_{1,3}^2-a_{2,3}^2\}/2\\
\{a_{1,2}^2+a_{1,3}^2-a_{2,3}^2\}/2 &
a_{1,3}^2
\end{pmatrix}
\right).
\end{multline}
Equivalently,
\begin{multline}
\label{eq:V123_scaled}
-V_{123}
=
\frac{1}{z_1^2z_2z_3}
\frac{1}{a_{1,2}a_{1,3}}
\varphi_2
\left(
\begin{pmatrix}
v_{1,2}(u_2)\\
v_{1,3}(u_3)
\end{pmatrix};
\Sigma_{R_1}
\right),
\end{multline}
where
\[
R_1
=
\frac{
d_{1,2}^{\alpha}+d_{1,3}^{\alpha}-d_{2,3}^{\alpha}
}
{
2(d_{1,2}d_{1,3})^{\alpha/2}
}.
\]
Since
\[
z_2=z_1e^{a_{1,2}u_2},
\qquad
z_3=z_1e^{a_{1,3}u_3},
\]
we have
\begin{multline}
\label{eq:log_V123}
\log(-V_{123})
=
-4\log z_1
-a_{1,2}u_2-a_{1,3}u_3
-\log(a_{1,2}a_{1,3})
\\
-\frac12\log(1-R_1^2)
-\log(2\pi)
-\frac12
\begin{pmatrix}
v_{1,2}(u_2)&v_{1,3}(u_3)
\end{pmatrix}
\Sigma_{R_1}^{-1}
\begin{pmatrix}
v_{1,2}(u_2)\\
v_{1,3}(u_3)
\end{pmatrix}.
\end{multline}

We first differentiate with respect to \(\sigma\), keeping
\(z_1,z_2,z_3\) and the locations fixed. Under this differentiation,
\[
\sigma\frac{\partial a_{i,j}}{\partial\sigma}=a_{i,j},
\qquad
\sigma\frac{\partial v_{i,j}(u_{i,j})}{\partial\sigma}
=
v_{i,j}(-u_{i,j}).
\]
Moreover, \(R_1\) does not depend on \(\sigma\). Differentiating
\eqref{eq:log_V123}, we obtain
\begin{multline}
\label{eq:sigma_der_logV123}
\sigma
\frac{\partial}{\partial\sigma}
\log(-V_{123})
=
-a_{1,2}u_2-a_{1,3}u_3-2
\\
-
\left[
\Sigma_{R_1}^{-1}
\begin{pmatrix}
v_{1,2}(u_2)\\
v_{1,3}(u_3)
\end{pmatrix}
\right]_1
v_{1,2}(-u_2)
-
\left[
\Sigma_{R_1}^{-1}
\begin{pmatrix}
v_{1,2}(u_2)\\
v_{1,3}(u_3)
\end{pmatrix}
\right]_2
v_{1,3}(-u_3).
\end{multline}
As \(\delta\to0\), we have \(a_{1,2}\to0\), \(a_{1,3}\to0\), and hence
\[
v_{1,2}(u_2)\to u_2,
\qquad
v_{1,3}(u_3)\to u_3,
\qquad
v_{1,2}(-u_2)\to -u_2,
\qquad
v_{1,3}(-u_3)\to -u_3.
\]
Therefore,
\begin{equation}
\label{eq:sigma_der_logV123_limit}
\sigma
\frac{\partial}{\partial\sigma}
\log(-V_{123})
=
Q_{R_1}(u_2,u_3)
+
O(\delta^{\alpha/2}).
\end{equation}

We now return to the full density. Taking logarithms in
\eqref{eq:triple_density_expansion} gives
\begin{equation}
\label{eq:log_triple_density_expansion}
\log f_{x_1,x_2,x_3}
=
\log(-V_{123})
-
V
+
\delta^{\alpha/2}
\frac{\Xi(u_2,u_3;d_{1,2},d_{1,3},d_{2,3})}{z_1}
+
o(\delta^{\alpha/2}).
\end{equation}
Differentiating this identity with respect to \(\sigma\) and using the same
operator as above yields the expansion
\begin{equation}
\label{eq:sigma_score_triple}
\frac{\partial}{\partial\sigma}
\log f_{x_1,x_2,x_3}(z_1,z_2,z_3;\sigma,\alpha)
=
\frac{1}{\sigma}Q_{R_1}(u_2,u_3)
+
\frac{\delta^{\alpha/2}}{z_1}
\Omega_\sigma(u_2,u_3;d_{1,2},d_{1,3},d_{2,3})
+
o(\delta^{\alpha/2}),
\end{equation}
where \(\Omega_\sigma\) is a measurable function collecting the first-order
contributions coming from \(-V\), from the non-dominant partition terms in
\eqref{eq:triple_density_partition}, and from the first-order expansion of
\(\log(-V_{123})\). This proves the expansion of the score with respect to
\(\sigma\).

We now consider the derivative with respect to \(\alpha\). Since
\[
a_{i,j}
=
\sigma\delta^{\alpha/2}d_{i,j}^{\alpha/2},
\]
we have
\[
\frac{\partial a_{i,j}}{\partial\alpha}
=
\frac12\log(\delta d_{i,j})\,a_{i,j}.
\]
Moreover, when \(z_1,z_2,z_3\) and the locations are kept fixed,
\[
u_{i,j}
=
\frac{\log(z_j/z_i)}{a_{i,j}},
\]
and therefore
\[
\frac{\partial u_{i,j}}{\partial\alpha}
=
-\frac12\log(\delta d_{i,j})\,u_{i,j}.
\]
Consequently,
\[
\frac{\partial v_{i,j}(u_{i,j})}{\partial\alpha}
=
\frac12\log(\delta d_{i,j})\,v_{i,j}(-u_{i,j}).
\]

In contrast with the derivative with respect to \(\sigma\), the correlation
coefficient
\[
R_1
=
\frac{
d_{1,2}^{\alpha}+d_{1,3}^{\alpha}-d_{2,3}^{\alpha}
}
{
2(d_{1,2}d_{1,3})^{\alpha/2}
}
\]
also depends on \(\alpha\). We write
\[
\dot R_1(\alpha)
=
\frac{\partial R_1}{\partial\alpha}.
\]
Explicitly,
\[
\dot R_1(\alpha)
=
\frac{
d_{1,2}^{\alpha}\log d_{1,2}
+
d_{1,3}^{\alpha}\log d_{1,3}
-
d_{2,3}^{\alpha}\log d_{2,3}
-
\frac12\log(d_{1,2}d_{1,3})
\left(
d_{1,2}^{\alpha}
+
d_{1,3}^{\alpha}
-
d_{2,3}^{\alpha}
\right)
}
{
2(d_{1,2}d_{1,3})^{\alpha/2}
}.
\]

Let
\[
\mathbf u=
\begin{pmatrix}
u_2\\
u_3
\end{pmatrix},
\qquad
A_{R_1}
=
\Sigma_{R_1}^{-1}\mathbf u
=
\begin{pmatrix}
A_{R_1,1}\\
A_{R_1,2}
\end{pmatrix}.
\]
Recall that
\[
Q_{R_1}(u_2,u_3)
=
\mathbf u^\top\Sigma_{R_1}^{-1}\mathbf u-2.
\]
We shall also use the notation
\[
B_{R_1}(u_2,u_3)
=
\frac{R_1}{1-R_1^2}
+
\frac{
(1+R_1^2)u_2u_3
-
R_1(u_2^2+u_3^2)
}
{(1-R_1^2)^2}.
\]
This quantity is the derivative of the \(R_1\)-dependent part of
\(\log(-V_{123})\) with respect to \(R_1\), evaluated at
\((v_{1,2}(u_2),v_{1,3}(u_3))=(u_2,u_3)\).

Differentiating \eqref{eq:log_V123} with respect to \(\alpha\), dividing by
\(\log\delta\), and using
\[
v_{1,2}(u_2)=u_2+O(\delta^{\alpha/2}),
\qquad
v_{1,3}(u_3)=u_3+O(\delta^{\alpha/2}),
\]
we obtain
\begin{multline}
\label{eq:alpha_der_logV123_refined}
\frac{1}{\log\delta}
\frac{\partial}{\partial\alpha}
\log(-V_{123})
=
\frac12 Q_{R_1}(u_2,u_3)
\\
+
\frac{1}{\log\delta}
\Bigg[
\frac12\log d_{1,2}
\left(
A_{R_1,1}u_2-1
\right)
+
\frac12\log d_{1,3}
\left(
A_{R_1,2}u_3-1
\right)
+
\dot R_1(\alpha)B_{R_1}(u_2,u_3)
\Bigg]
\\
+
O\!\left(\frac{\delta^{\alpha/2}}{|\log\delta|}\right).
\end{multline}

We now return to the full density. From
\eqref{eq:log_triple_density_expansion}, the terms coming from \(-V\), from
the non-dominant partition terms in
\eqref{eq:triple_density_partition}, and from the first-order correction in
\(\log(-V_{123})\) contribute at order \(\delta^{\alpha/2}\) after division
by \(\log\delta\). Hence there exists a measurable function
\(\Omega_\sigma\), depending on
\((u_2,u_3;d_{1,2},d_{1,3},d_{2,3})\), such that
\begin{multline}
\label{eq:alpha_score_triple_refined}
\frac{1}{\log\delta}
\frac{\partial}{\partial\alpha}
\log f_{x_1,x_2,x_3}
(z_1,z_2,z_3;\sigma,\alpha)
=
\frac12 Q_{R_1}(u_2,u_3)
\\
+
\frac{1}{\log\delta}
\Bigg[
\frac12\log d_{1,2}
\left(
A_{R_1,1}u_2-1
\right)
+
\frac12\log d_{1,3}
\left(
A_{R_1,2}u_3-1
\right)
+
\dot R_1(\alpha)B_{R_1}(u_2,u_3)
\Bigg]
\\
+
\frac{\sigma\delta^{\alpha/2}}{2z_1}
\Omega_\sigma(u_2,u_3;d_{1,2},d_{1,3},d_{2,3})
+
o(\delta^{\alpha/2})
+
O\!\left(\frac{\delta^{\alpha/2}}{|\log\delta|}\right).
\end{multline}
In particular,
\[
\frac{1}{\log\delta}
\frac{\partial}{\partial\alpha}
\log f_{x_1,x_2,x_3}
(z_1,z_2,z_3;\sigma,\alpha)
=
\frac12 Q_{R_1}(u_2,u_3)
+
O\!\left(\frac{1}{|\log\delta|}\right)
+
O(\delta^{\alpha/2}),
\]
and therefore
\[
\frac{1}{\log\delta}
\frac{\partial}{\partial\alpha}
\log f_{x_1,x_2,x_3}
(z_1,z_2,z_3;\sigma,\alpha)
=
\frac12 Q_{R_1}(u_2,u_3)+o(1).
\]

\subsection{Proof of Theorem \protect\ref{Prop:Asym_Prop_CL_Est}}

We prove first the pairwise results. Throughout this subsection, the true
parameter is denoted by $(\sigma_0,\alpha_0)$, with
$\sigma_0>0$ and $\alpha_0\in(0,1)$. For an edge
$e=(x_1,x_2)\in E_N$, write
\[
d_e=\|x_2-x_1\|,
\qquad
U_e
=
U_{x_1,x_2}^{(\eta)}
=
\frac{1}{\sigma_0 d_e^{\alpha_0/2}}
\log\left(\frac{\eta(x_2)}{\eta(x_1)}\right).
\]
We also set
\[
L_Z=\sum_{j\geq1}\sum_{k>j}L_{Z_{k\setminus j}}(0).
\]

We shall use the following consequence of Proposition
\ref{Prop_pairwise_pdf}. Uniformly for the values of the parameters in a
compact neighbourhood of the true parameter, and for
\[
u=\frac{\log(z_2/z_1)}{\sigma d^{\alpha/2}},
\]
we have, as $d\to0$,
\begin{equation}
\label{eq:pair_score_sigma_uniform}
\frac{\partial}{\partial\sigma}
\log f_{x_1,x_2}(z_1,z_2;\sigma,\alpha)
=
\frac{1}{\sigma}(u^2-1)
+
\frac{\omega(u)}{z_1}d^{\alpha/2}
+
r_\sigma(d,u,z_1),
\end{equation}
where the remainder satisfies, for some integer $k_0\geq1$,
\[
|r_\sigma(d,u,z_1)|
\leq
\frac{C}{z_1}
d^{\alpha}
\sum_{k=0}^{k_0}(1+|u|)^k .
\]
Similarly,
\begin{equation}
\label{eq:pair_score_alpha_uniform}
\frac{1}{\log d}
\frac{\partial}{\partial\alpha}
\log f_{x_1,x_2}(z_1,z_2;\sigma,\alpha)
=
\frac12(u^2-1)
+
\frac{\sigma\omega(u)}{2z_1}d^{\alpha/2}
+
r_\alpha(d,u,z_1),
\end{equation}
with
\[
|r_\alpha(d,u,z_1)|
\leq
\frac{C}{z_1}
d^{\alpha}
\sum_{k=0}^{k_0}(1+|u|)^k .
\]
The bounds above follow by continuing the Taylor expansion in Proposition
\ref{Prop_pairwise_pdf} one order further. Together with Proposition
\ref{prop:uniformbounds} and Lemma \ref{Le_Add_Bounds}, they imply that the
corresponding sums of remainders are negligible under the normalizations used
below.

We also use the following two consequences of Theorem
\ref{prop:BRtrajectories} and of the centering identity for the pairwise
score:
\begin{equation}
\label{eq:pair_V2_limit}
\frac{\sqrt{3}}{3}
N^{-(2-\alpha_0)/4}
V_{2,N}^{(\eta)}
\overset{\mathbb P}{\longrightarrow}
c_{V_2}L_Z,
\end{equation}
and
\begin{equation}
\label{eq:pair_bias_limit}
\frac{\sqrt{3}}{3}
N^{-(2-\alpha_0)/4}
\frac{1}{\sqrt{|E_N|}}
\sum_{e=(x_1,x_2)\in E_N}
\frac{\omega(U_e)}{\eta(x_1)}d_e^{\alpha_0/2}
\overset{\mathbb P}{\longrightarrow}
-\frac{c_{V_2}}{\sigma_0}\mathbb E[L_Z].
\end{equation}
Indeed, \eqref{eq:pair_bias_limit} is the deterministic centering term
associated with the exact identity
\[
\mathbb E\left[
\frac{\partial}{\partial\sigma}
\log f_{x_1,x_2}
\bigl(\eta(x_1),\eta(x_2);\sigma_0,\alpha_0\bigr)
\right]=0,
\]
combined with the local expansion
\eqref{eq:pair_score_sigma_uniform} and the same stabilization argument for
Poisson--Delaunay edges as in the proof of Theorem
\ref{prop:BRtrajectories}.

\subsubsection{Proof for $\widehat{\sigma}_{2,N}$}

Assume first that $\alpha_0$ is known. Let
$\widehat\sigma_{2,N}$ be a sequence of pairwise maximum composite likelihood
estimators of $\sigma_0$. It satisfies the first-order condition
\[
0
=
\frac{\partial}{\partial\sigma}
\ell_{2,N}(\widehat\sigma_{2,N},\alpha_0)
=
\sum_{e=(x_1,x_2)\in E_N}
\left.
\frac{\partial}{\partial\sigma}
\log f_{x_1,x_2}
\bigl(\eta(x_1),\eta(x_2);\sigma,\alpha_0\bigr)
\right|_{\sigma=\widehat\sigma_{2,N}} .
\]
At the parameter value $\widehat\sigma_{2,N}$, the normalized increment is
\[
u_e(\widehat\sigma_{2,N})
=
\frac{1}{\widehat\sigma_{2,N}d_e^{\alpha_0/2}}
\log\left(\frac{\eta(x_2)}{\eta(x_1)}\right)
=
\frac{\sigma_0}{\widehat\sigma_{2,N}}U_e.
\]
Using \eqref{eq:pair_score_sigma_uniform} and multiplying the score equation
by $\widehat\sigma_{2,N}$, we obtain
\begin{equation}
\label{eq:FOC_sigma_pair}
0
=
\sum_{e\in E_N}
\left[
\frac{\sigma_0^2}{\widehat\sigma_{2,N}^{\,2}}U_e^2-1
\right]
+
\widehat\sigma_{2,N}
\sum_{e=(x_1,x_2)\in E_N}
\frac{
\omega\left(\frac{\sigma_0}{\widehat\sigma_{2,N}}U_e\right)
}{\eta(x_1)}
d_e^{\alpha_0/2}
+
\widehat\sigma_{2,N}
\sum_{e=(x_1,x_2)\in E_N}
r_{\sigma,e}(\widehat\sigma_{2,N}),
\end{equation}
where the last term denotes the corresponding remainder.

Rearranging \eqref{eq:FOC_sigma_pair} gives
\begin{multline}
\label{eq:sigma_rearranged_pair}
\left(\widehat\sigma_{2,N}^{\,2}-\sigma_0^2\right)
\left[
\frac{1}{\widehat\sigma_{2,N}^{\,2}}
\frac{1}{|E_N|}\sum_{e\in E_N}U_e^2
\right]
=
\frac{1}{|E_N|}\sum_{e\in E_N}(U_e^2-1)
\\
+
\widehat\sigma_{2,N}
\frac{1}{|E_N|}
\sum_{e=(x_1,x_2)\in E_N}
\frac{
\omega\left(\frac{\sigma_0}{\widehat\sigma_{2,N}}U_e\right)
}{\eta(x_1)}
d_e^{\alpha_0/2}
+
\widehat\sigma_{2,N}
\frac{1}{|E_N|}
\sum_{e=(x_1,x_2)\in E_N}
r_{\sigma,e}(\widehat\sigma_{2,N}).
\end{multline}

The first term on the right-hand side is
\[
\frac{1}{|E_N|}\sum_{e\in E_N}(U_e^2-1)
=
\frac{1}{\sqrt{|E_N|}}V_{2,N}^{(\eta)}
=
O_{\mathbb P}(N^{-\alpha_0/4}),
\]
by \eqref{eq:pair_V2_limit} and $|E_N|=O_{\mathbb P}(N)$. The second term is
also $O_{\mathbb P}(N^{-\alpha_0/4})$ by
\eqref{eq:pair_bias_limit}. The remainder term is
$o_{\mathbb P}(N^{-\alpha_0/4})$ by the bound in
\eqref{eq:pair_score_sigma_uniform}, Proposition
\ref{prop:uniformbounds}, and Lemma \ref{Le_Add_Bounds}. Since
\[
\frac{1}{|E_N|}\sum_{e\in E_N}U_e^2
\overset{\mathbb P}{\longrightarrow}1,
\]
we first obtain
\[
\widehat\sigma_{2,N}^{\,2}
\overset{\mathbb P}{\longrightarrow}
\sigma_0^2.
\]

We can now replace
\[
\omega\left(\frac{\sigma_0}{\widehat\sigma_{2,N}}U_e\right)
\]
by $\omega(U_e)$ in the normalized score, since
$\widehat\sigma_{2,N}\to\sigma_0$ in probability and the function
$\omega$ has at most polynomial growth. Multiplying
\eqref{eq:sigma_rearranged_pair} by
\[
\frac{\sqrt{3}}{3}\sqrt{|E_N|}N^{-(2-\alpha_0)/4},
\]
and using \eqref{eq:pair_V2_limit}, \eqref{eq:pair_bias_limit}, and the
negligibility of the remainder, yields
\[
\frac{\sqrt{3}}{3}
\sqrt{|E_N|}
N^{-(2-\alpha_0)/4}
\left(\widehat\sigma_{2,N}^{\,2}-\sigma_0^2\right)
\overset{\mathbb P}{\longrightarrow}
c_{V_2}\sigma_0^2
\left(L_Z-\mathbb E[L_Z]\right).
\]
This proves the pairwise result for the scale parameter.

\subsubsection{Proof for $\widehat{\alpha}_{2,N}$}

Assume now that $\sigma_0$ is known. Let
$\widehat\alpha_{2,N}$ be a sequence of pairwise maximum composite likelihood
estimators of $\alpha_0$. It satisfies
\[
0
=
\frac{\partial}{\partial\alpha}
\ell_{2,N}(\sigma_0,\widehat\alpha_{2,N})
=
\sum_{e=(x_1,x_2)\in E_N}
\left.
\frac{\partial}{\partial\alpha}
\log f_{x_1,x_2}
\bigl(\eta(x_1),\eta(x_2);\sigma_0,\alpha\bigr)
\right|_{\alpha=\widehat\alpha_{2,N}} .
\]
At the parameter value $\widehat\alpha_{2,N}$, the normalized increment is
\[
u_e(\widehat\alpha_{2,N})
=
\frac{1}{\sigma_0 d_e^{\widehat\alpha_{2,N}/2}}
\log\left(\frac{\eta(x_2)}{\eta(x_1)}\right)
=
\rho_{e,N}U_e,
\]
where
\[
\rho_{e,N}
=
d_e^{-(\widehat\alpha_{2,N}-\alpha_0)/2}.
\]
Using \eqref{eq:pair_score_alpha_uniform}, the score equation becomes
\begin{multline}
\label{eq:FOC_alpha_pair}
0
=
\sum_{e=(x_1,x_2)\in E_N}
\log d_e
\left[
\rho_{e,N}^2U_e^2-1
+
\sigma_0
\frac{\omega(\rho_{e,N}U_e)}{\eta(x_1)}
d_e^{\widehat\alpha_{2,N}/2}
+
\widetilde r_{\alpha,e}(\widehat\alpha_{2,N})
\right],
\end{multline}
where the factor $1/2$ has been removed from the equation and
$\widetilde r_{\alpha,e}$ denotes the corresponding remainder.

We shall use the standard maximal edge bound for Poisson--Delaunay
triangulations,
\[
\max_{e\in E_N}\sqrt N\,d_e=O_{\mathbb P}(\log N),
\]
see, for instance, \cite{Henze82,Chenavier22}. Hence, uniformly over
$e\in E_N$,
\[
\log d_e
=
-\frac12\log N+O_{\mathbb P}(\log\log N).
\]
Dividing \eqref{eq:FOC_alpha_pair} by $-\frac12\log N$, we obtain
\begin{multline}
\label{eq:FOC_alpha_unweighted_pair}
0
=
\frac{1}{|E_N|}
\sum_{e\in E_N}
\left[
\rho_{e,N}^2U_e^2-1
+
\sigma_0
\frac{\omega(\rho_{e,N}U_e)}{\eta(x_1)}
d_e^{\widehat\alpha_{2,N}/2}
+
\widetilde r_{\alpha,e}(\widehat\alpha_{2,N})
\right]
\\
+
O_{\mathbb P}\left(\frac{\log\log N}{\log N}\right)
\left[
\frac{1}{|E_N|}
\sum_{e\in E_N}
\left|
\rho_{e,N}^2U_e^2-1
+
\sigma_0
\frac{\omega(\rho_{e,N}U_e)}{\eta(x_1)}
d_e^{\widehat\alpha_{2,N}/2}
+
\widetilde r_{\alpha,e}(\widehat\alpha_{2,N})
\right|
\right].
\end{multline}
The term in square brackets is tight, by Proposition
\ref{prop:uniformbounds}, the compactness of the parameter space, and the
remainder bounds. Therefore the last line in
\eqref{eq:FOC_alpha_unweighted_pair} is
$o_{\mathbb P}(1)$.

Since
\[
\frac{1}{|E_N|}\sum_{e\in E_N}U_e^2\overset{\mathbb P}{\longrightarrow}1,
\qquad
\frac{1}{|E_N|}\sum_{e\in E_N}(U_e^2-1)=O_{\mathbb P}(N^{-\alpha_0/4}),
\]
and the correction and remainder terms are $o_{\mathbb P}(1)$, the
first-order condition implies
\[
|\widehat\alpha_{2,N}-\alpha_0|\log N
\overset{\mathbb P}{\longrightarrow}0.
\]
Indeed,
\[
\rho_{e,N}^2
=
\exp\left\{-(\widehat\alpha_{2,N}-\alpha_0)\log d_e\right\},
\]
and the leading contribution of the average
$|E_N|^{-1}\sum_e(\rho_{e,N}^2U_e^2-1)$ is controlled by
\[
\exp\left\{\frac12(\widehat\alpha_{2,N}-\alpha_0)\log N
(1+o_{\mathbb P}(1))\right\}-1.
\]

Consequently,
\[
\rho_{e,N}^2
=
1-\left(\widehat\alpha_{2,N}-\alpha_0\right)\log d_e
+
o_{\mathbb P}
\left(|\widehat\alpha_{2,N}-\alpha_0||\log d_e|\right),
\]
uniformly over $e\in E_N$. Since
\[
-\log d_e
=
\frac12\log N+O_{\mathbb P}(\log\log N),
\]
we obtain
\[
\frac{1}{|E_N|}
\sum_{e\in E_N}
(\rho_{e,N}^2U_e^2-1)
=
\frac12(\widehat\alpha_{2,N}-\alpha_0)\log N
\frac{1}{|E_N|}\sum_{e\in E_N}U_e^2
+
\frac{1}{|E_N|}\sum_{e\in E_N}(U_e^2-1)
+
o_{\mathbb P}
\left(|\widehat\alpha_{2,N}-\alpha_0|\log N\right).
\]
Furthermore, since
$|\widehat\alpha_{2,N}-\alpha_0|\log N\to0$ in probability,
\[
\frac{1}{|E_N|}
\sum_{e=(x_1,x_2)\in E_N}
\frac{\omega(\rho_{e,N}U_e)}{\eta(x_1)}
d_e^{\widehat\alpha_{2,N}/2}
=
\frac{1}{|E_N|}
\sum_{e=(x_1,x_2)\in E_N}
\frac{\omega(U_e)}{\eta(x_1)}
d_e^{\alpha_0/2}
+
o_{\mathbb P}(N^{-\alpha_0/4}).
\]
The remainder term is also $o_{\mathbb P}(N^{-\alpha_0/4})$. Hence
\eqref{eq:FOC_alpha_unweighted_pair} yields
\begin{multline}
\label{eq:alpha_linearized_pair}
\frac12(\widehat\alpha_{2,N}-\alpha_0)\log N
\frac{1}{|E_N|}\sum_{e\in E_N}U_e^2
\\
=
-
\frac{1}{|E_N|}\sum_{e\in E_N}(U_e^2-1)
-
\sigma_0
\frac{1}{|E_N|}
\sum_{e=(x_1,x_2)\in E_N}
\frac{\omega(U_e)}{\eta(x_1)}
d_e^{\alpha_0/2}
+
o_{\mathbb P}(N^{-\alpha_0/4}).
\end{multline}

Multiplying \eqref{eq:alpha_linearized_pair} by
\[
\frac{\sqrt{3}}{3}\sqrt{|E_N|}N^{-(2-\alpha_0)/4},
\]
and using \eqref{eq:pair_V2_limit}, \eqref{eq:pair_bias_limit}, and
$|E_N|^{-1}\sum_eU_e^2\to1$ in probability, we get
\[
\frac{\sqrt{3}}{6}
\sqrt{|E_N|}
N^{-(2-\alpha_0)/4}
\log N
\left(\widehat\alpha_{2,N}-\alpha_0\right)
\overset{\mathbb P}{\longrightarrow}
-
c_{V_2}
\left(L_Z-\mathbb E[L_Z]\right).
\]
This proves the pairwise result for the smoothness parameter.

\subsubsection{Proof for $\widehat{\sigma}_{3,N}$ and $\widehat{\alpha}_{3,N}$}

We use the same strategy as in the proofs for
$\widehat{\sigma}_{2,N}$ and $\widehat{\alpha}_{2,N}$, replacing the
pairwise score expansion by the triplewise expansion of Proposition
\ref{Prop_triplewise_pdf}, and replacing the statistic
$V_{2,N}^{(\eta)}$ by $V_{3,N}^{(\eta)}$.

For a triangle $ \Delta=(x_1,x_2,x_3)\in DT_N$, write
\[
Q_{x_1,x_2,x_3}
=
Q_{R_{x_1,x_2,x_3}}
\left(
U_{x_1,x_2}^{(\eta)},
U_{x_1,x_3}^{(\eta)}
\right).
\]
Then
\[
V_{3,N}^{(\eta)}
=
\frac{1}{\sqrt{|DT_N|}}
\sum_{\Delta \in DT_N} Q_{x_1,x_2,x_3} .
\]
By Theorem \ref{prop:BRtrajectories},
\[
\frac{\sqrt{2}}{2}
N^{-(2-\alpha_0)/4}
V_{3,N}^{(\eta)}
\overset{\mathbb P}{\longrightarrow}
c_{V_3}L_Z,
\qquad
L_Z=\sum_{j\geq1}\sum_{k>j}L_{Z_{k\setminus j}}(0).
\]

When $\alpha_0$ is known, the first-order condition for
$\widehat{\sigma}_{3,N}$ is
\[
\frac{\partial}{\partial\sigma}
\ell_{3,N}(\widehat{\sigma}_{3,N},\alpha_0)=0.
\]
Using the expansion of the triplewise score with respect to $\sigma$ in
Proposition \ref{Prop_triplewise_pdf}, evaluated at
$\sigma=\widehat{\sigma}_{3,N}$, we obtain an identity of the same form as in
the pairwise case:
\[
\widehat{\sigma}_{3,N}^{\,2}-\sigma_0^2
=
\frac{\sigma_0^2}{|DT_N|}
\sum_{\Delta \in DT_N} Q_{x_1,x_2,x_3}
+
\text{centering term}
+
o_{\mathbb P}\left(N^{-\alpha_0/4}\right),
\]
where the centering term is generated by the first-order correction
involving $\Omega$. The exact score identity implies that this correction
centers the limit, namely
\[
\frac{\sqrt{2}}{2}
N^{-(2-\alpha_0)/4}
\frac{1}{\sqrt{|DT_N|}}
\sum_{\Delta \in DT_N}
\text{correction}_{\sigma,T}
\overset{\mathbb P}{\longrightarrow}
-c_{V_3}\mathbb E[L_Z].
\]
The remainder terms are negligible by the same stabilization and moment
estimates as in the proof for $\widehat{\sigma}_{2,N}$, together with the
uniform moment bounds of Proposition \ref{prop:uniformbounds}. Hence
\[
\frac{\sqrt{2}}{2}
\sqrt{|DT_N|}
N^{-(2-\alpha_0)/4}
\left(
\widehat{\sigma}_{3,N}^{\,2}-\sigma_0^2
\right)
\overset{\mathbb P}{\longrightarrow}
c_{V_3}\sigma_0^2
\left(
L_Z-\mathbb E[L_Z]
\right).
\]

We now consider $\widehat{\alpha}_{3,N}$, assuming that $\sigma_0$ is known.
The first-order condition is
\[
\frac{\partial}{\partial\alpha}
\ell_{3,N}(\sigma_0,\widehat{\alpha}_{3,N})=0.
\]
The proof is again parallel to the pairwise case. For each triangle
$T=(x_1,x_2,x_3)$, the normalized increments at the candidate value
$\widehat{\alpha}_{3,N}$ are obtained by multiplying
$U_{x_1,x_2}^{(\eta)}$ and $U_{x_1,x_3}^{(\eta)}$ by the factors
\[
\rho_{1,2,T}
=
d_{1,2}^{-(\widehat{\alpha}_{3,N}-\alpha_0)/2},
\qquad
\rho_{1,3,T}
=
d_{1,3}^{-(\widehat{\alpha}_{3,N}-\alpha_0)/2}.
\]
Since Delaunay edge lengths are of order $N^{-1/2}$, uniformly up to
logarithmic factors,
\[
\log d_{1,2}
=
-\frac12\log N+O_{\mathbb P}(\log\log N),
\qquad
\log d_{1,3}
=
-\frac12\log N+O_{\mathbb P}(\log\log N).
\]
As in the pairwise proof, the first-order condition first implies
\[
|\widehat{\alpha}_{3,N}-\alpha_0|\log N
\overset{\mathbb P}{\longrightarrow}0.
\]
Thus the triplewise score can be linearized around $\alpha_0$.

Using Proposition \ref{Prop_triplewise_pdf}, the leading term in the
linearized score is proportional to
\[
\frac{1}{|DT_N|}
\sum_{\Delta \in DT_N}Q_{x_1,x_2,x_3}.
\]
The first-order correction involving $\Omega$ provides the deterministic
centering term, while the additional term
\[
\frac{1}{\log\delta}\,
B_\alpha(u_2,u_3;d_{1,2},d_{1,3},d_{2,3})
\]
appearing in Proposition \ref{Prop_triplewise_pdf} is negligible after
summation and normalization. This follows from the same stabilization
arguments for Poisson--Delaunay triangles and from the fact that this term is
centered under the limiting Gaussian law associated with the triangle shape.
The remainder terms are also negligible under the normalization of the
theorem.

Consequently,
\[
\frac{\sqrt{2}}{4}
\sqrt{|DT_N|}
N^{-(2-\alpha_0)/4}
\log N
\left(
\widehat{\alpha}_{3,N}-\alpha_0
\right)
\overset{\mathbb P}{\longrightarrow}
-c_{V_3}
\left(
L_Z-\mathbb E[L_Z]
\right).
\]
This proves the triplewise part of Theorem
\ref{Prop:Asym_Prop_CL_Est}.

\section{Technical lemmas}

\label{sec:intermediary_results}

\subsection{Some additional bounds\label{Sec_Bounds_gamma}}

\begin{lemma}
\label{Le_Add_Bounds}
Let $\alpha\in(0,2)$. Then the following assertions hold.
\begin{enumerate}[(i)]
\item There exists a constant $C<\infty$, independent of $N$, such that, for
all sufficiently large $N$,
\[
\mathbb{E}\left[
\sum_{(x_1,x_2)\in E_N}
\|x_2-x_1\|^\alpha
\right]
\leq
C N^{1-\alpha/2}.
\]

\item As $N\to\infty$,
\[
\frac{1}{|E_N|}
\sum_{(x_1,x_2)\in E_N}
\frac{(\gamma(x_2)-\gamma(x_1))^2}{\|x_2-x_1\|^\alpha}
=
O_{\mathbb{P}}\left(N^{-1+\alpha/2}\right).
\]
\end{enumerate}
\end{lemma}

\begin{prooft}{Lemma \ref{Le_Add_Bounds}}
Throughout the proof, $C$ denotes a positive constant whose value may change
from line to line. We write
\[
\mathbf C_N:=N^{1/2}\mathbf C.
\]
Let $p_N(x_1,x_2)$ denote the probability that $x_1$ and $x_2$ are Delaunay
neighbours in the triangulation generated by $P_N\cup\{x_1,x_2\}$. By the
scaling property of the Poisson point process,
\[
p_N(x_1,x_2)=p_1(N^{1/2}x_1,N^{1/2}x_2).
\]
where, for all $y_1,y_2\in\mathbb R^2$,
\[
p_1(y_1,y_2)
\leq
\left(\pi \|y_2-y_1\|^2+4\right)
\exp\left\{-\frac{\pi}{4}\|y_2-y_1\|^2\right\},
\]
see e.g. Lemma 5 in \cite{Chenavier&Robert25a}.

\medskip
\noindent\textit{Proof of (i).}
By the Slivnyak--Mecke formula, and by dropping the lexicographic constraint,
which can only increase the integral, we have
\begin{align*}
\mathbb{E}\left[
\sum_{(x_1,x_2)\in E_N}
\|x_2-x_1\|^\alpha
\right]
&\leq
N^2
\int_{\mathbf C\times\mathbb R^2}
\|x_2-x_1\|^\alpha
p_N(x_1,x_2)
\,\mathrm{d}x_1\,\mathrm{d}x_2 .
\end{align*}
Using the change of variables $y_i=N^{1/2}x_i$, $i=1,2$, gives
\begin{align*}
\mathbb{E}\left[
\sum_{(x_1,x_2)\in E_N}
\|x_2-x_1\|^\alpha
\right]
&\leq
N^{-\alpha/2}
\int_{\mathbf C_N\times\mathbb R^2}
\|y_2-y_1\|^\alpha
p_1(y_1,y_2)
\,\mathrm{d}y_1\,\mathrm{d}y_2 .
\end{align*}
By stationarity of the Poisson--Delaunay triangulation,
$p_1(y_1,y_2)=p_1(0,y_2-y_1)$. Therefore
\begin{align*}
\mathbb{E}\left[
\sum_{(x_1,x_2)\in E_N}
\|x_2-x_1\|^\alpha
\right]
&\leq
N^{-\alpha/2}|\mathbf C_N|
\int_{\mathbb R^2}
\|h\|^\alpha p_1(0,h)\,\mathrm{d}h .
\end{align*}
Since $|\mathbf C_N|=N|\mathbf C|=N$ and
\[
p_1(0,h)
\leq
C\|h\|^2 e^{-\pi\|h\|^2/4},
\]
the last integral is finite. Hence
\[
\mathbb{E}\left[
\sum_{(x_1,x_2)\in E_N}
\|x_2-x_1\|^\alpha
\right]
\leq
C N^{1-\alpha/2}.
\]

\medskip
\noindent\textit{Proof of (ii).}
Set
\[
\Gamma_N
=
\frac{1}{N}
\sum_{(x_1,x_2)\in E_N}
\frac{(\gamma(x_2)-\gamma(x_1))^2}{\|x_2-x_1\|^\alpha}.
\]
Since
\[
\gamma(x)=\frac{\sigma^2}{2}\|x\|^\alpha,
\]
the same Slivnyak--Mecke formula and the change of variables
$y_i=N^{1/2}x_i$ yield
\begin{align*}
\mathbb{E}[\Gamma_N]
&\leq
C N^{-1-\alpha/2}
\int_{\mathbf C_N\times\mathbb R^2}
\frac{
\bigl(\|y_2\|^\alpha-\|y_1\|^\alpha\bigr)^2
}{
\|y_2-y_1\|^\alpha
}
p_1(y_1,y_2)
\,\mathrm{d}y_1\,\mathrm{d}y_2 .
\end{align*}
Writing $y_2=y_1+h$, this becomes
\[
\mathbb{E}[\Gamma_N]
\leq
C N^{-1-\alpha/2}
\int_{\mathbf C_N} I(y)\,\mathrm{d}y,
\]
where
\[
I(y)
=
\int_{\mathbb R^2}
\frac{
\bigl(\|y+h\|^\alpha-\|y\|^\alpha\bigr)^2
}{
\|h\|^\alpha
}
p_1(y,y+h)
\,\mathrm{d}h .
\]
We now bound $I(y)$. The estimate on $p_1$ gives
\[
p_1(y,y+h)
\leq
C\|h\|^2 e^{-\pi\|h\|^2/4}.
\]

First, if $\|y\|\leq 1$, then
\[
\bigl|\|y+h\|^\alpha-\|y\|^\alpha\bigr|
\leq
C(1+\|h\|^\alpha),
\]
and hence
\[
I(y)
\leq
C
\int_{\mathbb R^2}
\left(
\|h\|^{-\alpha}+\|h\|^\alpha
\right)
\|h\|^2 e^{-\pi\|h\|^2/4}
\,\mathrm{d}h
\leq C.
\]

Assume now that $\|y\|>1$. We split the integral defining $I(y)$ into the two
regions
\[
A_1(y)=\{\|h\|\leq \|y\|/2\},
\qquad
A_2(y)=\{\|h\|>\|y\|/2\}.
\]
On $A_1(y)$, the mean value theorem applied to
$x\mapsto \|x\|^\alpha$ gives
\[
\bigl|\|y+h\|^\alpha-\|y\|^\alpha\bigr|
\leq
C\|y\|^{\alpha-1}\|h\|.
\]
Therefore
\begin{align*}
\int_{A_1(y)}
\frac{
\bigl(\|y+h\|^\alpha-\|y\|^\alpha\bigr)^2
}{
\|h\|^\alpha
}
p_1(y,y+h)
\,\mathrm{d}h
&\leq
C\|y\|^{2\alpha-2}
\int_{\mathbb R^2}
\|h\|^{2-\alpha}
\|h\|^2 e^{-\pi\|h\|^2/4}
\,\mathrm{d}h  \\
&\leq
C\|y\|^{2\alpha-2}.
\end{align*}
On $A_2(y)$, we use
\[
\bigl|\|y+h\|^\alpha-\|y\|^\alpha\bigr|
\leq
C(\|h\|^\alpha+\|y\|^\alpha).
\]
Since $\|h\|>\|y\|/2$ on $A_2(y)$, it follows that
\[
\frac{
\bigl(\|y+h\|^\alpha-\|y\|^\alpha\bigr)^2
}{
\|h\|^\alpha
}
\leq
C\|h\|^\alpha.
\]
Thus
\[
\int_{A_2(y)}
\frac{
\bigl(\|y+h\|^\alpha-\|y\|^\alpha\bigr)^2
}{
\|h\|^\alpha
}
p_1(y,y+h)
\,\mathrm{d}h
\leq
C
\int_{\|h\|>\|y\|/2}
\|h\|^{\alpha+2}e^{-\pi\|h\|^2/4}
\,\mathrm{d}h.
\]
The last term is exponentially small as $\|y\|\to\infty$ and, in particular,
is bounded by $C\|y\|^{2\alpha-2}$ for $\|y\|>1$, after increasing $C$ if
necessary. Consequently,
\[
I(y)
\leq
C\mathbf 1_{\{\|y\|\leq 1\}}
+
C\|y\|^{2\alpha-2}\mathbf 1_{\{\|y\|>1\}}.
\]

We deduce that
\begin{align*}
\int_{\mathbf C_N} I(y)\,\mathrm{d}y
&\leq
C
+
C
\int_{\mathbf C_N\cap\{\|y\|>1\}}
\|y\|^{2\alpha-2}\,\mathrm{d}y  \\
&\leq
C
+
C
\int_1^{c\sqrt N}
r^{2\alpha-2}r\,\mathrm{d}r \\
&=
O(N^\alpha).
\end{align*}
Hence
\[
\mathbb{E}[\Gamma_N]
=
O\left(N^{-1-\alpha/2}N^\alpha\right)
=
O\left(N^{-1+\alpha/2}\right).
\]
Since $\Gamma_N$ is non-negative, Markov's inequality gives
\[
\Gamma_N
=
O_{\mathbb P}\left(N^{-1+\alpha/2}\right).
\]

Finally, the ergodic theorem for Poisson--Delaunay tessellations yields
\[
\frac{|E_N|}{N}
\overset{\mathbb P}{\longrightarrow}
\beta_1=3,
\]
up to negligible boundary effects (see e.g. Theorem 10.2.9 in
\cite{Schneider&Weil08}). In particular, $N/|E_N|=O_{\mathbb P}(1)$.
Therefore
\[
\frac{1}{|E_N|}
\sum_{(x_1,x_2)\in E_N}
\frac{(\gamma(x_2)-\gamma(x_1))^2}{\|x_2-x_1\|^\alpha}
=
\frac{N}{|E_N|}\Gamma_N
=
O_{\mathbb P}\left(N^{-1+\alpha/2}\right),
\]
which proves the second assertion.
\end{prooft}

\subsection{Uniform bounds for the moments of $|U_{x_{1},x_{2}}^{(\protect\eta )}|^{p}$}

\begin{proposition}
\label{prop:uniformbounds}
Let $\alpha\in(0,2)$ and $\sigma>0$ be fixed. For any $p\geq1$, there exists
a constant $C_p<\infty$ such that
\[
\sup_{\substack{x_1,x_2\in\mathbb{R}^{2}\\ x_1\neq x_2}}
\mathbb{E}
\left[
\left|U_{x_1,x_2}^{(\eta)}\right|^p
\right]
\leq
C_p .
\]
\end{proposition}

\begin{prooft}{Proposition \ref{prop:uniformbounds}}
Let
\[
d=\|x_2-x_1\|>0,
\qquad
a=\sigma d^{\alpha/2}.
\]
By Proposition \ref{Prop_cond_dist_U}, the marginal distribution function of
$U_{x_1,x_2}^{(\eta)}$ is
\[
F_a(u)
:=
\mathbb{P}\left\{
U_{x_1,x_2}^{(\eta)}\leq u
\right\}
=
\frac{\Phi(u+a/2)}
{
\Phi(u+a/2)+e^{-au}\Phi(-u+a/2)
},
\qquad u\in\mathbb{R}.
\]
In particular, the distribution of $U_{x_1,x_2}^{(\eta)}$ depends on
$x_1,x_2$ only through $a=\sigma\|x_2-x_1\|^{\alpha/2}$.

We first consider the case of small distances. Choose $\delta>0$ such that
\[
\sigma\delta^{\alpha/2}\leq1.
\]
Assume that $0<d\leq\delta$, so that $0<a\leq1$. For $t\geq1$, using the
explicit expression of $F_a$, we obtain
\[
\mathbb{P}\left\{
U_{x_1,x_2}^{(\eta)}>t
\right\}
=
1-F_a(t)
=
\frac{e^{-at}\Phi(-t+a/2)}
{
\Phi(t+a/2)+e^{-at}\Phi(-t+a/2)
}
\leq
\frac{\Phi(-t+a/2)}{\Phi(t+a/2)}.
\]
Since $a\leq1$ and $t\geq1$, we have $-t+a/2\leq -t/2$ and
$\Phi(t+a/2)\geq \Phi(1)\geq 1/2$. Hence
\[
\mathbb{P}\left\{
U_{x_1,x_2}^{(\eta)}>t
\right\}
\leq
2\Phi(-t/2).
\]
Similarly,
\[
\mathbb{P}\left\{
U_{x_1,x_2}^{(\eta)}<-t
\right\}
=
F_a(-t)
=
\frac{\Phi(-t+a/2)}
{
\Phi(-t+a/2)+e^{at}\Phi(t+a/2)
}
\leq
e^{-at}
\frac{\Phi(-t+a/2)}{\Phi(t+a/2)}
\leq
2\Phi(-t/2).
\]
Therefore, for all $t\geq1$ and all $0<d\leq\delta$,
\[
\mathbb{P}
\left\{
\left|U_{x_1,x_2}^{(\eta)}\right|>t
\right\}
\leq
4\Phi(-t/2).
\]
Using the identity
\[
\mathbb{E}\{|X|^p\}
=
p\int_0^\infty t^{p-1}\mathbb{P}\{|X|>t\}\,\mathrm{d}t,
\]
valid for any real random variable $X$, we deduce that
\[
\sup_{0<\|x_2-x_1\|\leq\delta}
\mathbb{E}
\left[
\left|U_{x_1,x_2}^{(\eta)}\right|^p
\right]
\leq
p\int_0^1 t^{p-1}\,\mathrm{d}t
+
4p\int_1^\infty t^{p-1}\Phi(-t/2)\,\mathrm{d}t
<\infty.
\]

It remains to consider the case $\|x_2-x_1\|>\delta$. In this case,
\[
\left|U_{x_1,x_2}^{(\eta)}\right|^p
=
\frac{
\left|\log\eta(x_2)-\log\eta(x_1)\right|^p
}{
\sigma^p\|x_2-x_1\|^{p\alpha/2}
}
\leq
\frac{
2^{p-1}
\left(
|\log\eta(x_1)|^p+|\log\eta(x_2)|^p
\right)
}{
\sigma^p\delta^{p\alpha/2}
}.
\]
Since $\eta(x)$ has a standard unit Fr\'echet distribution for every
$x\in\mathbb{R}^{2}$, the random variable $\log\eta(x)$ has a standard
Gumbel distribution. In particular,
\[
\mathbb{E}\left[|\log\eta(x)|^p\right]<\infty,
\qquad x\in\mathbb{R}^{2}.
\]
Thus
\[
\sup_{\|x_2-x_1\|>\delta}
\mathbb{E}
\left[
\left|U_{x_1,x_2}^{(\eta)}\right|^p
\right]
\leq
C_p'
\delta^{-p\alpha/2}
<\infty.
\]
Combining the bounds for $0<\|x_2-x_1\|\leq\delta$ and for
$\|x_2-x_1\|>\delta$ proves the result.
\end{prooft}

\end{document}